\documentclass[10pt,a4paper]{amsart}
\usepackage[utf8]{inputenc}
\usepackage[pagebackref,colorlinks,linkcolor=red,citecolor=blue,urlcolor=blue,hypertexnames=false]{hyperref}
\usepackage{enumerate}
\usepackage{csquotes}
\usepackage{graphicx}
\usepackage{amsmath,amsthm,amssymb,amsopn}
\usepackage{amsrefs}
\makeatletter
\newcommand{\makesynonym}[2]{%
  \AtEndDocument{\write\@auxout{\noexpand\@makesynonym{#1}{#2}}}%
}
\newcommand{\@makesynonym}[2]{%
  \global\expandafter\let\csname b@#1\expandafter\endcsname\csname b@#2\endcsname
}
\makeatother

\makesynonym{weibel1994introduction}{chuhom}
\usepackage{url}
\usepackage{enumerate}
\usepackage{tikz-cd}
\makeatletter
\usepackage{amscd}
\input xy
\xyoption{all}
\newdir{ >}{{}*!/-5pt/@{>}}
\usepackage{mathtools}
\let\overfence\overbrace % \overfence is similar to \overbrace
\let\downfencefill\downbracefill % match components of \overbrace
\patchcmd{\overfence}{\downbracefill}{\downfencefill}{}{}% patch \overfence...
\patchcmd{\downfencefill}{\braceru \bracelu}{}{}{}%... and \downfencefill

\usepackage{etoolbox}
\apptocmd{\sloppy}{\hbadness 10000\relax}{}{}

\newcommand{\comment}[1]{} 
\newcommand{\xto}{\xrightarrow}

\newcommand{\scst}{\scriptstyle}

\def\N{\mathbb{N}} 
\def\K{\mathbb{K}}
\def\Z{\mathbb{Z}} 
\def\Q{\mathbb{Q}} 
 
\def\C{\mathbb{C}} 
 
\def\F{\mathbb{F}} 
\def\S{\mathbb{S}}
\def\H{\mathbb{H}}
\def\B{\mathbb{B}}
\def\HH{\mathbb{HH}}
\def\HC{\mathbb{HC}}
\def\HP{\mathbb{HP}}
\def\HN{\mathbb{HN}}
\def\Tot{\operatorname{Tot}}

\def\res{\operatorname{res}}

\def\cyc{\operatorname{cyc}}
\def\adj{\operatorname{adj}}
\def\Cy{C^{\cyc}}
\def\cone{\operatorname{cone}}
\def\norm{\operatorname{nor}}
\def\Ho{\operatorname{Ho}}
\def\topdf{\texorpdfstring}
\def\alg{\mathrm{Alg}}

\def\ab{\operatorname{ab}}

\newcommand{\defeq}{\mathrel{:=}} % per Definition

\numberwithin{equation}{subsection}

\theoremstyle{plain}
\newtheorem{thm}[equation]{Theorem}
\newtheorem{lem}[equation]{Lemma}
\newtheorem{coro}[equation]{Corollary} 
\newtheorem{prop}[equation]{Proposition}

\theoremstyle{definition}
\newtheorem{conj}{Conjecture}

\newtheorem{defi}[equation]{Definition} 

\newtheorem{ex}[equation]{Example}

\theoremstyle{remark} 
\newtheorem{rmk}[equation]{Remark}
 \newtheorem{rem}[equation]{Remark}
 \newtheorem{nota}[equation]{Notation}

\newtheorem*{ack}{Acknowledgements}

\newcommand{\cA}{\mathcal A}

\newcommand{\cC}{\mathcal C}
\newcommand{\cD}{\mathcal D}
\newcommand{\cE}{\mathcal E}
\newcommand{\cF}{\mathcal F}
\newcommand{\cG}{\mathcal G}
\newcommand{\cH}{\mathcal H}
\newcommand{\cI}{\mathcal I}

\newcommand{\cK}{\mathcal K}

\newcommand{\cM}{\mathcal M}
\newcommand{\cN}{\mathcal N}

\newcommand{\cP}{\mathcal P}

\newcommand{\cR}{\mathcal R}
\newcommand{\cS}{\mathcal S}
\newcommand{\cT}{\mathcal T}
\newcommand{\cU}{\mathcal U}

\def\fB{\mathfrak{B}}

\def\fF{\mathfrak{F}}

\def\fX{\mathfrak{X}}

\newcommand{\BF}{\fB\fF}

\def\Mod{\operatorname{Mod}}

\def\ttimes{\widetilde{\times}}

\newcommand{\aha}{{\alg_{k}}}
\newcommand{\ahal}{{\alg_{\ell}}}

\def\ul{\underline}
\newcommand{\lra}{\longrightarrow}
\newcommand{\iso}{\overset{\cong}{\lra}}

\newcommand{\weq}{\overset{\sim}{\lra}}

\newcommand{\sg}{\mathrm{sg}}

\newcommand{\onto}{\twoheadrightarrow}

\newcommand{\ol}{\overline}
\newcommand{\olD}{\ol{D}}

\newcommand{\bomega}{\ol{\omega}}

\def\reg{\operatorname{reg}}
\def\sing{\operatorname{sing}}
\def\sink{\operatorname{sink}}
\def\inf{\operatorname{inf}}
\def\sour{\operatorname{sour}}

\def\triqui{\vartriangleleft}
\def\truqui{\vartriangleright}
\def\supp{\operatorname{Supp}}
\def\inc{\operatorname{inc}}
\def\tr{\operatorname{tr}}

\def\ad{\operatorname{ad}}

\def\id{\operatorname{id}}

\def\Iso{\operatorname{Iso}}

\newcommand{\coker}{{\rm Coker}}
\renewcommand{\ker}{{\rm Ker}}

\DeclarePairedDelimiterX{\setgiven}[2]{\{}{\}}{#1\,{:}\,\mathopen{}#2}% 

\newcommand{\im}{\mathrm{Im}}
\newcommand{\Dom}{\mathrm{Dom}}

\newcommand{\op}{\mathrm{op}}

\DeclareMathOperator*{\colim}{colim}

\def\tor{\operatorname{Tor}}
\def\Hom{\operatorname{Hom}}
\def\End{\operatorname{End}}

\def\otimesl{\otimes_\ell}

\newcommand{\bu}{\bullet}

\def\hoco{\operatorname{hocolim}}

\def\supp{\operatorname{Supp}}
\newcommand{\acts}{\bullet}
\newcommand{\cc}{{\mathcal{C}_c}}
\newcommand{\gring}{k}
\newcommand{\mspan}{\operatorname{span}_{\gring}}
\newcommand{\mspanl}{\operatorname{span}_\ell}
\newcommand{\WB}{\mathsf{WeakBool}}
\newcommand{\A}{{\mathcal{A}_{\gring}}}
\def\cyc{\operatorname{cyc}}
\def\tr{\operatorname{tr}}
\def\sing{\operatorname{sing}}
\def\reg{\operatorname{reg}}
\def\sink{\operatorname{sink}}
\def\sour{\operatorname{sour}}
\def\inf{\operatorname{inf}}
\def\Stab{\operatorname{Stab}}
\def\Or{\operatorname{Or}}
\newcommand{\bB}{\mathbb{B}}

\makeatletter
\@namedef{subjclassname@2020}{\textup{2020} Mathematics Subject Classification}
\makeatother

\begin{document}

\hfuzz=22pt
\vfuzz=22pt
\hbadness=2000
\vbadness=\maxdimen

\author{Guido Arnone}

\email{garnone@dm.uba.ar}
\address{Departamento de Matem\'atica/IMAS\\ Facultad de Ciencias Exactas y Naturales\\ Universidad de Buenos Aires\\ Ciudad Universitaria \\ (1428) Buenos Aires}

\author{Guillermo Cortiñas}

\email{gcorti@dm.uba.ar}
\address{Departamento de Matem\'atica/IMAS\\ Facultad de Ciencias Exactas y Naturales\\ Universidad de Buenos Aires\\ Ciudad Universitaria \\ (1428) Buenos Aires}
\author{Devarshi Mukherjee}

\email{dmukherjee@dm.uba.ar}

\address{University of M\"unster\\ Mathematics M\"unster\\ Einsteinstrasse 62\\ 48149 M\"unster}
\thanks{Arnone and Corti\~nas were supported by CONICET and partially supported by grants PIP 423 and UBACyT 206. Mukherjee was funded by a DFG Eigenestelle (project number 534946574) and the Cluster of Excellence: Groups, Geometry and Dynamics, Mathematics M\"unster. }
\title{Homology of Steinberg algebras}
\begin{abstract} We study homological invariants of the Steinberg algebra $\A(\cG)$ of an ample groupoid $\cG$ over a commutative ring $k$. For any ample Hausdorff groupoid $\cG$, we find that $H_*(\cG)$ is a direct summand of $HH_*(\A(\cG))$; using this and the Dennis trace we obtain a map $\bar{D}_*:K_*(\A(\cG))\to H_*(\cG)$. We study this map when $\cG$ is the (twisted) Exel-Pardo groupoid associated to a self-similar action of a group $G$ on a graph, and compute $HH_*(\A(\cG))$ and $H_*(\cG)$ in terms of the homology of $G$, and the $K$-theory of $\A(\cG)$ in terms of that of $k[G]$. 
\end{abstract}
\maketitle

\vspace{-2.5em}
\section{Introduction}\label{sec:intro}
\numberwithin{equation}{section}
A topological groupoid is a groupoid where the sets $\cG$ of arrows and $\cG^{(0)}$ of units are topological spaces and all structure maps are continuous; $\cG$ is \emph{étale} if the range and source maps $r$ and $s$ are local homeomorphisms, and \emph{ample} if in addition $\cG^{(0)}$ is 
Hausdorff and has a basis of compact open subsets. For a commutative ring $k$, we study \(K\)-theoretic and homological invariants of the \emph{Steinberg
algebra} $\A(\cG)$ \cites{steinappr}. This is the \(k\)-module $\cc(\cG)$ spanned by characteristic functions of compact open subsets, equipped with the convolution product. 
Write $\Cy(R/A)$ for the standard semicyclic module (see Example \ref{ex:ccycR}) of a ring $R$ which is an algebra over a $k$-algebra $A$, and $\Cy(R)=\Cy(R/k)$. For a $\cG$-space $W$, we write $\H(\cG,W)$ for the groupoid homology complex with coefficients in $\cc(W)$; $\H(\cG)=\H(\cG,\cG^0)$ is a semicyclic module, and we write $\HC(\cG)$, $\HP(\cG)$ and $\HN(\cG)$ for its cyclic complexes. Consider also the cyclic module $\H^{\cyc}(\cG)$ that results from applying the functor $\cc$ to the cyclic nerve of $\cG$. For $x\in\cG^{(0)}$, let $\cG_x^x=\{g\in\cG\colon r(g)=s(g)=x\}$. Put $\cG^{\Iso}=\bigcup_{x\in\cG^{(0)}}\cG_x^x$; this is a $\cG$-space containing $\cG^{(0)}$ where $\cG$ acts by conjugation. We say that $\cG$ is \emph{principal} if $\cG^{\Iso}=\cG^{(0)}$. The following is our first main theorem.
\begin{thm}\label{thm:burghintro}
Let $\cG$ be an ample groupoid. 
\begin{enumerate}[i)]
\item The natural map  $\Cy(\cA(\cG))\onto \Cy(\cA(\cG)/\cA(\cG^{(0)}))$ is a quasi-isomorphism. 
\item $H_*^{\cyc}(\cG)=H_*(\cG,\cC_c(\cG^{\Iso}))$ and  there are embeddings of semicyclic modules  $\Cy(\cA(\cG)/\cA(\cG^{(0)}))\overset{\iota}{\leftarrow}\H(\cG)\overset{\iota'}{\lra}\H^{\cyc}(\cG)$. If $\cG$ is principal, $\iota'$ is an isomorphism; if it is Hausdorff, it is a split monomorphism and there is a surjective homomorphism $\mu:\Cy(\cA(\cG)/\cA(\cG^{(0)}))\to \H^{\cyc}(\cG)$ such that $\mu\circ\iota=\iota'$.
\item There are quasi-isomorphisms
$$
\HC(\cG)\overset{\sim}{\to} \bigoplus_{n\ge 0}\H(\cG)[-2n],\,
 \HN(\cG)\overset{\sim}{\to}\prod_{n\ge 0}\H(\cG)[2n],\,\HP(\cG)\overset{\sim}{\to} \prod_{n\in\Z}\H(\cG)[2n].
 $$
\item Assume that $\cG$ is Hausdorff and $\cG^{\Iso}\setminus\cG^{(0)}$ discrete. Let $\cR$ be a full set of representatives of the orbits of the elements $x\in X$ with $\cG^x_x\ne\{x\}$. For each $x\in\cR$, choose a set $Z_x$ of representatives of the non-trivial conjugacy classes of $\cG^{x}_{x}$. We have a quasi-isomorphism of cyclic modules
\[
\H(\cG)\oplus \bigoplus_{x \in \cR} \bigoplus_{\eta \in Z_x} \H((\cG^{x}_{x})_\eta)\weq\H^{\cyc}(\cG).
\]
\end{enumerate}
\end{thm}
 
Theorem \ref{thm:burghintro} follows from Theorems \ref{thm:HC-main} and 
\ref{thm:H-cyc}, Propositions \ref{prop:twisted-hh-rel}, \ref{prop:agzcyc} and \ref{prop:ppal} , and Corollaries \ref{coro:twisted-retract}, \ref{coro:hadjG}, and \ref{coro:wsplit}. In the group case, the map of part 
i) of Theorem \ref{thm:burghintro} is the identity and the map 
$\mu$ of ii) is an isomorphism, while $\cC_c(G^{\Iso})=k[G]^{\adj}$ is the adjoint representation. Thus, the first assertion of part ii) recovers the well-known isomorphism (see Section \ref{sec:first})
\begin{equation}\label{eq:introadj}
HH_*(k[G])\cong H_*(G, k[G]^{\adj}).    
\end{equation}
The fact that $H_*(G)$ is a direct summand of $HH_*(k[G])$ and the computations of iii) when $\cG$ is a group go back at least to 
Karoubi's monograph 
\cite{karast}*{2.21-2.26}. In view of \eqref{eq:introadj}, part iv) specializes to Burghelea's theorem \cite{burghelea}*{Theorem I'} in the group case (Remark \ref{rem:burgh}).

When  $\cG$ is Hausdorff, there is a restriction map
$\res:\HH(\cA(\cG)/\cA(\cG^{(0)}))\to \H(\cG)$ that is left inverse to the embedding of part ii) of 
Theorem \ref{thm:burghintro}; composing it with the Dennis trace $K_n(\A(\cG))\to 
HH_n(\A(\cG))$ we get a map
\begin{equation}\label{map:bardintro}
    \olD_n=\res\circ D_n:K_n(\A(\cG))\to H_n(\cG).
\end{equation}
\comment{In fact we introduce a homology for pairs $(\cG,\omega)$ of an ample groupoid and a $2$-coycle with values in the invertible elements of $k$, and show it is a direct summand of the Hocschild homology of the twisted Steinberg algebra (Corollary \ref{coro:twisted-retract}). Thus a version of the map \eqref{map:bardintro} is also defined in the twisted case. 
}
Next we concentrate on Exel-Pardo groupoids, compute the Hochschild homology of their Steinberg algebras, and in the Hausdorff case use the splitting described above to also compute their groupoid homology, and the maps $\olD_*$ above to relate the latter to $K$-theory. 

In \cite{ep}, Ruy Exel and Enrique Pardo associate a groupoid  $\cG(G,E,\phi)$ to an action of a group $G$ on a (directed) graph $E$ by graph automorphisms and a $1$-cocycle $\phi:G\times E^1\to G$. For most results of the article we assume that $E$ is row-finite and that $G$ acts trivially on its set of vertices $E^0$. As in \cite{eptwist}, we additionally consider another $1$-cocycle taking values in the group of invertible elements of $k$, 
$c:G\times E^1\to\cU(k)$; the latter induces a groupoid $2$-cocycle 
$\bomega=\bomega_c:\cG(G,E,\phi)^{(2)}\to\cU(k)$  and we write $\cG(G,E,\phi_c)$ for the pair $(\cG(G,E,\phi),\bomega)$. The \emph{twisted Exel-Pardo $k$-algebra}
$L(G,E,\phi_c)$ is the twisted Steinberg $k$-algebra of $\cG(G,E,\phi_c)$. To better capture the effect of the cocycle $c$, which takes values in $k$, and so as to let \eqref{map:bardintro} be nontrivial on elements coming from $K_*(k)$, we consider Hochschild homology over a subring $\ell\subset k$ such that  $k/\ell$ is a flat ring extension (e.g. we could take $k=\C$ and $\ell=\Z$ or
$\Q$). Theorem \ref{thm:introep1} computes the Hochschild homology of $L=L(G,E,\phi_c)$ as an $\ell$-algebra, $HH_*(L/\ell)$, its homotopy algebraic $K$-theory $KH_*(L)$ and, under further assumptions, also its 
(Quillen) $K$-theory and the twisted groupoid homology $H_\ast(\cG(G,E,\phi_c),k/\ell)$ relative to the extension $k/\ell$. The latter is defined by a complex $\H(-,k/\ell)$ introduced in Definition 
\ref{defi:twistedhom}. When the cocycle is trivial, $H_*(\cG,k/\ell)=HH_*(k/\ell)\otimes_\ell H_*(\cG(G,E,\phi),\ell)$, the tensor product of graded $\ell$-modules. Proposition \ref{prop:twisted-hh-rel} shows that $
\HH(\cA_k(\cG)/\cA_\ell(\cG^{(0)}))$ is quasi-isomorphic to $\HH(\cA_k(\cG)/\ell)$ and Corollary \ref{coro:twisted-retract} proves that if $\cG$ is Hausdorff then $\H(\cG,k/\ell)$ is a direct summand of $\HH(\cA_k(\cG)/
\cA_\ell(\cG^{(0)}))$   . 

Since $L$ is $\Z$-graded, we have a weight decomposition 
\[
\HH(L)=\bigoplus_{m\in\Z}{}_m\HH(L)
\]
into a direct sum of chain complexes. 
Let $\reg(E)\subset E^0$ be the set of vertices that emit a finite nonzero number of edges. We introduce a $k[G]$-bimodule $S^c_m$ for each $0\ne m\in\Z$, and define chain maps

\begin{gather}
\sigma_m:\HH(k[G]/\ell,S^c_m)\to \HH(k[G]/\ell,S^c_m)\,\, (m\ne 0)\nonumber\\
\sigma_0:\HH(k[G]/\ell)^{(\reg(E))}\to \HH(k[G]/\ell)^{(E^0)}\label{maps:introsigmatau}\\
\tau:\H(G,k/\ell)^{(\reg(E))}\to \H(G,k/\ell)^{(E^0)}\nonumber
\end{gather}
given by explicit formulas that encompass information about the graph and the cocycles $\phi$ and $c$. 
Similarly, we define a map of spectra \eqref{map:elPhi}
\[
\Phi^t:KH(k[G])^{(\reg(E))}\to KH(k[G])^{(E^0)}
\]
induced by a zig-zag of explicit algebra homomorphisms. 
We say that a twisted Exel-Pardo triple $(G,E,\phi_c)$ is \emph{pseudo-free} if
$g(e)=e$ with $g\ne 1$ and $e\in E^1$ implies that $\phi(g,e)\ne 1$; in this case $\cG(G,E,\phi)$ is Hausdorff \cite{ep}*{Proposition 12.1}. If in addition $k[G]$ is regular supercoherent (e.g. if it is regular Noetherian) and $G$ acts trivially on $E^0$, then 
$L(G,E,\phi_c)$ is $K$-regular by \cite{eptwist}*{Corollary 8.17}), and thus the canonical map $K_*(L(G,E,\phi_c))\to KH_*(L(G,E,\phi_c))$ is an isomorphism.

The following is another main theorem of this article. It includes Theorem \ref{thm:hhep} (rephrased using the isomorphisms of Remark \ref{rem:smclosed}), Theorem \ref{thm:hep}, Corollary \ref{coro:kep} and Lemma \ref{lem:conmutrace}. Therein and elsewhere we shall abuse notation and write $I$ both for the identity map and for the $E^0\times \reg(E)$-matrix that results from the identity matrix after eliminating the columns corresponding to non-regular vertices.

\begin{thm}\label{thm:introep1} 
Assume that $E$ is row-finite and that $G$ acts trivially on $E^0$. Let $L=L(G,E,\phi_c)$ and $\cG=\cG(G,E,\phi_c)$, and let $\sigma_m$ and $\tau$ be as in \eqref{maps:introsigmatau}.
\begin{enumerate}[i)]
\item For each $m\in\Z$ there is a natural zig-zag of quasi-isomorphisms

\goodbreak

\noindent $\cone(I-\sigma_m)\weq {}_m\HH(L/\ell)$. In particular, there  is a long exact sequence
\begin{equation}\label{seq:hhintro}
\xymatrix{{}_0HH_{n+1}(L/\ell)\ar[r]&HH_n(k[G]/\ell)^{\reg(E)})\ar[d]^{I-\sigma_0}
\\
{}_0HH_n(L/\ell)&HH_n(k[G]/\ell)^{(E^0)}\ar[l]}
\end{equation}
\item We have a long exact sequence of homotopy algebraic $K$-theory groups
\begin{equation}\label{seq:khintro}
\xymatrix{KH_{n+1}(L)\ar[r]&KH_n(k[G])^{(\reg(E))}\ar[d]^{I-\Phi^t}
\\
KH_n(L)&KH_n(k[G])^{(E^0)}\ar[l]}
\end{equation}

\item If $(G,E,\phi)$ is pseudo-free and $k[G]$ is regular supercoherent, then we may substitute $K$ for $KH$ in the sequence \eqref{seq:khintro} and we have a commutative diagram with exact rows
\begin{equation}\label{diag:k&hintro}
\xymatrix{K_{n+1}(L)\ar[r]\ar[d]^{\olD_{n+1}}& K_n((k[G])^{(\reg(E))})\ar[d]^{\olD_n}\ar[r]^(.55){I-\Phi^t}&K_n(k[G])^{(E^0)}\ar[d]^{\olD_n}\ar[r]&K_n(L)\ar[d]^{\olD_n}\\
H_{n+1}(\cG,k/\ell)\ar[r]&H_n(G,k/\ell)^{(\reg(E))}\ar[r]^{I-\tau}&H_n(G, k/\ell)^{(E^0)}\ar[r]& H_n(\cG,k/\ell)}
\end{equation}
\end{enumerate}
\end{thm}
When the group $G$ is trivial, $L$ is the Leavitt path algebra of $E$, and parts i) and ii) of the theorem above
recover \cite{aratenso}*{Theorem 4.4} and \cite{abc}*{Theorem 8.6}. By Theorem \ref{thm:burghintro}, when $\cG$ is 
Hausdorff, each of the terms in the sequence of twisted groupoid homology relative to $k/\ell$ appearing in part iii) 
is a direct summand of the corresponding term in the exact sequence for Hochschild homology of part i). The pseudo-freeness hypothesis 
(which implies Hausdorffness) was used to show that the Hochschild homology exact sequence restricts 
to an exact sequence of twisted groupoid homology, and in particular that restriction gives a map of exact 
sequences between the two. This is used to show that the diagram of part iii) commutes. Using different techniques, A. 
Miller and B. Steinberg have recently shown \cite{ms} that in the untwisted case with $k=\ell$, the exact sequence of 
groupoid homology of part iii) holds without the pseudo-freeness hypothesis and even more generally for groupoids 
associated to self-similar actions of groupoids instead of groups. In the particular case of Katsura groupoids the 
latter sequence had been obtained by Ortega in \cite{homology-katsura}. Part ii) of Theorem \ref{thm:introep1} uses the 
computations of \cite{eptwist}*{Proposition 6.2.3 and Theorem 6.3.1}. The main novelty of this part of the theorem is 
the explicit description of the map $\Phi^t$ for general twisted Exel-Pardo groupoids (see \eqref{map:elPhi}); the 
particular case of twisted Katsura groupoids had been worked out in \cite{eptwist}*{Theorem 7.3}.

Several consequences of Theorem \ref{thm:introep1} are studied  in Section \ref{subsec:dtrace}. Theorem \ref{thm:introep2} below illustrates some of them. It includes all or part of Theorem \ref{thm:k01ep}, Proposition \ref{prop:Dcomm} and Corollaries \ref{coro:hep} and \ref{coro:matui}. 

Recall that the \emph{reduced incidence matrix} of a graph $E$ is the matrix $A=A_E\in\N_0^{(\reg(E)\times E^0)}$ whose $(v,w)$ entry is the number of edges $e$ with source $v$ and range $w$. The \emph{Bowen-Franks} group of $E$ is
\[
\BF(E)=\coker(I-A_E^t).
\]
\begin{thm}\label{thm:introep2}
Let $k$ be a field or a PID, $G$ a torsionfree group satisfying the Farrell-Jones conjecture such that $k[G]$ is regular supercoherent, $E$ a row-finite graph, and $(G,E,\phi_c)$ a pseudo-free Exel-Pardo tuple where $G$ acts trivially on $E^0$. Put $L=L(G,E,\phi_c)$. 

\begin{enumerate}[i)]
\item $K_0(L)=\BF(E)$, and $D_0$ is the composite of the inclusion and the scalar extension
\[
D_0:K_0(L)=\BF(E)\to \BF(E)\otimes_{\Z}k=H_0(\cG^{\bomega},k/\ell)\subset HH_0(L/\ell).
\]
In particular $\olD_0$ induces an isomorphism $K_0(L)\otimes_\Z k\iso H_0(\cG^{\bomega},k/\ell)$.

\smallskip

\item If $c$ is trivial and $k/\Z$ is flat, then there is a short exact sequence
\[
0\to \cU(k)\otimes_\Z\BF(E)\otimes k \to K_1(L)\otimes_\Z k\overset{\olD_1}{\lra}H_1(\cG,k)\to 0.
\]
\end{enumerate}
\end{thm}

In the case $G=\{1\}$, the identity $K_0(L)=\BF(E)$ is of course classical (see e.g. \cite{lpabook}*{Theorem 6.1.9}), the identity $H_0(\cG)=\BF(E)\otimes_{\Z}k$ is a particular case of Ortega's calculations \cite{homology-katsura}, the fact that it embeds into $HH_0(L)$ follows from \cite{aratenso}*{Theorem 4.4}, and the fact that $D_0$ is the scalar extension follows by direct computation. 

Motivated by part i) of Theorem \ref{thm:introep2} and the Bass trace conjecture for groups \cite{loday}*{8.5.2}, we propose the following. 

\begin{conj} Let $\cG$ be an ample Hausdorff groupoid. Then the image of the Dennis trace $D_0:K_0(\cA_{\Z}(\cG))\to HH_0(\cA_{\Z}(\cG))$ is contained in the direct summand $H_0(\cG,\Z)\subset HH_0(\cA_{\Z}(\cG))$.
\end{conj}    

In \cite{xlinotes}, Xin Li formulates a version of the Farrell-Jones conjecture for Steinberg algebras of torsionfree ample groupoids over noetherian regular coefficient rings. We explain in \ref{rem:compali} that part ii) of Theorem \ref{thm:introep2} is evidence in favor of that conjecture. Further connections with \cite{xlinotes} are discussed in Section \ref{sec:disc}, where another conjecture, Conjecture \ref{conjdisc}, pertaining to discretization invariance, is formulated.

The rest of this paper is organized as follows. Section \ref{sec:prelis} recalls basic definitions, facts and notation; it also contains the elementary technical Lemmas \ref{lem:res-sec}, \ref{lem:slice-tuples} and \ref{lem:hhsep}. In Section \ref{sec:hh-stein}  we establish Proposition \ref{prop:twisted-hh-rel} which says that if $(\cG,\omega)$ is a twisted ample groupoid with unit space $\cG^{(0)}=X$ and $k/\ell$ is a flat ring extension, then the canonical surjective homomorphism $\Cy(\cA_k(\cG,\omega)/\ell)\onto \Cy(\cA_k(\cG,\omega)/\cA_\ell(X))$ is a quasi-isomorphism. Part i) of Theorem \ref{thm:burghintro} follows from this by specializing to the case when $k=\ell$ and $\omega$ is trivial. Then in Definition \ref{defi:twistedhom} we introduce the twisted groupoid homology complex $\H(\cG^\omega,k/\ell)$, and show in Corollary \ref{coro:twisted-retract} that it embeds as a sub-semi-cyclic module of $\Cy(\cA_k(\cG)/\cA_\ell(X))$ which is a direct summand whenever $\cG$ is Hausdorff. Again, specializing to $k=\ell$ and trivial $\omega$ we get the embedding $\iota$ of part ii) of Theorem \ref{thm:burghintro} and the fact that it is split when $\cG$ is Hausdorff. Section \ref{sec:hc} contains Theorem \ref{thm:HC-main}, which establishes part iii) of Theorem \ref{thm:burghintro}. The rest of Theorem \ref{thm:burghintro} is proved in Section \ref{sec:first}. The identity $H^{\cyc}_*(\cG)=H_*(\cG,\cG^{\Iso})$ is established in Corollary \ref{coro:hadjG}. The maps $\mu$ and $\iota'$ of part ii) of the theorem are introduced in Proposition \ref{prop:agzcyc} and Corollary \ref{coro:wsplit}; the latter also establishes their relation with the map $\iota$. In the principal case, $\iota'$ is the isomorphism of Proposition \ref{prop:ppal}. Part iv) of Theorem \ref{thm:burghintro} follows from Theorem \ref{thm:H-cyc}.  Subsection \ref{subsec:sparse} specializes Theorem \ref{thm:H-cyc} to the case of ample Hausdorff transport groupoids 
$\cS\ltimes X$ associated to an action with sparse fixed points of an inverse semigroup $\cS$ on a locally compact Hausdorff space $X$. Section \ref{sec:ep} contains the proofs of 
Theorems \ref{thm:introep1} and \ref{thm:introep2}. Subsection \ref{subsec:graphs} recalls basic definitions, facts and notation on graphs and 
(twisted) Exel-Pardo groupoids. Subsection \ref{subsec:ept} contains two basic useful lemmas; Lemma \ref{lem:cohn=leav} and Lemma 
\ref{lem:EPcoli}. The first of these pertains to the (twisted) Steinberg algebra of the universal groupoid of the inverse semigroup $\cS(G,E,
\phi)$ associated to an Exel-Pardo tuple, and shows, among other things, that it coincides with the twisted Cohn algebra of \cite{eptwist}; 
this lemma is used later on in Subsection \ref{subsec:dtrace}, to establish the commutativity of the diagram of part iii) of Theorem 
\ref{thm:introep1}. 
The second lemma says that if $E$ is row-finite (each vertex emits finitely many edges) then the Exel-Pardo algebra $L(G,E,\phi_c)$ can be written as a colimit of EP-algebras over finite graphs; this is used in the Subsection \ref{subsec:hhep} to prove part i) of Theorem \ref{thm:introep1}. Subsection \ref{subsec:l0} studies the homogeneous component of degree $0$ of $L(G,E,\phi_c)$. The latter is an increasing union of subalgebras $L_0=\bigcup_{n\ge 0}L_{0,n}$ where $L_{0,n}$ is isomorphic to  sum of matrix algebras, indexed by the vertices $v\in E^0$, where the $v$-component consists of matrices with entries in $R_v$, the image of the map $k[G]\to L=L(G,E,\phi_c)$, $g\mapsto vg$. In general this map has a nonzero kernel $I_v$. However Proposition \ref{prop:Iv} gives useful technical information about $I_v$ and shows that $L_0$ can also be described as an ultamatricial algebra with coefficients in $k[G]$. In the next subsection we introduce the chain map 
$$
\sigma_m:\HH(k[G]/\ell,S_m^{\reg})\to \HH(k[G]/\ell,S_m)
$$ 
and show in Theorem \ref{thm:hhep} that ${}_m\HH(L/\ell)$ is quasi-isomorphic to the cone of $I-\sigma_m$. Part i) of Theorem 
\ref{thm:introep1} follows from this. For this result we use a description of the Hochschild homology of a twisted Laurent polynomial algebra 
associated to a corner isomorphism, proved in Appendix \ref{sec:appa}. The main result of Subsection \ref{subsec:hep} is Theorem \ref{thm:hep}, 
which says that if $(G,E,\phi_c)$ is pseudo-free, then for the twisted groupoid $\cG=\cG(G,E,\phi_c)$, $\H(\cG,k/\ell)$ is quasi-isomorphic to 
the cone of the restriction $I-\tau$ of $I-\sigma_0$ to the subcomplex $\H(G,k/\ell)^{\reg(E)}\subset \HH(k[G]/\ell)$. The exactness of the 
sequence of (twisted) groupoid homology groups of part iii) of Theorem \ref{thm:introep1} follows from this, and implies that 
$H_0(\cG,k/\ell)=\BF(E)\otimes_\Z k$ (Corollary \ref{coro:h0hep}). The next subsection contains Corollary \ref{coro:kep}, which 
establishes part ii) of Theorem \ref{thm:introep1}, and also the exact sequence of $K$-groups of iii), since under the hypothesis therein we 
can subsitute $K$ for $KH$ by \cite{eptwist}*{Corollary 8.17}. Corollary \ref{coro:kep} is deduced from Theorem \ref{thm:kep}, which says that 
if $\cT$ is a triangulated category and $\cH\aha\colon\to\cT$ is a homotopy invariant, excisive functor which is matricially stable and commutes with direct sums of sufficiently high number of summands (depending on $E$), then there is a distinguished triangle
\begin{equation}\label{seq:introtri}
\cH(k[G])^{(\reg(E))}\overset{I-\cH(\Phi^t)}\lra \cH(k[G])^{(E^0)}\to \cH(L_k(G,E,\phi_c)).
\end{equation}
Theorem \ref{thm:k01ep} of the same subsection says 
that under the hypothesis of part i) of Theorem 
\ref{thm:introep2}, we have $K_0(L)=\BF(E)$, and gives 
a short exact sequence computing $K_1(L)$. Theorem 
\ref{thm:kepz} describes the map $I-\Phi^t$ of part ii) 
of Theorem \ref{thm:introep1} in the particular case 
when $G=\Z$, and recovers the computation of $KH$ of twisted Katsura algebras \cite{eptwist}*{Theorem 7.3}.
Subsection \ref{subsec:dtrace} is concerned with the 
map \eqref{map:bardintro}. Lemma \ref{lem:conmutrace} 
shows that the diagram of part iii) of Theorem 
\ref{thm:introep1} commutes, concluding the proof of 
that theorem. Proposition \ref{prop:Dcomm} says that 
under the hypothesis of Theorem \ref{thm:introep2}, 
$D_0(K_0(L))\subset 
H_0(\cG,k/\ell)=\BF(E)\otimes k
\subset HH_0(L/\ell)$ and that $\olD_0$ induces an 
isomorphism $K_0(L)\otimes k\cong  H_0(\cG,k/
\ell)$, which completes the proof of part i) of Theorem 
\ref{thm:introep2}. The proposition also contains a 
description of the diagram of part iii) of Theorem 
\ref{thm:introep1} for $n=1$ which is used in Corollary 
\ref{coro:matui} to establish part ii) of Theorem 
\ref{thm:introep2}. Section \ref{sec:disc} concerns the universal groupoid $\cG_u(\cS)$ of an inverse semigroup $\cS$, and its discretization
$\cG_d(\cS)$. Xin Li's groupoid version of the Farrell-Jones conjecture mentioned above implies that if $\cG_u(\cS)$ is torsionfree and $k$ Noetherian regular, then $K_*(\A(\cG_u(\cS)))\cong K_*(\A(\cG_d(\cS)))$. Let $\cT$ be a triangulated category and $\cH:\aha\to \cT$ a functor. Assuming that $\cH$ is matricially stable on algebras with local units, we define a natural map  
\begin{equation}\label{map:rodéintro}
\tilde{\rho}_d:\cH(\A(\cG_d(\cS)))\lra\cH(\A(\cG_u(\cS))).    
\end{equation}

We call $\cH$ \emph{discretization invariant} if the latter map is an isomorphism for all $\cS$. We show in Proposition \ref{prop:nodisc} that $HH$ is not discretization invariant. Proposition \ref{prop:sidisc} says that if $\cH$ satisfies the hypothesis of \eqref{seq:introtri} and $(G,E,\phi)$ is an Exel-Pardo tuple, then \eqref{map:rodéintro} is an isomorphism for $\cS=\cS(G,E,\phi)$. Based on this we conjecture (Conjecture \ref{conjdisc}) that any functor $\cH:\ahal\to\cT$ that is excisive, homotopy invariant, matricially-stable and infinitely additive must be discretization invariant. 

Finally, Appendix \ref{sec:appa} is about the Hochschild homology of the twisted Laurent polynomial algebra $S=R[t_+,t_-,\phi]$ associated to a corner isomorphism $\phi:R\iso \phi(1)R\phi(1)$, introduced in \cite{fracskewmon}. Proposition \ref{prop:skewcorner} shows that for each $m\in\Z$, ${}_m\HH(S)$ is quasi-isomorphic to the cone of a certain endomorphism of $\HH(R,S)$. For example, the Exel-Pardo algebra $L=L(G,E,\phi_c)$ with $E$ finite without sources and $G$ acting trivially on $E^0$ is a twisted Laurent polynomial over $L_0$; Proposition \ref{prop:skewcorner} is used in the proof of Theorem \ref{thm:hhep}, which establishes part i) of Theorem \ref{thm:introep1}. 

\begin{ack}
 The second named author wishes to thank Xin Li for sharing his article \cite{xlinotes} and for useful email interchanges and several (in person and online) discussions. Thanks also to Pere Ara and Valentín Nico for their comments on previous versions of this article, and to Alistair Miller for pointing us towards his article with Ben Steinberg \cite{ms}, cited above.
\end{ack}

\section{Preliminaries}\label{sec:prelis}

We write
$\N = \{1, 2, 3, \dots\}$ and $\N_0 = \{0\} \cup \N$. 
Throughout the text we fix a commutative unital ring $\gring$. 
A $\gring$-bimodule $M$ is \emph{symmetric} if $\lambda x=x\lambda$ for all $x\in M$ and $\lambda\in k$. By an \emph{algebra} over $\gring$ we understand an associative ring $A$ with a structure of symmetric $\gring$-bimodule so that the multiplication map $A \otimes_\Z A\to A$, $a\otimes b\to ab$ induces a $\gring$-bimodule homomorphism $A\otimes_k A\to A$.

In this article, a \emph{compact} topological space is a Hausdorff space in which every open cover has a finite subcover. 

Let $f:X\to Y$ be a continous function. We say that $f$ is \emph{étale} if it is a local homeomorphism, and \emph{proper} if $f^{-1}(K)$ is compact for every compact subspace $K\subset Y$. 

If $\sigma:E\to X\leftarrow F:\tau$ are continuous maps we write
$$E\times F\supset E{ _\sigma\times}_{\tau}F=\{(e,f)\colon \sigma(e)=\tau(f)\}$$ 
for the pullback. 

\subsection{Groupoids}\label{subsec:gpd}
\numberwithin{equation}{subsection}
A (topological) \emph{groupoid} $\cG$ is a topological space together with a distinguished subspace $\cG^{(0)} \subset \cG$ of \emph{units} or objects, continuous \emph{source and range} maps $r,s\colon \cG\to \cG^{(0)}$, and composition and inverse maps 
\begin{gather*}
  \cG^{(2)} := \cG{_s \times_r } \cG\to \cG,\, (g,h)\mapsto gh,\\
\cG\to \cG,\, g\mapsto g^{-1},
\end{gather*}
satisfying the expected compatibility conditions. Groupoid
homomorphisms are continuous maps preserving compositions. We refer to \cite{xlispectra}*{Sections 2.1 and 2.2} for a succint introduction to topological groupoids; see also \cite{steinappr}*{Section 3} and \cite{exel}*{Section 3}. Throughout this text, the unit space $\cG^{(0)}$ will often be called $X$ and will always be 
assumed to be Hausdorff.
We say that a groupoid is \emph{étale} if the source (and, equivalently, the range) map is étale. A \emph{bisection} (or \emph{slice}) is a subset $U \subset \cG$ such that $s|_U$ and $r|_U$ are injective.
An étale groupoid is \emph{ample} if its compact open bisections form a basis of its topology. 

For a subset $Z \subset X$, we write $\cG^Z = s^{-1}(Z)$ and $\cG_Z = r^{-1}(Z)$. When $Z$ is a singleton, we omit the braces; we write $\cG^z = \cG^{\{z\}}$, $\cG_z = \cG_{\{z\}}$ and $\cG_z^z = \cG_z \cap \cG^z$. 
Observe that $\cG_z^z$ is a group with neutral element $z$; we call it the \emph{isotropy group} of $\cG$ at $z$. We say that $z$ has \emph{trivial isotropy} if $\cG_z^z=\{z\}$.
The \emph{isotropy} of $\cG$ is the subgroupoid
\[
\cG\supset\cG^{\Iso} = \{\eta \in \cG: s(\eta) = r(\eta)\} = \bigsqcup_{x \in X} \cG_x^x.
\]

Let $\Lambda$ be a discrete abelian group. A $\Lambda$-\emph{graded} groupoid is a groupoid $\cG$ together with a continuous groupoid homomorphism $|\cdot| \colon \cG \to \Lambda$
called the \emph{grading} or \emph{cocycle}.  

\subsection{\topdf{$\cG$}{G}-spaces}\label{subsec:gspace}

Let $\cG$ be an étale groupoid. A \emph{left $\cG$-space} is a topological space $Z$ together with a continuous map $\tau \colon Z \to X$, 
called the \emph{anchor map}, and a continuous \emph{action map} $\acts \colon \cG{_s \times_\tau } Z \to Z$ such that
\begin{enumerate}[i)]
    \item $\tau(g \acts z) = r(g)$ for each $z \in Z$ and $g \in \cG^{\tau(z)}$;
    \item $\tau(z) \acts z = z$ for all $z \in Z$;
    \item $g \acts (h \acts z) = gh \acts z$ for each $z \in Z$ and 
    each composable pair $(g,h) \in \cG^{(2)}$ such that $h \in \cG^{\tau(z)}$.
\end{enumerate}
The notion of right $\cG$-space is defined analogously.

If $\cG$ comes equipped with a $\Lambda$-grading, we 
define a \emph{graded} (left) $\cG$-space as a $\cG$-space $Z$ together with a continuous \emph{grading} $| \cdot | \colon Z \to \Lambda$ such that
$|g \bullet z| = |g| + |z|$ for each $z \in Z$ and $g\in \tau^{-1}(z)$.

\begin{ex} Any groupoid $\cG$ acts on itself by left multiplication, i.e. $g \acts h = gh$
for each pair of composable arrows. 
\end{ex}

\begin{ex} A groupoid $\cG$ acts on $\cG^{\Iso}$ by conjugation: we define $\tau = s \colon \cG^{\Iso} \to  X$
and $g \acts \eta = g\eta g^{-1}$.
\end{ex}

Given a $\cG$-space $Z$, the relation $x \sim y$ if
$x = g \acts y$ for some $g \in \cG$ is an equivalence relation on $Z$; we write $Z/\cG$ for the resutling quotient space. The \emph{orbit} of $x \in Z$ is its equivalence class with respect to this relation, denoted by $\cG\acts x$.

\subsection{Compactly supported functions}\label{subsec:ccx}

All spaces considered in this paper are locally compact. Such a space is \emph{weakly Boolean} if its compact open subsets form a basis of the topology, and \emph{generalized Boolean} if, in addition, it is Hausdorff.In \cite{steinappr}, generalized Boolean spaces are called locally compact Boolean.
For a weakly Boolean space $X$, we define
\[
\cc(X)=\mspan\{\chi_K\colon X\supset K \text{ compact open}\} \subset k^X.
\]
Remark that if in addition $X$ is Hausdorff, and we give $\gring$ the discrete topology, then $\cc(X)$ identifies with the set of compactly supported continuous functions $X\to \gring$, and the pointwise operations make the latter into a $\gring$-subalgebra of $\gring^X$.

We now recall how the construction $\cc(-)$ is functorial for proper maps and for étale maps. If $f \colon X \to Y$ is proper, composition with $f$ defines a $\gring$-linear map:
\begin{equation} \label{def:pullback}
f^\ast \colon \cc(Y) \to \cc(X), \qquad \chi_K \mapsto \chi_{f^{-1}(K)}.
\end{equation}

If $f \colon X \to Y$ is étale, then the following is a well-defined
$\gring$-linear map
\begin{equation}\label{def:pushforward}
f_\ast \colon \cc(X) \to \cc(Y), \qquad f_\ast(\phi)(x) = \sum_{z \in f^{-1}(x)} \phi(z). 
\end{equation}

\begin{ex} If $F \subset X$ is a closed subspace,
then the inclusion $i \colon F \to X$ is proper. If $X$ is weakly Boolean, the induced map will be denoted
$\res_{X,F} \colon \cc(X) \to \cc(F)$ since it maps $\chi_K$ to $\chi_{K \cap F}$
for each compact open subset of $X$. The subindices on $\res_{X,F}$ 
will be omitted when they can be deduced from the context.
\end{ex}

\begin{rmk} \label{rmk:pushforward}
Notice that if $f \colon X \to Y$ is étale
and $K \subset X$ a compact open such that $f$ is injective on $K$, i.e., 
such that $f|_K \colon K \to f(K)$ is a homeomorphism, then $f_\ast(\chi_K) = \chi_{f(K)}$.
\end{rmk}

The argument of \cite{steinappr}*{Proposition 4.3} also proves the lemma below; we include a proof for completeness.

\begin{lem} \label{lem:gen-int}
Let $X$ be a weakly Boolean space and $\mathcal{B}$ a
basis of compact open sets; then we have the following. 
\begin{enumerate}[i)]
\item 
$\cc(X) = \mspan\{ \chi_{\cap_{i=1}^n B_i} : B_i \in \mathcal{B} \text{ and $\cup_{i=1}^n B_i \subset Y\subset X$
with $Y$ Hausdorff}\}.$
\item If for every $B_1, \dots, B_n \in \mathcal{B}$ such that $\bigcup_{i=1}^nB_i$ is contained in a Hausdorff subspace of $X$ their intersection $B_1 \cap \cdots \cap B_n$ lies in $\mathcal{B}$, 
then $\cc(X) = \mspan\{\chi_B : B \in \mathcal{B}\}$.
\end{enumerate}
\end{lem}
\begin{proof} Item ii) follows directly from i); we prove the latter.
It suffices to prove that, for a compact open subset $O \subset X$,
the element $\chi_O \in \cc(X)$ lies in the span of the generators described in (i). Since $O$ is open, it is a union of elements of $\mathcal{B}$; further, since it is also compact, there exists finitely many $B_1, \ldots, B_n$ such that $O = B_1 \cup \cdots \cup B_n$. 
By the inclusion-exclusion principle,
\[
\chi_O = \chi_{B_1\cup \cdots \cup B_n}
= \sum_{i=1}^n (-1)^i \sum_{I\subset \{1, \ldots, n\},\ |I| = i} \chi_{\bigcap_{j\in I} B_j}.
\]
Given that $B_1, \dots B_n$ are contained in $O$, which is Hausdorff, each finite intersection in the right hand side is compact. Thus $\chi_{\bigcap_{j\in I} B_j} \in \cc(X)$ for all $I \subset \{1,\ldots, n\}$; this concludes the proof.
\end{proof}

\begin{rmk}\label{rmk:gen-int} We may apply Lemma \ref{lem:gen-int} ii), for example, to the basis of all compact open subsets
of a weakly Boolean space. It also applies to the set of all compact open slices of an ample groupoid.
\end{rmk}

\begin{lem} \label{lem:res-sec} Let $X$ be a generalized Boolean space and  $F\subset X$ a closed subspace.
Put $U = X \setminus F$ and let $i \colon U \to X$ be the inclusion. There
is a short exact sequence
\[
0 \to \cc(U) \xto{i_\ast} \cc(X) \xto{\res_{X,F}} \cc(F) \to 0.
\]
\end{lem}
\begin{proof} We have the formulas
\[
\res_{X,F}(\varphi) = \varphi|_F, \qquad i_\ast(\varphi)(x) = \begin{cases}
\varphi(x) & \text{if $x \in U$}\\
0 & \text{otherwise}
\end{cases}
\]
from which it follows that $\res_{X,F} \circ i_\ast = 0$ and that $i_\ast$ is injective.
Let $\varphi \in \cc(X)$. If $\varphi|_F = 0$, then the support of $\varphi$ is contained in $U$ and $\varphi|_U \in \cc(U)$. Because $X$ is Hausdorff, this implies that $\varphi = i_\ast(\varphi|_U)$, proving exactness at the middle of the sequence. Finally we turn to proving that $\res_{X,F}$ is surjective. Let $\mathcal{B}$ be a basis of compact open subsets of $X$; then $S = \{F \cap B : B \in \mathcal{B}\}$ is a basis of compact open subsets of $F$. Since $X$ is Hausdorff, so is $F$, hence $S$ lies in the hypothesis
of Lemma \ref{lem:gen-int} ii) and
$\cc(F) = \mspan \{\chi_{B \cap F} : B \in \mathcal{B}\} = \im(\res_{X,F})$.
\end{proof}

\subsection{Steinberg algebras}\label{subsec:steinalg} For an ample groupoid $\cG$, its \emph{Steinberg algebra}
(\cite{steinappr}, \cite{CFST}) is the $\gring$-module $\A(\cG) := \cc(\cG)$ equipped with the product
\[
(f_1 \ast f_2)(g) = \sum_{g = \alpha\beta} f_1(\alpha) f_2(\beta) \qquad (g \in \cG).
\] 
By \cite{steinappr}*{Proposition 4.3} $\cA(\cG)$ is generated as a $\gring$-module by the indicator functions of all of compact open bisections (see also Remark \ref{lem:gen-int}).
If $\cG$ is $\Lambda$-graded, there is an induced grading on $\A(\cG)$ via
\[
\A(\cG)_l = \{f \in \A(\cG) : \supp(f) \subset | \cdot |^{-1}(l) \} \quad (l \in \Lambda). 
\]

 \subsection{(Graded) \topdf{$\cG$}{G}-modules}\label{subsec:graded} Recall that a (left) module $M$ over a not necessarily unital ring $R$ is called \emph{unital} if $RM=M$. For an ample groupoid $\cG$, we shall study unital $\A(\cG)$-modules, which we will refer to as \emph{$\cG$-modules}. We write $\Mod_{\A(G)}$ for the category of $\cG$-modules. In this section we concentrate on left $\cG$-modules; right $\cG$-modules
are defined symmetrically.
A large family of examples stems from $\cG$-spaces; for any $\cG$-space $X$ 
with anchor map $\tau \colon X \to \cG^{(0)}$ the $\gring$-module $\cc(X)$
can be equipped with a $\cG$-module structure via
\[
\chi_U \cdot \chi_K := \chi_{UK}, \qquad UK = \{u \acts k : k \in K, u \in \cG^{\tau(k)}\cap U\}
\]
for any compact open sets $U \subset \cG$, $K \subset X$. When $\cG$ is $\Lambda$-graded and $X$
is a graded $\cG$-space, 
then $\cc(X)$ is $\Lambda$-graded via 
$\cc(X)_l = \{f \in \cc(X) : \supp(l) \subset |\cdot|^{-1}(l)\}$.

\subsection{Simplicial and cyclic weakly Boolean spaces}\label{subsec:weakbool}

Equipping weakly Boolean spaces with proper 
(resp. étale) maps, we obtain a contravariant (resp. covariant)
functor $X \mapsto \cc(X)$ taking values in $\gring$-modules. Write $\WB$
for the category of weakly Boolean spaces and étale maps.

A simplicial weakly Boolean space is a functor $X \colon \Delta^\acts \to \WB^\op$. 
It induces a simplicial $\gring$-module $\cc(X)$, and, in particular, a complex of $\gring$-modules with differentials
\[
\partial_n = \sum_{i=0}^n (-1)^i (d_i)_\ast.
\]
In this paper we will mainly be interested in two examples of this concept, associated to any ample groupoid $\cG$, that we proceed to describe below.

\begin{ex}[Nerve of a groupoid] \label{ex:nerve}
For each $n \ge 1$, write
\begin{equation}
    \cN(\cG)_n = \cG^{(n)}= \{(g_1,\dots,g_n) \in \cG^{n} \colon s(g_i)=r(g_{i+1}) \, \forall 1\le i\le n-1\}.\label{gn}
\end{equation} 
for the $n$-tuples of composable arrows of $\cG$, equipped with the subspace topology of the cartesian product $\cG^n$. Write also $\cN(\cG)_0 = \cG^{(0)}$. Because $\cG^{(0)}$ is Hausdorff, $\cG^{(n)}\subset \cG^n$ is closed. In particular,
if $A_1, \ldots, A_n\subset\cG$ are compact open bisections, 
the open subset
\begin{equation}\label{slice:bracket}
    [A_1 | \cdots | A_n] := (A_1 \times \cdots \times A_n) \cap \cG^{(n)}
\end{equation}
is also compact. These compact open subsets form a basis of $\cG^{(n)}$, proving that the latter space is weakly Boolean. 
For each $n \ge 0$ and $i \in \{0,\ldots, n\}$, put
\begin{align*}
&d_i \colon \cG^{(n)} \to \cG^{(n-1)}, \qquad d_i(g_1,\ldots,g_n) = \begin{cases}
(g_2,\ldots,g_n) & \text{if $i = 0$}\\
(g_1, \ldots, g_{n-1}) & \text{if $i=n$}\\
(g_1, \ldots, g_i g_{i+1}, \ldots, g_n) & \text{otherwise}
\end{cases}\\
&s_i \colon \cG^{(n)} \to \cG^{(n+1)}, s_i(g_1, \ldots, g_n) = 
\begin{cases}
(r(g_1),g_1,\ldots,g_n) & \text{if $i = 0$}\\
(g_1, \ldots, g_n,s(g_n)) & \text{if $i=n$}\\
(g_1, \ldots, g_{i}, r(g_{i+1}) , g_{i+1}, \ldots, g_n) & \text{otherwise}.
\end{cases}
\end{align*}

Further, one verifies 
that
\begin{align} \label{eq:d-s-slice}
d_0[A_1 | \cdots | A_n] &= [s(A_1)A_2 | \cdots | A_n],\\ 
d_n[A_1 | \cdots | A_n] &= [A_1 | \cdots | A_{n-1}r(A_n)],\notag\\ 
d_i[A_1 | \cdots | A_n] &= [A_1 | \cdots | A_i A_{i+1} | \cdots A_n],\notag\\
s_0[A_1 | \cdots | A_n] &= [r(A_1) | A_1 | \cdots | A_n],\notag\\
s_n[A_0 | \cdots | A_n] &= [A_0 | \cdots  | A_n | s(A_n)],\notag\\ 
s_i[A_1 | \cdots | A_n] &= [A_1 | \cdots | A_i | r(A_{i+1})  | A_{i+1} | \cdots A_n]\notag
\end{align}
and that $d_i$ and $s_i$ restricted to $[A_0 | \cdots | A_n]$ are injective, proving that all faces and degeneracies are étale maps. Hence $\cN(\cG)$ is a 
simplicial weakly Boolean space  in the sense defined above.

As a simplicial set $\cN(\cG)$ is isomorphic to 
the nerve $N(\cG)$ of $\cG$ viewed as a category. Since in the standard convention (see e.g. \cite{goejar}) maps point in the opposite direction as ours (which are oriented as in \cite{bouka}), the isomorphism must invert the maps. It
is given by the natural bijections
\[
(g_1, \ldots, g_n) \in \cN(\cG) \mapsto (g_1^{-1}, \ldots, g_n^{-1}) \in N(\cG). 
\]
As we shall recall below, the complex $\H(\cG) = \H(\cc(\cG^{(\acts)}))$ computes the homology 
of $\cG$ with coefficients in $\gring$.
\end{ex}

\begin{ex}[Cyclic nerve of a groupoid] \label{ex:cyc}
For each $n \ge 0$, we can consider the \emph{cyclically composable}
arrows
\[
\cG^{(n+1)}\supset\cG^n_{\cyc}=\{ (g_0,\dots,g_n)\in\cG^{(n+1)}\colon s(g_n)=r(g_0)\}.
\]
equipped with the subspace topology. This is a closed subspace because $\cG^{(0)}$ is Hausdorff. Each space $\cG_{\cyc}^n$ has a basis of compact open subsets
given by
\begin{equation}\label{slice:parent}
(A_0 | \cdots | A_n) = (A_0 \times \cdots \times A_n) \cap \cG_{\cyc}^n
\end{equation}
where $A_0, \ldots, A_n \subset \cG$ are compact open bisections. 
For each $n \ge 0$ and $i \in \{0,\ldots, n\}$, put
\begin{align*}
&d_i \colon \cG_{\cyc}^{n} \to \cG_{\cyc}^{n-1}, \qquad d_i(g_0,\ldots,g_n) = \begin{cases}
(g_ng_0,g_1,\ldots,g_{n-1}) & \text{if $i = n$}\\
(g_0, \ldots, g_i g_{i+1}, \ldots, g_n) & \text{otherwise}
\end{cases}\\
& s_i \colon \cG_{\cyc}^{n} \to \cG_{\cyc}^{n+1}, \qquad s_i(g_0, \ldots, g_n) = (g_0, \ldots, g_i, s(g_i), g_{i+1},\ldots, g_n).
\end{align*}
The maps $d_i$ and $s_i$ interact with a basic compact open set \eqref{slice:parent} in a way analogous to the identities \eqref{eq:d-s-slice}; hence they are étale.
We thus have a simplicial weakly Boolean space $\cG_{\cyc}$. 
In an abuse of notation, we write $\H^{\cyc}(\cG) = C_c(\cG_{\cyc}^{\acts})$ for both the associated simplicial $k$-module and its associated chain complex.  Starting in Example \ref{ex:cycobj-hg} below we shall further abuse notation and use the same name for the associated cyclic module. 
\end{ex}

\begin{rmk}\label{rmk:cyc-graded}
We point out that if $\cG$ is $\Lambda$-graded, then $\cG^n_{\cyc}$ is $\Lambda$-graded
with grading $|(g_0,\dots,g_n)| = |g_0| + \cdots + |g_n|$ and, since $\Lambda$ is assumed to be abelian, 
all face and degeneracy maps of the cyclic nerve construction are compatible with the grading.
Hence $\H^{\cyc}(\cG)$ is a simplicial $\Lambda$-graded $k$-module with all face and degeneracy maps homogeneous of degree zero.
\end{rmk}

We record the following straighforward lemma.

\begin{lem} \label{lem:slice-tuples}
Let $A_0, A_1, \ldots, A_n \subset \cG$ be compact open bisections and $U \subset \cG^{(0)}$ a 
compact open subset. We have the following equalities:
\begin{enumerate}[i)]
    \item $[A_1 | \ldots | A_i U | A_{i+1} | \cdots | A_n] = 
[A_1 | \ldots | A_i  |  UA_{i+1} | \cdots | A_n]$;
    \item $(A_0 | \ldots | A_i U | A_{i+1} | \cdots | A_n) = 
(A_0 | \ldots | A_i  |  UA_{i+1} | \cdots | A_n)$;
\item $(U A_0 | \cdots | A_n) = 
(A_1 | \cdots | A_0 U)$.
\end{enumerate}
\qed
\end{lem}

\subsection{Groupoid homology}\label{subsec:hgpd}

We now come to the definition of groupoid 
homology. We follow the presentation of 
\cite{miller-corre}*{Section 2}; see also \cite{xlispectra}*{2.3}. Fix an ample
groupoid $\cG$. 
Let $n\ge 0$; the $n^{th}$-\emph{homology} of $\cG$ with coefficients
in a $\cG$-module $M$ relative to $\gring$ is defined as
\[
H_n(\cG,M) = \tor^{\A(\cG)}_n(\cc(\cG^{(0)}),M). 
\]
We also write $H_\ast(\cG) := H_\ast(\cG,\cc(G^{(0)}))$
and $H_\ast(\cG,Z) = H_\ast(\cG,\cc(Z))$ for each $\cG$-space $Z$.
As observed in \cite{miller-corre}*{Section 2}
and the references therein, 
we shall use the fact that $\tor$ can be computed
via flat resolutions. Namely, if $P_\acts$
is a flat resolution of $M$, then 
$H_\ast(\cG,M)$ is the homology of $\cc(\cG^{(0)}) \otimes_{\A(\cG)} P_\acts$; likewise if we resolve $\cc(\cG^{(0)})$ by flat right $\cG$-modules and then tensor by $M$. We shall revise the construction of a concrete complex that computes groupoid homology using this fact.

First, we recall some useful results 
from \cite{miller-corre} on flatness and tensor product of  $\cG$-modules. A left $\cG$-space $Z$
is said to be \emph{basic} if the map
\[
\cG \times_{\cG^{(0)}} Z \to Z \times_{Z/\cG} Z, \qquad (g,x) \to (g \acts x, x).
\]
is a homeomorphism, and étale if its anchor map is étale.

It is straightforward to verify that 
$\cG^{(n)}$ is basic and étale for each $n \ge 1$. 
Our interest
in basic $\cG$-spaces lies in the following result.

\begin{prop}[\cite{miller-corre}*{Proposition 2.8}] \label{prop:basic-flat}
Let $\cG$ be an ample groupoid and let $Y$ be a basic
étale $\cG$-space. Then $\cc(Y)$ is a flat $\cG$-module.
\qed
\end{prop}

We abbreviate $\otimes_\cG:= \otimes_{\A(\cG)}$.
Given a left $\cG$-space $Z$ and a right $\cG$-space
$Y$, we may form the pullback $Y \times_{\cG^{(0)}} Z$ along their respective anchor maps; its quotient 
by the relation $(y \acts g,z) \sim (y,g \acts z)$
will be denoted $Y \times_\cG Z$.

\begin{prop}[\cite{miller-corre}*{Proposition 2.9}] \label{prop:fiber-tensor}
Let $\cG$ be an ample groupoid, let $Y$ be a basic étale right $\cG$-space with anchor map $\sigma\colon Y \to G^{(0)}$ let $Z$ be a totally disconnected left $\cG$-space. 
Then $Y \times_\cG Z$ is totally disconnected and locally compact, and there is an isomorphism $\kappa \colon \cc(Y) \otimes_{\cG} \cc(Z) \iso \cc(Y \times_{\cG} Z)$ given by 
\begin{equation}\label{def:kappa}
    \kappa(\xi \otimes \eta)([y,z]) = \sum_{g \in \cG^{\sigma(y)}} \xi(y \acts g)\eta(g^{-1} \acts z).
\end{equation}
\qed
\end{prop}

\begin{rmk} \label{rmk:pullback-tensorX}
Let $\cG$ be an ample groupoid, $Y$ an étale right $\cG$-space, and $Z$ 
a totally disconnected left $\cG$-space. Then $Y$ is basic and étale as a $\cG^0$-space. Hence Proposition \ref{prop:fiber-tensor} applied to $\cG^{(0)}$ in place of $\cG$ says that $\cc(Y) \otimes_{\cG^{(0)}} \cc(Z) \cong \cc(Y \times_{\cG^{(0)}} Z)$.
\end{rmk}

\begin{rmk} \label{rmk:pullback-graded}
In Proposition \ref{prop:fiber-tensor}, if $\cG$, $Y$ and $Z$ are $\Lambda$-graded, then 
$Y \times_{\cG} Z$ can be equipped with a $\Lambda$-grading via $|[y,z]| = |y|+|z|.$
With this grading the map $\kappa$ becomes homogeneous of degree zero: if $\xi\in \cc(Y)_l$ and
$\eta \in \cc(Z)_{l'}$ for some $l,l' \in \Lambda$, then for $\kappa(\xi \otimes \eta)([y,z])$
to be non-zero there must exist $g \in \cG^{\sigma(y)}$ such that $y \bullet g \in \supp(\xi)$
and $g^{-1} \bullet z \in \supp(\eta)$. Hence $|y|+|g| = l$, $-|g|+|z|=l'$, and thus 
$|[y,z]| = l+l'$. It follows that $\supp(\kappa(\xi \otimes \eta))$ is contained in $|\cdot|^{-1}(l+l')$
and thus $\kappa(\xi \otimes \eta) = l+l' = |\xi| + |\eta| = |\xi \otimes \eta|$ as claimed.
\end{rmk}

\begin{coro}\label{coro:agzbimo}
Let $\cG$ be an ample groupoid and $Z$ a topological space with right and left $\cG$-space structures. If $Z$ is totally disconnected, then the map
\begin{gather*}
\mu:\A(\cG)\otimes_{\cG^{(0)}}\cc(Z)\otimes_{\cG^{(0)}} \A(\cG)\to \cc(\cG\times_{\cG^{(0)}} Z \times_{\cG^{(0)}}\cG),\\
\mu(\phi_0\otimes\psi\otimes\phi_1)(g_0,z,g_1)=\phi_0(g_0)\psi(z)\phi_1(g_1).
\end{gather*}
is an isomorphism of bimodules.
\qed
\end{coro}

\begin{ex}[Bar and standard resolution] \label{ex:bar}
Write $B_n(\cG) = \cG^{(n+1)}$ for each $n \ge -1$,
and for each $n \ge 0$ define
\begin{align*}
d_i(g_0,\dots, g_n) &= (g_0, \dots, g_{i-1}g_{i}, \dots, g_n), \qquad 0 < i \le n,\\
d_0(g_0,\dots,g_n) &= (g_1\dots,g_{n}) \qquad n > 0,\\
s_i(g_0,\ldots,g_n) &= (g_0,\ldots,r(g_i),g_i,\ldots,g_n).
\end{align*}
At the level of $B_0(\cG)$
we define $d_0(g_0) = s(g_0)$.
A similar analysis as the one done for \eqref{ex:nerve}
shows that these are étale $\cG$-equivariant maps. 
We then have an associated complex $(\cc(B_\acts(\cG)),b_\acts)_{n \ge -1}$
with boundary
$b_n = \sum_{0 \le i \le n} (-1)^i (d_i)_\ast$.
Consider $h_n  \colon B_n(\cG) \to B_{n+1}(\cG)$, $h_n(g_0,\ldots,g_n) =(g_0,\ldots,g_n,s(g_n))$ and also the open inclusion
$h_{-1} \colon B_{-1}(\cG) \to B_0(\cG)$. These maps satisfy the relations
\[
d_{i} h_n = h_{n-1} d_{i}, \quad
d_{n+1} h_n = \id,
\quad d_0 h_{-1} = \id \quad (0 \le i \le n).
\]
It follows
that $\{(-1)^{n+1}(h_{n})_\ast\}_{n \ge -1}$ is a contracting homotopy of the complex 
$(\cc(B_\acts(\cG)),b_\acts)_{n \ge -1}$; whence the latter is (pure) exact. 
Thus by Proposition \ref{prop:basic-flat} have a flat resolution $\bB(\cG) := \cc(B_n(\cG))_{n \ge 0}$ 
of $\bB(\cG)_{-1} := \cc(\cG^{(0)})$.

It follows that the homology of $\bB(\cG) \otimes_\cG M$ computes $H_\ast(G,M)$.
When $M = \cc(Z)$ for some totally disconnected 
$\cG$-space $Z$, the using Proposition 
\ref{prop:fiber-tensor} for the first isomorphism, we have
\[
\bB(\cG)_n \otimes_\cG\cc(Z) \cong \cc(\cG^{(n+1)} \times_\cG Z) \cong \cc(\cG^{(n)} \times_{\cG^{(0)}} Z).
\]
Furthermore, the maps $(d_i)_\ast \otimes_\cG\cc(Z)$ are induced by 
the maps $\delta_i \colon \cG^{(n)} \times_{\cG^{(0)}} Z \to \cG^{(n-1)} \times_{\cG^{(0)}} Z$
given by 
\[
\begin{cases}
\delta_0(g_1,\ldots,g_n,z)  &= (g_2,\ldots,g_n,z)\\
\delta_i (g_1,\ldots,g_n,z) &=  (g_1,\ldots,g_ig_{i+1},\ldots, g_n,z) \, i < n\\
\delta_n(g_1, \ldots,g_n,z) &= (g_1, \ldots,g_{n-1},g_n \acts z).
\end{cases}
\] 
We write $\H(\cG,Z)$ for the resulting complex. As observed, its homology computes $H_\ast(\cG,\cc(Z))$ as defined above.
For $Z=\cG^{(0)}$, the complex $\H(\cG,\cG^{(0)})$ can be identified with the complex $\H(\cG)$ 
associated to the nerve of $\cG$ described in Example \ref{ex:nerve}.
\end{ex}

\subsection{Hochschild homology}\label{subsec:hh}

Let $A$ be a $\gring$-algebra. A \emph{system of local units} in $A$ is a set $\mathcal{E}\subset A$ of idempotent elements such that the set $\{pAp\colon p\in\mathcal{E}\}$, ordered by inclusion, is filtered and satisfies $\bigcup_{p\in\mathcal{E}}pAp=A$. We say that $A$ has \emph{local units} if it has a system of local units. 

Assume that $A$ has local units. Consider the $\N_0$-graded complex $\HH(A/\gring)$ given by the $\gring$-modules $\HH(A/\gring)_n=A^{\otimes_\gring n+1}$ together with boundary maps
\begin{multline}\label{map:hochb}
b(a_0\otimes\cdots\otimes a_n)= \sum_{i=0}^{n-1} (-1)^i a_0\otimes\cdots\otimes a_ia_{i+1}\otimes\cdots\otimes a_n  + (-1)^{n}a_{n}a_0\otimes a_1\otimes \cdots\otimes a_n.     
\end{multline}

We call $\HH(A/\gring)$ the \emph{Hochschild complex} and its homology $HH_\ast(A/\gring)$ the \emph{Hochschild homology} of $A$ (relative to $\gring$). 

\begin{rem}\label{rem:naive}
In \cite{loday}*{Section 1.4.3} the complex $\HH(A/\gring)$ is denoted $C^{\mathrm{naiv}}(A/\gring)$ and called the \emph{naive} Hochschild complex. For general $A$, its homology may differ from Hochschild homology as defined in \cite{loday}*{Section 1.4.1}; however both definitions agree when $A$ has local units, by \cite{loday}*{Propositions 1.4.4 and 1.4.8}.
\end{rem}

For a given $A$-bimodule $M$, we write $[M,A]$ for the $\gring$-linear span of all commutators $[m,a]=ma-am$ and 
\begin{equation}\label{eq:natural}
M_\sharp=M/[M,A]
\end{equation}
for the quotient $\gring$-module. Viewing $M$ as an left module over the enveloping algebra $A \otimes_k A^\op$, we have an isomorphism of $k$-modules \[M_\sharp \cong A \otimes_{A \otimes A^\op} M.\]
Let $B$ be another $k$-algebra such that $A\subset B$ is a subalgebra. We shall assume that $A$ contains a system of local units of $B$ (and thus also of $A$). Regard $B^{\otimes_{A} n+1}$ ($n\ge 0$) as an $A$-bimodule in the obvious way and put 
$$
\HH(B/A)_n=B^{\otimes_A n+1}_\sharp \cong A \otimes_{A \otimes A^\op} B^{\otimes_A n+1}.
$$
The Hochschild boundary map \eqref{map:hochb} descends to a map $b:\HH(B/A)_{*+1}\to \HH(B/A)_*$ that makes $\HH(B/A)$ into a chain complex. 

\begin{rmk} \label{rmk:hh-graded}
If $B$ is a $\Lambda$-graded algebra then $B^{\otimes_k n+1}$ is a graded
$k$-module. If $A \subset B_0$, then 
$B^{\otimes_A n+1}$ is also a graded $k$-module. In both cases the grading
is given on elementary tensors of homogeneous elements by $|b_0 \otimes \cdots \otimes b_n|
= |b_0| + \cdots + |b_n|$. The grading on $B^{\otimes_A n+1}$ descends 
to one on $B^{\otimes_A n+1}_\sharp$.
Hence both $\HH(B/A)$ and $\HH(B/k)$ are complexes
of $\Lambda$-graded modules with boundary maps that are homogeneous of degree zero,
and the canonical comparison map $\HH(B/k) \to \HH(B/A)$ is compatible with the
respective gradings. 
\end{rmk}

\begin{lem}\label{lem:hhsep}
Let $B$ be a $\gring$-algebra and $A\subset B$ a commutative $\gring$-subalgebra. 
Let $\mathcal{F}$ be the set of all finite sets of orthogonal idempotent elements of $A$. Assume that
\begin{enumerate}[i)]
    \item for each $a_1,\cdots, a_n\in A$, there exists $F\in\mathcal{F}$ such that $\{a_1,\cdots,a_n\}\subset \mspan F$;
    \item $A$ contains a system of local units of $B$. 
\end{enumerate}
Then the canonical projection 
$$\HH(B/\gring)\onto \HH(B/A)$$ 
is a quasi-isomorphism. 
\end{lem}
\begin{proof} 
For $F\in\mathcal{F}$, $kF := \mspan F\subset A$ is a unital subalgebra with unit $p_F=\sum_{p\in F}p$. Hypothesis i) implies that the system $\{kF\}_{F\in\mathcal{F}}$ is 
filtered and that $\bigcup_{F\in\mathcal{F}}kF=A$. By ii), there exists $\cE\subset A$ that is a system of local units for $B$; in particular $B=\bigcup_{p\in\cE}pBp$. By what we have just seen, for every $p\in\mathcal{E}$ there exists $F\in\mathcal{F}$ such that $p\in kF$; hence $p\in p_FBp_F$ and thus $pBp\subset p_FBp_F$. It follows that $B=\bigcup_{F\in \mathcal{F}}p_FBp_F$, so $\{p_F\colon F\in\mathcal{F}\}$ is a system of local units for $B$. Hence the inclusion $A\subset B$ 
is the colimit over $F\in\mathcal{F}$ of the inclusions $kF\subset B_F:=p_FBp_F$, the latter are unital $\gring$-algebra homomorphisms, and the map of the 
proposition is the colimit over $F\in\mathcal{F}$ of the projections $\HH(B_F/\gring)\to \HH(B_F/\gring F)$. Hence we may assume that $F$ is finite, $A=kF$ is a finite 
direct sum of copies of $\gring$, and the inclusion $A\subset B$ is a unital homomorphism of unital $\gring$-algebras. Under these assumptions, the statement of 
the lemma is a particular case of \cite{loday}*{Theorem 1.2.13}.
\end{proof}

\subsection{Cyclic homology}\label{subsec:hc}

In this section we give a brief account 
on the cyclic homology of (semi-) cyclic modules, following \cite{loday}*{Section 2.5}.

A \emph{cyclic $\gring$-module}
is a simplicial $k$-module $M_\acts$ equipped together with a~$\Z/(n+1)\Z$-action on
$M_n$ for each $n \ge 0$, given by homomorphisms $t_n \colon M_n \to M_n$ subject to the following
compatibility conditions:
\begin{align}
t_n^{n+1} &= \id,\label{id:tn=1}\\ 
d_i t_n &= -t_{n-1}d_{i-1} \text{ for } 1 \le i \le n, \label{id:cycfaces1}\\ 
d_0 t_n &= (-1)^n d_n,\label{id:cycfaces2}\\ 
s_i t_n &= -t_{n+1}s_{i-1} \text{ for } 1 \le i \le n,\nonumber\\
s_0 t_n &= (-1)^n t_{n+1}^2 s_n.\nonumber
\end{align}
A \emph{semicyclic $k$-module} (called precyclic module in \cite{loday}*{page 77}) is a semisimplicial $k$-module $M$ with operators $t_n$ as above, satisfying identities \eqref{id:tn=1}, \eqref{id:cycfaces1} and \eqref{id:cycfaces2}.

By definition every cyclic module is a semicyclic module. 
Our motivation to consider the latter stems from the following example.

\begin{ex}\label{ex:ccycR}
Let $R$ be a unital $k$-algebra. The standard cyclic $k$-module $C^{\cyc}(R)$ associated to $R$ 
\cite{loday}*{Proposition 2.5.4}
is the simplicial module underlying $\HH(R)$ together with the $\Z/(n+1)\Z$-action on $\HH(R)_n=R^{\otimes n+1}$ via permutation of tensors. The definition of the degeneracy operators depends upon the fact that $R$ is unital. 
For a non-unital algebra $A$, we can define the face maps and cyclic operators in the same fashion, thus making $C^{\cyc}(A)$ a semicyclic module. 
\end{ex} 

\begin{ex}\label{ex:b/acyclic}
Let $A\subset B$ be $k$-algebras as in Lemma \ref{lem:hhsep}. Then $B=\bigcup_{p\in\cF}pBp$ is a filtering union, and each corner $pBp$ with $p\in\cF$ is unital, so $C^{\cyc}(pBp)$ is a cyclic module, with degeneracies defined by inserting a $p$ in the appropriate place. If $p,q\in\cF$ and $pBp\subset qBq$, then for $a_0,\dots,a_n\in pBp$ we have
\begin{align*}
a_0\otimes\cdots\otimes a_i\otimes q\otimes a_{i+1}\otimes\cdots\otimes a_n=&a_0\otimes\cdots\otimes a_ip\otimes q\otimes a_{i+1}\otimes\cdots\otimes a_n\\
&=a_0\otimes\cdots\otimes a_i\otimes pq\otimes a_{i+1}\otimes\cdots\otimes a_n\\
&=a_0\otimes\cdots\otimes a_i\otimes p\otimes a_{i+1}\otimes\cdots\otimes a_n.
\end{align*}
Hence degeneracies are well-defined on $C^{\cyc}(B)=\colim_{p\in\cF}C^{\cyc}(pBp)$, and give it a cyclic module structure.
\end{ex}
Given a semicyclic module $M$, we define operators $b,b' \colon M_n \to M_{n-1}$ and $N \colon M_n \to M_n$ by $b = \sum_{i = 0}^n (-1)^n d_i$, $b' = \sum_{i = 0}^{n-1} (-1)^n d_i$ and $N = \sum_{i = 0}^n t^i$, which satisfy the relations $b(1-t) = (1-t)b'$ and $b'N = Nb$, 
thus assembling into a bicomplex $CC(M)$ with anticommuting differentials as follows:
\begin{equation} \label{cc-bicomplex}
\begin{tikzcd}
	{} & {} & {} & {} \\
	{M_2} & {M_2} & {M_2} & {M_2} & {} \\
	{M_1} & {M_1} & {M_1} & {M_1} & {} \\
	{M_0} & {M_0} & {M_0} & {M_0} & {}
	\arrow["b"', dotted, from=1-1, to=2-1]
	\arrow["{-b'}"', dotted, from=1-2, to=2-2]
	\arrow["b"', dotted, from=1-3, to=2-3]
	\arrow["{-b'}"', dotted, from=1-4, to=2-4]
	\arrow["b"', from=2-1, to=3-1]
	\arrow["{1-t}"', from=2-2, to=2-1]
	\arrow["{-b'}"', from=2-2, to=3-2]
	\arrow["N"', from=2-3, to=2-2]
	\arrow["b"', from=2-3, to=3-3]
	\arrow["{1-t}"', from=2-4, to=2-3]
	\arrow["{-b'}"', from=2-4, to=3-4]
	\arrow["N"', dotted, from=2-5, to=2-4]
	\arrow["b"', from=3-1, to=4-1]
	\arrow["{1-t}"', from=3-2, to=3-1]
	\arrow["{-b'}"', from=3-2, to=4-2]
	\arrow["N"', from=3-3, to=3-2]
	\arrow["b"', from=3-3, to=4-3]
	\arrow["{1-t}"', from=3-4, to=3-3]
	\arrow["{-b'}"', from=3-4, to=4-4]
	\arrow["N"', dotted, from=3-5, to=3-4]
	\arrow["{1-t}"', from=4-2, to=4-1]
	\arrow["N"', from=4-3, to=4-2]
	\arrow["{1-t}"', from=4-4, to=4-3]
	\arrow["N"', dotted, from=4-5, to=4-4]
\end{tikzcd}
\end{equation}

The Hochschild homology  $H_*(M)$ is that of the complex $\HH(M)=(M,b)$. The \emph{cyclic homology} of $M$ is the homology of the totalization of $CC(M)$,
\begin{equation}\label{def:HC}
\HC(M) = \Tot(CC(M)), \qquad HC_n(M) = H_n(\HC(M)).
\end{equation}
Remark that te bicomplex above can be extended, by repeating columns infinitely to the left, to obtain an upper half-plane bicomplex $CC^{\mathrm{per}}(M)$, of which the second quadrant truncation is a subcomplex $CC^-(M)$. The \emph{periodic} and \emph{negative cyclic} complexes of $M$ are the direct product totalisations 
$\HP(M)=\Tot (CC^{\mathrm{per}}(M))$ and $\HN(M)=\Tot (CC^-(M))$ of the upper half plane and second quadrant bicomplexes, respectively. A homomorphism $\phi:M\to N$ of semi-cyclic complexes is a \emph{quasi-isomorphism} if it induces an isomorphism in Hochschild homology. This implies that it also induces an isomorphism in cyclic homology and its variants. 

\begin{ex}\label{ex:naivecyc}
As in \ref{rem:naive}, we remark that if $A$ is a ring with local units, then the complex $\HC(A) \defeq \HC(C^{\cyc}(A))$, computes the cyclic homology of $A$; the same holds for its negative and periodic variants.
\end{ex}

\begin{ex} \label{ex:cycobj-cyc}
Let $\cG$ be an ample groupoid and consider the simplicial weakly Boolean space $\cG_{\cyc}$ of Example \ref{ex:cyc}. Notice that we have a $\Z/(n+1)\Z$-action on $\cG_{\cyc}^n$ given by cyclic permutations,
\[
\tau_n \colon \cG_{\cyc}^n \to \cG_{\cyc}^n, \qquad \tau_n(g_0,g_1, \ldots, g_n) = (g_n,g_0,\ldots, g_{n-1}).
\]
Hence each module $\cc(\cG_{\cyc}^n)$ carries a $\Z/(n+1)\Z$ action given by $t_n = (-1)^{n}(\tau_n)_\ast$. These maps
are compatible with the simplicial structure and make $\H^{\cyc}(\cG)=\cc(\cG_{\cyc}^\acts)$
into a cyclic module.
\end{ex}

\begin{ex}\label{ex:cycobj-bar}
The modules $\{\B_n(\cG)\}_{n \ge 0}$ of Example \ref{ex:bar} together with the face and degeneracy maps defined therein assemble into a simplicial $\gring$-module $\B(\cG)$. Let 
\begin{gather*}
\tau_n \colon B_n(\cG) \to B_n(\cG)\\
\tau_n(g_0,\dots,g_n) =\left\{\begin{matrix} ((g_0,\dots,g_{n-1})^{-1},g_0,\dots,g_{n-2},g_{n-1} g_n) & n\ge 1\\
                                               g_0& n=0.\end{matrix} \right.
\end{gather*}
One checks that the simplicial module $\B(\cG)$ together with the maps $t_n = (-1)^n (\tau_n)_\ast$ is a cyclic module. 
\end{ex}

\begin{ex} \label{ex:cycobj-hg}
Let 
\begin{gather*}
\tau_n:\cG^{(n)}\to \cG^{(n)}    \\
\tau_n(g_1,\dots,g_n) =  ((g_1,\dots,g_n)^{-1},g_1,\dots,g_{n-1}).
\end{gather*}
The complex $\H(\cG)$, regarded as a simplicial module,
futher equipped with the maps $t_n = (-1)^n(\tau_n)_\ast$, is a cyclic module.
\end{ex}

\begin{rmk} \label{rmk:CC-res}
Let $M$ be a cyclic $k$-module and write $M_{-1} = \coker(b \colon M_1 \to M_0)$. By the argument of \cite{loday}*{2.5.7}, the complex $(M_\bullet, b')$ is always contractible. If we assume that $(M_\bullet,b)$ is (pure) exact in positive degrees, then we obtain a bicomplex with (pure) exact columns 
\[\begin{tikzcd}
	{} & {} & {} & {} \\
	{M_2} & {M_2} & {M_2} & {M_2} & {} \\
	{M_1} & {M_1} & {M_1} & {M_1} & {} \\
	{M_0} & {M_0} & {M_0} & {M_0} & {} \\
	{M_{-1}} & 0 & {M_{-1}} & 0 & {} \\
	0 & 0 & 0 & 0
	\arrow["b"', dotted, from=1-1, to=2-1]
	\arrow["{-b'}"', dotted, from=1-2, to=2-2]
	\arrow["b"', dotted, from=1-3, to=2-3]
	\arrow["{-b'}"', dotted, from=1-4, to=2-4]
	\arrow["b"', from=2-1, to=3-1]
	\arrow["{1-t}"', from=2-2, to=2-1]
	\arrow["{-b'}"', from=2-2, to=3-2]
	\arrow["N"', from=2-3, to=2-2]
	\arrow["b"', from=2-3, to=3-3]
	\arrow["{1-t}"', from=2-4, to=2-3]
	\arrow["{-b'}"', from=2-4, to=3-4]
	\arrow["N"', dotted, from=2-5, to=2-4]
	\arrow["b"', from=3-1, to=4-1]
	\arrow["{1-t}"', from=3-2, to=3-1]
	\arrow["{-b'}"', from=3-2, to=4-2]
	\arrow["N"', from=3-3, to=3-2]
	\arrow["b"', from=3-3, to=4-3]
	\arrow["{1-t}"', from=3-4, to=3-3]
	\arrow["{-b'}"', from=3-4, to=4-4]
	\arrow["N"', dotted, from=3-5, to=3-4]
	\arrow[from=4-1, to=5-1]
	\arrow["{1-t}"', from=4-2, to=4-1]
	\arrow[from=4-2, to=5-2]
	\arrow["N"', from=4-3, to=4-2]
	\arrow[from=4-3, to=5-3]
	\arrow["{1-t}"', from=4-4, to=4-3]
	\arrow[from=4-4, to=5-4]
	\arrow["N"', dotted, from=4-5, to=4-4]
	\arrow[from=5-1, to=6-1]
	\arrow[from=5-2, to=5-1]
	\arrow[from=5-2, to=6-2]
	\arrow[from=5-3, to=5-2]
	\arrow[from=5-3, to=6-3]
	\arrow[from=5-4, to=5-3]
	\arrow[from=5-4, to=6-4]
	\arrow[dotted, from=5-5, to=5-4]
\end{tikzcd}\]
\end{rmk}
\begin{rem}\label{rem:mixed}
Let $C$ be a cyclic complex and let $s:C\to C[1]$ the extra degeneracy, so that $1=sb'+b's$. Set $B:C\to C[1]$,
$B=(1-t)sN$. Then  $M(C)=(C,b,B)$ is what is called a \emph{mixed complex}; this means that $b^2=B^2=bB+Bb=0$. One can 
define the cyclic, periodic cyclic and negative cyclic bicomplexes of a mixed complex \cite{kasmix}. Their totalizations 
are the graded modules given in degree $n$ by $\HP(M)_n=\prod_{j\in\Z}M_{n+2j}\supset \HN(M)_n=\prod_{j\ge 0}M_{n+2j}$ 
and $\HC(M)=\bigoplus_{j\ge 0}M_{n-2j}$, with boundary maps induced by $b+B$. In the case of $M(C)$, the totalization of 
each of these is quasi-isomorphic to that of the corresponding complex defined above for $C$. An explicit formula for a 
quasi-isomorphism $\HC(M(C))\to \HC(C)$ is given in \cite{lq}*{Proposition 1.5}. The same formula works also for $\HN$ 
and $\HP$. If $M$ and $N$ are mixed complexes and we write $b$ and $B$ for their descending and ascending boundary maps, 
then an \emph{$S$-map} $G^\bu:M\to N$ is a sequence of homogeneous linear maps $G^n:M\to N[2n]$, $n\ge 0$, such that 
$[G^0,b]=0$ and such that $[G^{n+1},b]=-[G^n,B]$ for all $n\ge 0$. If $G^\bu$ is an $S$-map, then 
$G^{\infty}=\sum_{n\ge 0}G^n:\HP(M)\to\HP(N)$ is a chain map, which sends $\HN(M)\to \HN(N)$ and thus induces a chain map
$\HC(M)\to \HC(N)$. Each of these chain maps is a quasi-isomorphism whenever $G^0$ is one. 
\end{rem}

\section{Hochschild complexes for (twisted) Steinberg algebras} \label{sec:hh-stein}
\numberwithin{equation}{section}
In this section we set out to give a concrete description
of the Hochschild homology complex of a twisted Steinberg algebra. Throughout the section we fix an ample groupoid
$\cG$ with unit space $X$. Recall that we consider the ground ring $\gring$ as a discrete topological ring; we give the units $\cU(\gring)\subset \gring$ the subspace topology, which is also discrete. A \emph{(continuous) $2$-cocycle} on
$\cG$ over a commutative ring $k$ is a continuous map
\[
\omega \colon \cG^{(2)} \to \cU(k), 
\]
satisfying 
\begin{align}
   \omega(\alpha,\beta)\omega(\alpha\beta,\gamma) &= \omega(\alpha,\beta\gamma)\omega(\beta,\gamma);\label{eq:cocy1} \\
   \omega(r(\alpha),\alpha) &= \omega(\alpha,s(\alpha)) = 1.\label{eq:cocy2}
\end{align}

The \emph{twisted} Steinberg algebra (\cite{twisted-stein}*{page 5}) $\A(\cG, \omega)$ is the $k$-module $\cc(\cG)$ equipped with the product
given by:
\begin{equation}\label{eq:convotwistel}
(\eta \ast_\omega \mu)(\gamma) = \sum_{\gamma = \alpha\beta} \omega(\alpha, \beta)\eta(\alpha)\mu(\beta).
\end{equation}
We fix a $2$-cocyle $\omega$ in $\cG$. In the following lemma and elsewhere we adopt the following notation; for $A,B\subset\cG$, we put
\[
A\ttimes B:=(A\times B)\cap\cG^{(2)}.
\]

\begin{lem}\label{lem:convotwistchi}
Let $K,L\subset \cG$ be compact open bisections. Then there exist $n\ge 1$, $K_1,\dots,K_n\subset K$, $L_1,\dots,L_n\subset L$ compact open subspaces 
and $u_1,\dots,u_n\in\cU(k)$ such that $s(K_i)=r(L_i)$ for all $1\le i\le n$, $\omega$ is constantly equal to $u_i$ on $K_i\ttimes L_i$, and 
\begin{equation}\label{eq:convotwistchi}
\chi_K\ast_\omega\chi_L=\sum_{i=1}^nu_i\chi_{K_i\cdot L_i}.
\end{equation}
\end{lem}
\begin{proof}
Because $\omega$ is continuous, it is locally constant. Because the compact open bisections form a basis of the topology of $\cG$, the sets of the form $A\ttimes B$ with $A$ and $B$ compact open bisections form a basis of $\cG^{(2)}$, and are compact since $\cG^{(2)}\subset\cG\times\cG$ is closed because $\cG^{(0)}$ is Hausdorff. Hence the compact open subset $K\ttimes L$ admits a finite covering $\{K'_i\ttimes L'_i\colon 1\le i\le n\}$, 
such that for each $1\le i\le n$, $\omega$ is a constant $u_i\in\cU(k)$ on $K'_i\ttimes L'_i$. Since both intersections and differences of compact open subspaces of $K \ttimes L$ remain compact and open, and furthermore they can be written as disjoint unions of compact opens of the form $A\ttimes B$, we may assume that the covering is disjoint. Applying the formula \eqref{eq:convotwistel} we get the identity 
\[
\chi_K\ast_\omega\chi_L=\sum_{i=1}^nu_i\chi_{K'_i\cdot L'_i}.
\]
Next remark that for $K_i=s^{-1}(r(L'_i))\cap K'_i$ and $L_i=r^{-1}(s(K_i))\cap L'_i$, we have $s(K_i)=r(L_i)$ and $K_i\cdot L_i=K'_i\cdot L'_i$; in particular, the identity \eqref{eq:convotwistchi} holds. This finishes the proof.  
\end{proof}
\begin{coro}\label{coro:convotwistchi}
 The identity map of $\cC_k(\cG)$ defines an $\A(X)$-bimodule isomorphism $\cA_k(\cG)\iso \cA_k(\cG,\omega)$.
\end{coro}
\begin{proof}
In view of Lemma \ref{lem:convotwistchi}, it suffices to show that if $K\subset \cG$ is a compact open bisection, then $\omega$ is constantly equal to $1$ on $K\ttimes s(K)$ and $r(K)\ttimes K$. This is immediate from \eqref{eq:cocy2}.
\end{proof}

Fix a flat ring extension $\ell \subset \gring$.  Notice that the submodule $\cA_\ell(X) \subset \A(\cG, \omega)$ is
in fact a commutative $\ell$-subalgebra; Proposition \ref{prop:twisted-hh-rel} below shows that the Hochschild homology of $\A(\cG)$ over $\ell$ can be instead computed over $\cA_\ell(X)$. We need the following lemma.

\begin{lem}\label{lem:cxortho}
Let $\mathcal{F}$ be the set of all finite sets of orthogonal idempotents of $\cA_\ell(X)$ and let $n\ge 1$. Then for every $a_1,\dots, a_n\in \cA_\ell(X)$ there exists $F\in\mathcal{F}$ such that $a_1,\dots,a_n\in kF$.
\end{lem}
\begin{proof}
It suffices to show that the condition of the lemma holds when each $a_i$ is the characteristic function of some compact open subset of $X$, since the latter span $\cA_\ell(X)$. Let 
$A_1,\dots, A_n\subset X$ be compact open. For each subset $I\subset [n]_+=\{1,\dots,n\}$ let $I^c=[n]_+\setminus I$, $A_I=\bigcap_{i\in I}A_i\setminus\bigcup_{j\in I^c}A_j$. Because $X$ is Hausdorff, each subset 
$A_I$ is compact open, so $\chi_{A_I}\in\A_\ell(X)$. Moreover we have $A_I\cap A_J=0$ for $I\ne J$ and for all $i\in [n]_+$, $A_i=\bigsqcup_{i\in I}A_I$. Thus 
$F=\{\chi_{A_I}\colon I\subset[n]_+\}\in\cF$ and we have $\chi_{A_i}=\sum_{i\in I}\chi_{A_I}\in kF$.
\end{proof}

\begin{prop} \label{prop:twisted-hh-rel}
The projection map is a quasi-isomorphism $\Cy(\A(\cG,\omega)/\ell) \weq \Cy(\A(\cG, \omega) / \cA_\ell(X))$.
\end{prop}
\begin{proof}
Lemma \ref{lem:cxortho} implies that $A = \cA_\ell(X)$ satisfies part i) of Lemma \ref{lem:hhsep}. Moreover the elements $\chi_K$ with $K\subset X$ compact open form a system of local units for $B=\A(\cG, \omega)$, so part ii) of the latter lemma also holds. Hence the corollary follows from Lemma \ref{lem:hhsep}.
\end{proof}

The proposition above describes the Hochschild
homology complex of $\A(\cG,\omega)$ over $A_\ell(X)$ in terms of the nerve of $\cG$.

\begin{prop} \label{prop:hh=nerve-com}
For each $n \ge 0$ there are $\ell$-module isomorphisms
\[
\HH(\A(\cG,\omega)/\cA_\ell(X))_n \cong \cA_\ell(X) \otimes_{\cA_\ell(X) \otimes_\ell \cA_\ell(X)^{\op}} \cC_c(\cG^{(n+1)}, \gring^{\otimes_\ell n+1}).
\]
\end{prop}
\begin{proof}
By definition, $\A(\cG,\omega)$, like $\A(\cG)$, has the same underlying $k$-module $\cC_c(\cG,k)$. Moreover by Corollary \ref{coro:convotwistchi}, these two algebras have the same underlying $\A(X)$-bimodule structure, and hence also the same $\cA_\ell(X)$-bimodule structure. From this we see that $\HH(\A(\cG, \omega) / \cA_\ell(X))_n=\HH(\A(\cG) / \cA_\ell(X))_n$ as $k$-modules. In the identities below we write $M_\# = M/[M,\cA_{\ell}(X)] = \cA_\ell(X) \otimes_{\cA_\ell(X) \otimes_\ell \cA_\ell(X)^{\op}} M$. Using
Corollary \ref{coro:agzbimo} and the fact that $\cA_k(\cG)=k\otimesl \cA_{\ell}(\cG)$ we then obtain 
\begin{align*}
\HH(\A(\cG, \omega) / \cA_\ell(X))_n &= \left((k\otimes_\ell \cA_{\ell}(\cG))^{\otimes_{\cA_\ell(X)}n+1}\right)_\#\\
                                    &= k^{\otimesl n+1}\otimesl (\cA_{\ell}(\cG)^{\otimes_{\cA_\ell(X)}n+1})_\#\\
                                    &= k^{\otimesl n+1}\otimesl\cC_c(\cG^{(n+1)},\ell)_{\#}\\
                                    &=\left(k^{\otimesl n+1}\otimesl\cC_c(\cG^{(n+1)},\ell)\right)_\#\\
                                    &=\cC_c(\cG^{n+1},k^{\otimesl n+1})_\#.
\end{align*}    
\end{proof}

\begin{rmk}\label{rmk:hh-gens}
By Proposition
\ref{prop:hh=nerve-com} we get that $\HH(\A(\cG, \omega) / \cA_\ell(X))_n$ is $k$-linearly spanned by elements of 
the form $a_0\otimes\cdots\otimes a_n\otimes\chi_{[K_0|\cdots|K_n]}$ with $K_i$ a compact open bisection such that $s(K_i)=r(K_{i+1})$ for all $0\le i\le n-1$ and such that $\omega$ is a constant $\omega_i\in\cU(k)$ on $K_i\ttimes K_{i+1}$ for all 
$0\le i\le n$ (where $n+1$ is taken modulo $n$).  
For these elements, we have
\begin{gather}
b(a_0\otimes a_n\otimes\chi_{[K_0|\cdots|K_n]})=\sum_{i=0}^{n-1}
(-1)^ia_0\otimes\cdots\otimes \omega_ia_ia_{i+1} \otimes \cdots a_n\otimes 
\chi_{[K_0|\cdots|K_iK_{i+1}|\cdots|K_n]}\label{eq:btwist}\\
+(-1)^n\omega_n a_na_0\otimes\cdots\otimes a_{n-1}\chi_{[K_nK_0|\cdots|
K_{n-1}]}\nonumber
\end{gather}
\end{rmk}

This description of the Hochschild complex allows us
to see that if $\cG$ is Hausdorff then (twisted)
groupoid homology is a direct summand of $\HH(\A(\cG,\omega)/\cA_\ell(X))$, as we proceed to show.

\begin{lem}\label{lem:twistsplit}
Let $n \ge 0$ and
\[
\Gamma(\cG,X)_n =\{(g_0,\cdots,g_n)\colon g_0\cdots g_n\in X\} \subset \cG^{(n+1)}.
\]
The open inclusions $\Gamma(\cG,X)_\bullet \subset \cG^{(n+1)}$ 
define a sub-semicycic module

\goodbreak

\noindent $\cc(\Gamma(\cG,X)_\bullet, k^{\otimesl n+1})$ of $\Cy(\A(\cG, \omega)/\cA_\ell(X))$.
If $\cG$ is Hausdorff, then the inclusion
\[
\inc \colon \cc(\Gamma(\cG,X)_\bullet, k^{\otimesl n+1}) \to \Cy(\A(\cG, \omega)/\cA_\ell(X)). 
\]
is a split monomorphism of semicyclic modules, whose
left inverse is given by the restriction-induced homomorphism 
\[
\res \colon \Cy(\cA_k(\cG, \omega)/\cA_\ell(X))\to \cC_c(\Gamma(\cG,X)_\bu, k^{\otimesl \bu+1}).
\]
In particular, for every Hausdorff ample groupoid we have that $\cC_c(\Gamma(\cG,X)_\bu, k^{\otimesl \bu+1})$ is a direct summand of $\Cy(\cA_k(\cG,\omega)/\cA_\ell(X))$. 
\end{lem}
\begin{proof} 
The existence of $\inc$ is straightfoward from the functoriality of $\cc$ with respect to étale maps. When $\cG$ is Hausdorff, it follows that $X\subset \cG$ is clopen and thus $\Gamma(\cG,X)_n\subset \cG^{(n+1)}$ is clopen. 
Hence $\cC_c(\Gamma(\cG,X)_n,k^{\otimesl n+1})$ is a direct 
summand of $\cC_c(\cG^{(n+1)},k^{\otimesl n+1})$, with the inclusion split by the restriction map. In light of Remark \ref{rmk:hh-gens},
the module $C_c(\Gamma(\cG,X)_n,k^{\otimesl n+1})$ is spanned by the elements of the form
\begin{gather*}
a_0\otimes\cdots\otimes a_n\otimes\chi_{[K_0|\cdots|K_n]\cap \Gamma(\cG,X)_n}
\end{gather*}
where $a_0,\dots,a_n \in\gring$ and $K_0,\dots, K_n$ are as in Remark \ref{rmk:hh-gens}. Consider the compact open set $K'=(K_0\cdots K_n)\cap X$. Let 
$K'_n=s^{-1}(K)\cap K_n$; for $1\le i\le n$, set $K'_{i-1}=K_{i-1}\cap s^{-1}
(r(K'_i))$. Then each $K'_i$ is a compact open bisection, $s(K'_i)=r(K'_{i+1})$ for 
all $0\le i\le n-1$, $\omega$ is constant on $K'_i\ttimes K'_{i+1}$ for all $0\le 
i\le n$, and $\Gamma(\cG,X)_n \cap [K_0|\cdots|K_n]=[K'_0|\cdots|K'_n]$. Hence 
$C_c(\Gamma(\cG,X)_n,k^{\otimesl n+1})$ is spanned by the elements $a_0\otimes 
a_n\otimes\chi_{[K_0|\cdots|K_n]}$ with $[K_0|\cdots |K_n]\subset \Gamma(\cG,X)_n$. One checks 
that if $[K_0|\cdots |K_n]\subset \Gamma(\cG,X)_n$ then $[K_n|K_0|\cdots|K_{n-1}]$ and $[K_nK_0|\cdots|K_{n-1}]\subset 
\Gamma(\cG,X)_{n-1}$. By \eqref{eq:btwist}, this goes to show that
$\cC_c(\Gamma(\cG,X)_n,\gring^{\otimes_\ell \bu+1})$ is a sub-semicyclic module of $\Cy(\A(\cG) / \cA_\ell(X))$. Moreover one checks that 
the restriction-induced map 
\[
\res\colon\Cy(\A(\cG) / \cA_\ell(X))\to \cC_c(\Gamma(\cG,X)_\bu,\gring^{\otimes_\ell \bu+1})
\]
is a homomorphism of semicyclic modules. 
\end{proof}
    
\begin{defi}[Twisted groupoid homology]\label{defi:twistedhom}
For each $0 \le j \le n$, write 
\[
\sigma_j^n \colon k \to k^{\otimes_\ell n+1}, \qquad \lambda \mapsto 1 \otimes \cdots \otimes \overbrace{\lambda}^i \otimes \cdots \otimes 1.
\]
for the degeneracy maps and $\delta^n_j \colon k^{\otimes_\ell n+1} \to k^{\otimes_\ell n}$ for the face maps of $\HH_\ast(k/\ell)$.
We define the \emph{twisted groupoid homology complex} as $\H(\cG^\omega, k/\ell) = \cC_c(\cG^{(\bu)}, \gring^{\otimes_\ell \bu+1})$ with boundary maps $\partial = \sum_{i=0}^n (-1)^i d_i$ where
\begin{align*}
d_i(f)(g_1,\dots,g_n) &= \sum_{g_i = \alpha\beta} \sigma^n_i(\omega(\alpha,\beta))\delta^n_i(f(g_0,\ldots, g_{i-1},\alpha,\beta,g_{i+1},\dots,g_n)), \, 0 < i < n, \\
d_{0}(f)(g_1,\dots,g_n) &= \sum_{s(\beta) = r(g_1)} \sigma^n_0(\omega((g_1\cdots g_n)^{-1}\beta^{-1},\beta))\delta^n_{0}(f(\beta, g_1,\dots,g_{n})),\\
d_n(f)(g_1,\dots,g_n) &= \sum_{r(\beta) = s(g_n)} \sigma^n_0(\omega(\alpha, \alpha^{-1} (g_1\cdots g_n)^{-1}))\delta^n_{n}(f(g_1,\dots,g_{n},\alpha)).
\end{align*}

We call the homology $H_\ast(\cG^\omega,k/\ell)$ of $\H(\cG^\omega, k/\ell)$ the \emph{twisted groupoid homology} of $\cG$ with respect to $\omega$ and the ring extension $\ell \subset k$.  
\end{defi}
\begin{rem}\label{rem:htwistiscyc}
The complex $\H(\cG^{\omega},k/\ell)$ has a semi-cyclic $\ell$-module structure, with 
\[
t_n(f)(g_1,\cdots,g_n)=t_n(f((g_1\cdots g_n)^{-1},g_1,\cdots, g_{n-1}),
\]
where the second $t_n$ is the signed cyclic permutation coming from the cyclic module $k^{\otimesl \bu+1}$.
\end{rem}
\begin{rmk}\label{rem:trivtwist}
When the $2$-cocycle $\omega$ is trivial, we obtain the semicyclic module $\H(\cG, k/\ell) = \H(\cG) \boxtimes \HH(k/\ell)$ arising from the tensor product of the standard semi-cyclic modules. Since $k$ is flat over $\ell$, by the Eilenberg-Zilber theorem and Künneth's formula, we have
that
\[
H_n(\cG, k/\ell) = 
\bigoplus_{i+j=n} H_i(\cG) \otimes_\ell HH_j(k/l).
\]
\end{rmk}

\begin{prop} \label{prop:twisted=twisted}
The isomorphisms $\cG$-spaces
\[
\cG^{(n)} \to \Gamma(\cG,X)_n, \qquad (g_1,\dots,g_n) \mapsto ((g_1\cdots g_n)^{-1}, g_1,\dots, g_n).
\]
induce an isomorphism of semicyclic modules $\H(\cG^\omega, k/\ell) \cong \H(\cC_c(\Gamma(\cG,X)_\bu, \gring^{\otimes_\ell \bu+1})$.
\qed
\end{prop}

\begin{coro} \label{coro:twisted-retract}
The twisted group homology complex $\H(\cG^\omega, k/\ell)$ can be identified with a subcomplex of $\HH(\A(\cG, \omega)/\cA_\ell(X))$. If $\cG$ is Hausdorff, then it is moreover a direct
summand.
\qed
\end{coro}

\section{Cyclic groupoid homology}\label{sec:hc}

In this section we discuss, for an ample groupoid $\cG$, the cyclic homology of the cyclic module $\H(\cG)$ of Example \ref{ex:cycobj-hg}. We write $\HC(\cG)$, $\HN(\cG)$ and $\HP(\cG)$ for the cyclic, negative cyclic and periodic cyclic complexes of $\H(\cG)$. 
If $\cG$ is Hausdorff, then by Corollary \ref{coro:twisted-retract} each of these complexes is a direct summand of the cyclic, negative cyclic and periodic cyclic homology of
$\A(\cG)$ respectively. 

In what follows 
we shall use some tools and terminology from relative homological algebra. An \emph{extension} of $\A(\cG)$-modules
is a kernel-cokernel pair
\[
K \stackrel{i}{\rightarrowtail} E \stackrel{p}{\twoheadrightarrow} Q.
\]
We say that it is \emph{semi-split} if $p$ has an $\cA(X)$-linear section. An $\A(\cG)$-module $P$ is 
\emph{relatively projective} if $\hom_{\A(\cG)}(P,-)$ maps semi-split extensions to exact sequences, and \emph{relatively free} if $P \cong \A(\cG) \otimes_{\cA(X)} N$ for some $\cA(X)$-module $N$. A \emph{(relatively) projective resolution} of an $\A(\cG)$-module $M$ is an exact complex of $\A(\cG)$-modules
\[
\cdots \to P_2 \to P_1 \to P_0 \to M \to 0
\]
that admits an $\A(X)$-linear contracting homotopy, and in which
each module $P_i$ is relatively projective. As in the setting of classical homological algebra, we recall that relatively free modules are relatively free. We shall also use that maps between $\A(\cG)$-modules extend to chain maps between projective
resolutions, and that two such extensions are unique up to 
an $\A(\cG)$-linear chain homotopy.

\begin{lem}\label{lem:rel-proj}
 Let $\cG$ be an ample groupoid with unit space $X$ and let $n \ge 1$. The unital $\A(\cG)$-module $\cc(\cG^{(n)})$ is relatively free with respect to $\A(X)$; in particular, it is relatively projective.
\end{lem}
\begin{proof} We have $\cc(\cG^{(n)}) \cong 
\cc(\cG^{(n-1)}) \otimes_{\A(X)} \A(\cG)$.
\end{proof}

\begin{thm}\label{thm:HC-main}
  Let \(\cG\) be an ample groupoid. We have quasi-isomorphisms
  $$
\HC(\cG)\overset{\sim}{\to} \bigoplus_{n\ge 0}\H(\cG)[-2n],\,
 \HN(\cG)\overset{\sim}{\to}\prod_{n\ge 0}\H(\cG)[2n],\,\HP(\cG)\overset{\sim}{\to} \prod_{n\in\Z}\H(\cG)[2n].
 $$

  of complexes of \(k\)-modules. Consequently, we obtain isomorphisms 
  \[
  HC_*(\cG) \cong \bigoplus_{i \geq 0} H_{* - 2i}(\cG),\, HN_*(\cG) \cong \prod_{i \geq 0} H_{* + 2i}(\cG), \, HP_*(\cG) \cong \prod_{i \in \Z} H_{* + 2i}(\cG).
  \] 
  for all \(*\).  
\end{thm}
\begin{proof} Let $\bB(\cG)$ be the cyclic module of Example \ref{ex:cycobj-bar}. Observe that $(\bB(\cG),b)$ is a resolution of $\bB(\cG)_{-1}=\cc(\cG^{(0)})$ by relatively free $\A(\cG)$-modules. Hence for $n>0$, any chain map $\bB(\cG)\to\bB(\cG)[n]$ is $\A(\cG)$-linearly chain homotopic to zero, since it lifts the zero map $\cc(\cG^{(0)})\to 0$. As in Remark \ref{rem:mixed}, we consider the associated mixed complex $M=(\bB(\cG),b,B)$. Set $N=(\bB(\cG),b,0)$ and define an $S$-map $G^\bu:M\to N$ recursively as follows. Set $G^0=\id_{\bB(\cG)}$; as remarked above $B$ is homotopic to zero, so there is an $\cA(\cG)$-linear maps $G^1$ so that $[G^1,b]=-B=-G^0B$. Let $n\ge 1$ and assume $G^n$ defined so that $[G^n,b]=-G^{n-1}B$. Then $G^nBb=-G^nbB=(-bG^n+G^{n-1}B)B=-bG^nB$, so $G^nB$ is an $\cA(\cG)$-linear chain map 
$\bB(\cG)\to\bB(\cG)[2n+1]$, and is therefore homotopic to zero. Hence we can find $G^{n+1}:\bB(\cG)\to\bB(\cG)[2(n+1)]$
with $[G^{n+1},b]=-G^nB$. Then $$\hat{G}^\bu=\cc(\cG^{0})\otimes_{\A(\cG)}G^\bu:\bar{M}=(\H(\cG),b,B)\to \bar{N}=(\H(\cG),b,0)$$ is an $S$-map with $\bar{G}^0=\id_{\H(\cG)}$ and therefore induces quasi-isomorphisms at the level of $\HC$, $\HP$ and $\HN$.
\end{proof}

\section{Cyclic nerve computations}\label{sec:first}

This section is dedicated to studying the cyclic nerve complex $\H^{\cyc}(\cG)$. Regard $\cG^{\Iso}$ as a $\cG$-space with the conjugation. We show in Corollary \ref{coro:hadjG} that $\H^{\cyc}(\cG)\cong \H(\cG,\cG^{\Iso})$. In particular, $H_*^{\cyc}(\cG)=H_*(\cG,\cC_c(\cG^{\Iso}))$, which, when $G$ is a group, equals $HH_*(k[G])$. More generally, if $M$ is a $k[G]$-bimodule and $M^{\adj}$ is the $k$-module $M$ equipped 
with the left $G$-action $g\cdot m=gmg^{-1}$, then 
\begin{equation*}
HH_*(k[G],M)=H_*(G,M^{\adj}).
\end{equation*} 
This follows from the fact that the functor $M\mapsto M^{\adj}$ is exact and preserves projectives, and that $M_\sharp=k\otimes_{k[G]}M^{\adj}$. Remark that the definition of the adjoint action on a $k[G]$-bimodule uses the homomorphism $k[G]\to k[G]\otimes k[G]^{\op}$, $g\mapsto g\otimes g^{-1}$. There is no analogue of this algebra map for general ample groupoids. 
We show however that for any ample groupoid $\cG$, the standard semicyclic module for $\cA(\cG)$ maps surjectively onto $\H^{\cyc}(\cG)$ (Proposition \ref{prop:agzcyc}) and that if $\cG$ is Hausdorff both contain the semicyclic module $H(\cG)$ as a direct summand (Corollary \ref{coro:wsplit}). Furthermore, under additional hypothesis, we obtain a decomposition of the latter complex (see Theorems \ref{thm:H-cyc} and \ref{thm:sparse}) which recovers Burghelea's decomposition for Hochschild and cyclic homology in the group case \cite{burghelea}. 

\comment{
We have a groupoid analogue of the isomorphism above replacing Hochschild homology by the cyclic nerve complex, where the isotropy subspace plays the role of the adjoint representation. Before describing this result in more detail, we record the existence a comparison map between these two complexes.
}
\subsection{Mapping \topdf{$\HH$}{HH} to \topdf{$\H^{\cyc}$}{Hcyc}}
\numberwithin{equation}{subsection}
\begin{lem}\label{lem:Xtensor-n=cyc}
Let $\cG$ be an ample groupoid with unit space $X$ 
and $n \ge 0$. There is an $\cA(X)$-bimodule epimorphism
\[
\A(X) \otimes_{\A(X) \otimes_k \A(X)^{\op}} \cc(\cG^{(n+1)}) \onto \cc(\cG_{\cyc}^n), \, \chi_U \otimes \chi_{[A_0 | \cdots | A_n]}\mapsto \chi_{(UA_0| \cdots | A_n U)}.
\]
\end{lem}
\begin{proof} Since $X$ is Hausdorff, we know that $\cG^n_{\cyc} = \{(g_0,\dots,g_n) \in \cG^{(n+1)} : r(g_0) = s(g_n)\}$ is a closed subspace 
of $\cG^{(n+1)}$. Hence we have a restriction morphism $\res \colon \cc(\cG^{(n+1)}) \to \cc(\cG_{\cyc}^n)$ mapping $\chi_{[A_0 | \cdots | A_n]}$ to $\chi_{(A_0| \cdots | A_n)}$.
It suffices to see that this descends to the commutator quotient $\cc(\cG^{(n+1)})_\sharp$. This is equivalent to showing that $\res(\chi_U\cdot  \chi_{[A_0 | \cdots | A_n]}) = \res(\chi_{[A_0 | \cdots | A_n]} \cdot \chi_U )$
for every compact open subset $U \subset X$ and compact open bisections $A_0,\dots, A_n \subset \cG$, which follows from Lemma \ref{lem:slice-tuples} iii).
\end{proof}

\begin{prop}\label{prop:agzcyc}
If $\cG$ is a Hausdorff ample groupoid, there is a levelwise surjective 
map of semicylic modules
\begin{align}\label{mor:mu}
&\mu\colon C^{\cyc}(\A(\cG)/\A(X))\to \H^{\cyc}(\cG) = \cc(\cG^{\acts}_{\cyc}),\\     
&\mu(\phi_0\otimes\cdots\otimes\phi_n)(g_0,\dots,g_n) :=\phi_0(g_0)\cdots\phi_n(g_n)\notag
\end{align}
Further, if $\cG$ is $\Lambda$-graded, 
then $\mu$ is a homogeneous map of degree zero between $\Lambda\times \N_0$-graded $k$-modules.
\end{prop}
\begin{proof} By Proposition \ref{prop:hh=nerve-com}, we have isomorphisms $\HH(\A(\cG)/\A(X))_n \cong \cc(\cG^{(n+1)})_\sharp$.
Applying Lemma \ref{lem:Xtensor-n=cyc} we obtain the 
desired map. It is a morphism of semicyclic modules because it is induced by the inclusion of semicyclic spaces $\cG^{\cyc}\subset \cG^{\bu}$.
\end{proof}

\begin{rmk}\label{rem:burgh} We point out that the map of Proposition \ref{prop:agzcyc}
is an isomorphism when $\cG$ is discrete. In particular, 
our description of $\H^{\cyc}(\cG)$ on Theorem \ref{thm:H-cyc}
below recovers
Burghelea's theorem for Hochschild 
and cyclic homology of group algebras.
\end{rmk}

\subsection{Invariant subspaces of \topdf{$\cG^{\Iso}$}{Iso(G)}
and direct summands of \topdf{$\H(\cG)$}{H(G)}}\label{subsec:preburghe}

Fix an ample groupoid $\cG$.
We say that $W \subset \cG^{\Iso}$ is \emph{invariant}
if 
\[
\cG \acts W = \{gwg^{-1} : s(g) = r(w), w \in W\} \subset W.
\]

Such a subspace defines a cyclic subobject of $\cG_{\cyc}^\acts$; namely, 
\[
\Gamma(\cG,W)_n = \{(g_0,\ldots, g_n) \in \cG^n_{\cyc} : g_0 \cdots g_n \in W\}.
\]
Notice also that each space $\Gamma(\cG,W)_n$ is open (resp. closed) whenever $W$ is open (resp. closed), since $\Gamma(\cG, W)_n$ is the preimage of $W$ under the product map $\cG^n_{\cyc} \to \cG^{\Iso}$. 
If $\cG$ is $\Lambda$-graded, the restriction of the degree map makes $W$ into a $\Lambda$-graded 
$\cG$-space.

\begin{lem} \label{lem:gamma=homol}
The assignment 
\[
\cG^{(n)} \times_{\cG^{(0)}} W\to\Gamma(\cG,W)_n, \, ((g_1,\ldots, g_n),w) \mapsto (w(g_1\ldots g_n)^{-1},g_1,\ldots,g_n)
\]
is a homeomorphism with inverse $(g_0,\ldots,g_n) \mapsto ((g_1,\ldots,g_n),g_0g_1\cdots g_n)$.
If we equip $W$ with the left $\cG$-space structure given by conjugation, then 
the above map defines an isomorphsm of simplicial spaces between $\Gamma(\cG,W)$
and the simplicial space $\cG^{(\acts)} {}_s \times_r W$ associated to the groupoid homology of $\cG$ with coefficients in $W$.
The cyclic structure on $\Gamma(\cG,W)$
corresponds on the left hand side to that given by
\[
t ((g_1,\ldots, g_n),w) =
((w(g_1\cdots,g_n)^{-1},g_1,\ldots,g_{n-1}),g_n w g_n^{-1}).
\]
In particular, we have an isomorphism of cyclic modules

\[
\H(\cc(\Gamma(\cG,W))) \cong \H(\cG,W).
\]
\end{lem}
\begin{proof} Straightforward.
\end{proof}

\begin{coro}\label{coro:hadjG} For every ample groupoid $\cG$
we have $\H^{\cyc}(\cG) \cong \H(\cG, \cG^{\Iso})$.
\qed
\end{coro}

\begin{rmk} \label{rmk:gamma=homol-gr}
If $\cG$ is $\Lambda$-graded, and we equip $\cG^{(n)}$
with the trivial grading, and $W$ with its canonical grading as a subspace of $\cG$, 
then the homeomorphism of Lemma \ref{lem:gamma=homol} is compatible 
with the grading of $\Gamma(\cG, W)$ induced by the one on $\cG^n_{\cyc}$. 
\end{rmk}

Recall that a groupoid is called \emph{principal} if $\cG^{\Iso} = \cG^{(0)}$. 

\begin{prop}\label{prop:ppal}
If $\cG$ is a principal groupoid, then $\H^{\cyc}(\cG) \cong \H(\cG)$ as cyclic modules.
\end{prop}
\begin{proof} Because $\cG$ is principal, $\cG_{\cyc}^\acts=\Gamma(\cG,\cG^{(0)})$. Now use Lemma \ref{lem:gamma=homol}.
\end{proof}

\begin{lem}\label{lem:wsplit} 
If $W$ is an open (resp. clopen) subspace of $\cG^{\Iso}$, then $\H(\cG,W)$ is a subcomplex (resp. direct summand) of $\H^{\cyc}(\cG)$. Furthermore, if $W$ is clopen in $\cG$, then the inclusion $\Gamma(\cG, W) \subset \cG^{(n+1)}$ induces a split-monomorphism $\H(\cG,W) \to \HH(\A(\cG)/\A(X))$. 
\end{lem}
\begin{proof} The proof of the first assertion is immediate from Lemma \ref{lem:gamma=homol} and the fact that if $Z$ is a clopen subspace of a space $Y$, then $\cc(Y) \cong \cc(Z) \oplus \cc(Y\setminus Z)$. Since  $\Gamma(\cG,W)$ can be also thought of as the preimage of $W$ under the product map $\cG^{(n+1)} \to \cG$, it is a clopen subpace of $\cG^{(n+1)}$. Hence there is an (open) inclusion induced map $\cc(\Gamma(\cG,W)) \to \cc(\cG^{(n+1)})$, which can be composed with the projection onto $HH(\A(\cG)/ \A(X))_n = \cc(\cG^{(n+1)})_\sharp$. Its right inverse is given by \eqref{mor:mu} followed by the direct summand projection $\H^{\cyc}(\cG) \to \cc(\Gamma(\cG, W))$. 
\end{proof}

\begin{coro}[cf. Corollary \ref{coro:twisted-retract}]\label{coro:wsplit} For every ample groupoid $\cG$ the semicylic module $\H(\cG)$ is a sub-semicylic module of $\H^{\cyc}(\cG)$. If $\cG$ is Hausdorff, then $\H(\cG)$ is a direct summand of both $\H^{\cyc}(\cG)$ and $\HH(\A(\cG)/\A(X))$, 
and the following diagram commutes
\[
\begin{tikzcd}
\HH(\A(\cG)/\A(X)) \arrow{r}{\eqref{mor:mu}}&  \H^{\cyc}(\cG) \\ \H(\cG) \arrow[hook]{u} \arrow[hook]{ur} &
\end{tikzcd}
\]
\qed
\end{coro}

\subsection{Homology with coefficients on discrete orbits of \topdf{$\cG^{\Iso}$}{Iso(G)}}\label{subsec:burghe}

Let $\eta \in \cG^{\Iso}$ and assume that $\cG \acts \eta$
is discrete. Put $s(\eta) = r(\eta) = x$ and write 
$\cH := (\cG^x_x)_\eta$ for the centralizer subgroup of $\eta$. 

Notice that since $s \colon \cG \to \cG^{(0)}$ is an étale map, 
the fiber $\cG^x$ over $x$ is discrete, 
and so are any subspace such as $\cG^x_x$ and all of its centralizer subgroups.
In particular 
$s$ makes $\cG^x$ into an étale $\cH$-space.

\begin{lem} \label{lem:flat-fiber-right}
$\cc(\cG^x)$ is flat as a right $\cH$-module.
\end{lem}
\begin{proof} By Proposition \ref{prop:basic-flat}, and the fact that 
$\cG^x$ is an étale $\cH$-space, it suffices 
to show that the action $\cG^x \curvearrowleft \cH$ 
is basic; that is, it suffices to see that the map
\[
\cG^x \times \cH \to \cG^x \times_{\cG^x/\cH} \cG^x, \,
(\alpha,h) \mapsto (\alpha,\alpha h)
\]
is a homeomorphism. Since this map is a bijection between discrete 
spaces, the conclusion follows.
\end{proof}

\begin{lem}\label{lem:orbit-quotient}
There is a homeomorphism 
\[
\cG^x/\cH \to \cG \acts \eta, \quad 
[g] \mapsto g\eta g^{-1}.
\]
\end{lem}
\begin{proof} Both spaces are discrete and the map above is a bijection.
\end{proof}

\begin{prop} \label{prop:stab-hom}
$H_\ast(\cG,\cG^x/\cH) \cong H_\ast(\cH)$.
\end{prop}
\begin{proof} We adapt the proof of Shapiro's Lemma \cite{miller-corre}*{Lemma 2.19} to the present setting. 
We consider the canonical flat resolution of $\cc(\cH^\bullet) \to \cc(\cH^{(0)}) = \gring$ as an $\A(\cH)$-module,
dually to Example \ref{ex:bar}. 
By Lemma \ref{lem:flat-fiber-right} we have that
$P_\bullet = \cc(\cG^x) \otimes_{\A(\cH)} \cc(\cH^\bullet)$ is a flat resolution of $\cc(\cG^x) \otimes_{\A(\cH)} \cc(\cH^{(0)})$. 
By Proposition \ref{prop:fiber-tensor}, the latter is $\cc(\cG^x) \otimes_{\A(\cH)} \cc(\cH^{(0)}) 
\cong \cc(\cG^x \times_\cH \cH^{(0)}) = \cc(\cG^x/\cH)$.
Hence we may compute $H_\bullet(\cG,\cG^x/\cH)$
as the homology of the complex $\cc(\cG^{(0)}) \otimes_{\A(\cG)} P_\bullet$. Since 
$\cc(\cG^{(0)}) \otimes_{\A(\cG)} \cc(\cG^x)
 = \cc(\cG^{(0)} \times_{\cG} \cG^x) = \cc(\cH^{(0)})$, 
it follows that $\cc(\cG^{(0)}) \otimes_{\A(\cG)} P_\bullet\cong \cc(\cH^{(0)})\otimes_{\A(\cH)}\cc(\cH^{\bu})$ 
which computes $H_*(\cH)$.
\end{proof}

\begin{thm} \label{thm:H-cyc}
Let $\cG$ be an ample, Hausdorff groupoid. Set $X=\cG^{(0)}$. Assume that $\cG^{\Iso} \setminus X$ is discrete. Choose $\cR\subset X$ such that each element of $\cR$ has nontrivial isotropy and such that each element of $X$ with nontrivial isotropy is isomorphic in $\cG$ to exactly one element of $\cR$. For each $x\in\cR$, choose a set $Z_x$ of representatives of the non-trivial conjugacy classes of $\cG^{x}_{x}$. We have a quasi-isomorphism of cyclic modules
\[
\H(\cG)\oplus \bigoplus_{x \in \cR} \bigoplus_{\eta \in Z_x} \H((\cG^{x}_{x})_\eta)\weq\H^{\cyc}(\cG).
\]

Further, if $\cG$ is $\Lambda$-graded, then under the quasi-isomorphism above, the homogeneous component of degree $m$ of
$\H^{\cyc}(\cG)$
corresponds to 
\[
\bigoplus_{x \in \cR}\, \bigoplus_{\substack{\eta \in Z_x,\\ |\eta| = m}} \H((\cG^{x}_{x})_\eta)
\]
if $m\in\Lambda\setminus\{0\}$ and to 
\[
\H(\cG) \oplus\, \bigoplus_{x\in \cR} \bigoplus_{\substack{\eta \in Z_x,\\ |\eta| = 0}} \H((\cG^{x}_{x})_\eta).
\]
if $m=0$. 
\end{thm}
\begin{proof} We have a decomposition into clopen invariant sets of the form 
\[
\cG^{\Iso} = X \sqcup \bigsqcup_{x\in \cR} \bigsqcup_{\eta \in Z_x} \cG \acts \eta.
\]
Hence
\[
\H^{\cyc}(\cG) \cong \H(\cG,X) \oplus \bigoplus_{x \in \cR} \bigoplus_{\eta \in Z_x} \H(\cG,\cG \acts \eta).
\]
Now apply Lemma  \ref{lem:orbit-quotient} and Proposition \ref{prop:stab-hom}.
\end{proof}

\subsection{Semigroup actions with sparse fixed points}\label{subsec:sparse}

We now give a reformulation of Theorem \ref{thm:H-cyc} for the 
groupoid of germs of a semigroup action.
Let $\cS$ be an inverse semigroup, that is, a semigroup such that for every element $s\in\cS$ there is a unique element $s^*$ which is \emph{inverse} to $s$, in the sense that $ss^\ast s = s^\ast$ and $s^\ast s s^\ast = s$. The subset $\cS\supset \mathcal{E}(\cS)$ of its idempotent elements forms a commutative subsemigroup \cite{pater}*{Proposition 2.1.1}.
Let $X$ be a locally compact Hausdorff space. The set
\[
\cI(X) = \{f \colon U \to V : U,V \subset X \text{ open subsets and $f$ a homeomorphism} \}.
\]
is an inverse semigroup with the operations of partial inverses and partial composition. An action $\cS \curvearrowright X$ is a semigroup homomorphism $\phi  \colon \cS \to \cI(X)$. We write 
$\Dom(s)$ for the domain of $\phi(s)$ and $s \cdot x = \phi(s)(x)$. The \emph{orbit} of $x\in X$ is 
\[
\Or(x) = \{s \cdot x : s \in \cS, \Dom(s) \ni x\}.
\]
The latter are equivalence classes of the relation induced by the action; write $X/\cS$ for the associated quotient set.
The \emph{stabilizer} of $x \in X$ is 
\[
\Stab(x) = \{s \in \cS : \Dom(s)\ni x, \, s \cdot x = x\}/\sim 
\]
where if $s,t\in\cS$ and $x\in \Dom(s)\cap\Dom(t)$, then $s\sim t$ if there is $p\in\cE(\cS)$ such that $x\in\Dom(p)$ and $sp=tp$. The action $\cS \curvearrowright X$ gives rise to a groupoid $\cS \rtimes X$, the \emph{groupoid of germs} or \emph{transformation groupoid} of the action \cite{exel}*{Section 4}. This is defined as the quotient 
of $\cS \times X$ by the equivalence relation $(s,x) \sim (t,y)$ if $x = y$ and there exists
$p \in \mathcal{E}(\cS)$ such that $x\in\Dom(p)$ and $sp = tp$.
Units are given by $[p,x]$ with $p \in \mathcal{E}(\cS)$ and $x \in \Dom(p)$. As recalled above, idempotents in an inverse semigroup commute; hence given $e,f \in \mathcal{E}$ and $x \in \Dom(e) \cap \Dom(f)$ we have
$[e,x] = [ef,x] = [fe,x] = [f,x]$; thus $(\cS \rtimes X)^{(0)}$ 
can be homeomorphically identified with $X$ via $[p,x] \mapsto x$.
Sources and ranges are given by
$s([t,x]) = x$, $s([t,x]) = t \cdot x$,
composition by $[t',t \cdot x][t,x] = [t't,x]$
and inverses by $[t,x]^{-1} = [t^\ast, x
]$. Conditions for $\cS \rtimes X$ to be ample and Hausdorff are
given in \cite{steinappr}*{Definition 5.2 and Proposition 5.13} and \cite{steinappr}*{Theorem 5.17} 
respectively.
\begin{rem}\label{rem:stab=gxx}
    Remark that if $x\in X$ then for $\cG=\cS\ltimes X$ we have a bijection $\cG_x^x\cong \Stab(x)$, $[s,x]\mapsto [s]$.
    It follows that the product of $\cS$ makes $\Stab(x)$ into a group. 
\end{rem}
\begin{rmk}
Any ample groupoid $\cG$ arises as a germ groupoid construction via the action of the semigroup
$\mathcal{B}(\cG)$ of compact open bisections
on its unit space; if $U$ is a compact open bisection, then $\Dom(U) = s(U)$
and $U \cdot x = y$ if $r(s^{-1}(x) \cap U) = \{y\}$.
\end{rmk}

\begin{rmk} If $\Lambda$ is an abelian group and $c \colon S \to \Lambda$
a semigroup homomorphism, then $\cS \rtimes X$ is graded by $|[s,x]| = c(s)$.
\end{rmk}

\begin{defi} \label{defi:sparse}
We say that a semigroup action $\cS \curvearrowright X$
has \emph{sparse fixed points} if for each $s \in \cS \setminus \mathcal{E}(\cS)$
there exists at most one point $x \in \Dom(s)$ such that $s \cdot x = x$. 
\end{defi}

\begin{rmk} In Section \ref{sec:ep} below, we concentrate
on Exel-Pardo algebras, a particular family of Steinberg algebras arising from partial semigroup actions. We point out that for the subfamily consisting of Leavitt path algebras, the corresponding action has sparse fixed points; see Remark \ref{rmk:lpa-sparse} for further details. 
\end{rmk}

\begin{lem}\label{lem:sparse-disc-iso}
If $\cS \curvearrowright X$ is an inverse semigroup
action on a locally compact Hausdorff space
that has sparse fixed points, 
then $\Iso(\cS \rtimes X) \setminus X$ is discrete.
\end{lem}
\begin{proof}
Let $[s,x] \in \Iso(\cS \rtimes X) \setminus X$; in particular $s \not \in \mathcal{E}(\cS)$. Since the action has sparse fixed points, the subset
\[
[s,\Dom(s)] \cap \Iso(\cS \rtimes X) = 
\{[s,y] : y \in \Dom(s), s \cdot y = y\} = \{[s,x]\}.
\]
is open in $\Iso(\cS \rtimes X)$.
\end{proof}

\begin{thm} \label{thm:sparse}
Let $\cS$ be an inverse semigroup and $X$
a locally compact Hausdorff space. Suppose that $\cS \curvearrowright X$ is an action with sparse fixed points and that $\cS \rtimes X$ is both ample and Hausdorff. 
Fix a family $\cR\subset X$
of representatives for $X/\cS$ and for each $x\in\cR$ a set $Z_x$ 
of representatives of the non-trivial conjugacy classes of
$\Stab(x)$. Then there are quasi-isomorphisms of cyclic modules

\[
\H(\cS \rtimes X)\oplus \bigoplus_{x \in \cR} \bigoplus_{\eta \in Z_x} \H(\Stab(x)_{\eta}) \weq \H^{\cyc}(\cS \rtimes X).
\]
The grading on $\cS \rtimes X$ induced by a semigroup homomorphism $c \colon \cS \to \Lambda$
yields a decomposition
\[
{ _m \H}^{\cyc}(\cS \rtimes X)\sim \begin{cases}
\H(\cS \rtimes X) \oplus \bigoplus_{x \in \cR} \bigoplus_{\eta \in Z_x, c(\eta) = 0} \H(\Stab(x)_{\eta}) & m  = 0 \\
\bigoplus_{x\in X} \bigoplus_{\eta \in Z_x, , c(\eta) = m} { _m \H}(\Stab(x)_{\eta}) & \text{otherwise.}
\end{cases}
\]
\end{thm}
\begin{proof} In view of Remark \ref{rem:stab=gxx}, it suffices to point out that, by Lemma \ref{lem:sparse-disc-iso}, we are in position to apply Theorem \ref{thm:H-cyc}.
\end{proof}

\section{Exel-Pardo groupoids}\label{sec:ep}

In this section we concentrate on the Exel-Pardo groupoid $\cG$ associated to a self-similar action of a group $G$ on a directed graph $E$. We combine the results of the previous sections and some further results from \cite{aratenso} and \cite{eptwist} to describe the groupoid and Hochschild homology $\cG$ and of its Steinberg algebra of $\cG$, that is, the Exel-Pardo algebra of the action, and more generally of the twisted Steinberg algebra of groupoid cocycle twists of $\cG$, called a twisted Exel-Pardo algebra. In addition, we compute the $K$-theory of twisted Exel-Pardo algebras and relate it to groupoid homology. 

\numberwithin{equation}{subsection}

\subsection{Graphs}\label{subsec:graphs}
A (directed) \emph{graph} $E$ consists of sets $E^0$ and $E^1$ of \emph{vertices} and 
\emph{edges}, and \emph{source} and \emph{range} maps $s,r:E^1\to E^0$. A vertex $v$ 
\emph{emits} an edge $e$ if $v=s(e)$, and \emph{receives} it if $v=r(e)$. We say that 
$v$ is a \emph{sink} if it emits no edges, a \emph{source} if it receives no edges, 
and an \emph{infinite emitter} if it emits infinitely many edges. We write $\sink(E)$, 
$\sour(E)$ and $\inf(E)$ for the sets of sinks, sources and infinite emitters. The 
union $\sing(E)=\inf(E)\cup\sink(E)$ is the set of \emph{singular} vertices. 
Nonsingular vertices are called \emph{regular}; we write 
$\reg(E)=E^0\setminus\sing(E)$. We say that $E$ is regular if $E^0=\reg(E)$, \emph{row-finite} if $\inf(E)=\emptyset$ and \emph{finite} if both $E^0$ and $E^1$ are finite. 

A \emph{morphism of graphs} $f:E\to F$ consists of functions $f^i:E^i\to F^i$, $i=0,1$ such that $s\circ f^1=f^0\circ s$ and $r\circ f^1=f^0\circ r$. A \emph{subgraph} of a graph $E$
is a graph $F$ with $F^i\subset E^i$ such that the inclusions define a graph homomorphism $F\to E$, that is, if the source and range maps of $F$ are the restrictions of those of $E$. We say that a subgraph $F\subset E$ is \emph{complete} if $s^{-1}\{v\}\subset F^1$  for all $v\in\reg(F)\cap\reg(E)$. 

The \emph{reduced incidence matrix} of a graph $E$ is the matrix $A=A_E\in\N_0^{(\reg(E)\times E^0)}$ with coefficients
\[
A_{v,w}=\sharp\{e\in E^1\colon s(e)=v,\, r(e)=w\}.
\]
Let
\[
I\in\Z^{(E^0\times\reg(E))},\, I_{v,w}=\delta_{v,w}.
\]
The \emph{Bowen-Franks} group of $E$ is
\[
\BF(E)=\coker(I-A_E^t).
\]
A \emph{path} in a graph $E$ is a (finite or infinite) sequence $\alpha=e_1e_2\cdots$ such that $r(e_i)=s(e_{i+1})$ for all $i$. The source of $\alpha$ is $s(\alpha)=s(e_1)$; if $\alpha$ is finite of length $n$, we put $r(\alpha)=r(e_n)$ and $|\alpha|=n$. Vertices are considered as paths of length $0$. If $\alpha$ and $\beta$ are paths with $|\alpha|<\infty$, and $r(\alpha)=s(\beta)$, then we write $\alpha\beta$ for their \emph{concatenation}. If $\gamma$ is another path, we say that $\alpha$ \emph{precedes} $\gamma$ if $\gamma=\alpha\gamma_1$ for some path $\gamma_1$.

We write $\cP(E)$ for the set of all finite paths in $E$, which maybe regarded as the edges of a graph with $E^0$ as vertex set and the maps $s$ and $r$ defined above as source and range maps. If $v$ and $w$ are vertices and $n\in \N_0$, we consider the following subsets of $\cP(E)$
\begin{gather*}
\cP(E)_w=r^{-1}\{w\},\, \, 
\cP(E)^v=s^{-1}\{v\},\,\,\cP(E)^v_w=\cP(E)^v\cap\cP(E)_w,\\ \cP(E)_{n}=\{\alpha\in\cP(E)\colon |\alpha|=n\}, \cP(E)_{w,n}=\cP(E)_w\cap\cP(E)_{n}
\end{gather*}
and so on. Whenever $E$ is understood, we drop it from the notation and write $\cP$ for $\cP(E)$. For $n=1$ we use special notation; we put
\[
vE^1w=\cP^v_{w,1}.
\]

\subsection{Exel-Pardo tuples, twists and algebras}\label{subsec:ept}
Let $G$ be a group acting on a graph $E$ by graph automorphisms and $\phi:G\times E^1\to G$ a map satisfying
\begin{gather}
\phi(gh, e)=\phi(g,h(e))\phi(h,e),\label{eq:cocy}\\
\phi(g,e)(v)=g(v)\label{eq:cocynorma}    
\end{gather}
for all $g,h\in G$, $e\in E^1$ and $v\in E^0$. The first condition says that $\phi$ is a $1$-cocyle. We call the data $(G,E,\phi)$ an \emph{Exel-Pardo tuple} or simply an \emph{EP-tuple}. 

\begin{lem}[\cite{ep}*{Proposition 2.4}]\label{lem:epextend}
Let $(G,E,\phi)$ be an Exel-Pardo tuple. Then the $G$-action on $E$ and the cocycle $\phi$ extend respectively to a
$G$-action and a $1$-cocycle on the path graph $\cP(E)$ satisfying all four conditions below.
\begin{itemize}
\item[i)] $\phi(g,v)=g$ for all $v\in E^0$.
\item[ii)] $|g(\alpha)|=|\alpha|$ for all $\alpha\in\cP(E)$. 
\end{itemize}
The next two conditions hold for all concatenable $\alpha$, $\beta\in \cP(E)$. 
\begin{itemize}
\item[iii)] $g(\alpha\beta)=g(\alpha)\phi(g,\alpha)(\beta)$ 
\item[iv)] $\phi(g,\alpha\beta)=\phi(\phi(g,\alpha),\beta)$. 
\end{itemize}

Moreover, such an extension is unique.
\end{lem}

Any EP-tuple $(G,E,\phi)$ has an associated pointed inverse semigroup $\cS(G,E,\phi)$ \cite{ep}*{Definition 4.1}. Its nonzero elements are triples $\alpha g\beta^*$ where $g\in G$, $\alpha$ and $\beta$ are finite paths, $*$ is a (concatenation order reversing) involution, 
\begin{gather*}
\beta^*\gamma =\left\{\begin{matrix}\gamma_1& \gamma=\beta\gamma_1\\
\beta_1^*& \beta=\gamma\beta_1\\
0& \text{ else}\end{matrix}\right.\\
vg\cdot \alpha=\delta_{v,g(s(\alpha))}g(\alpha)\phi(g,\alpha),\text{ and }\alpha^*vg=\delta_{v,s(\alpha)}\phi(g,g^{-1}(\alpha))g^{-1}(\alpha)^*. 
\end{gather*}
The idempotent subsemigroup of $\cE=\cE(\cS(G,E,\phi))=\cE(\cS(E))=\{\alpha\alpha^*\colon \alpha\in\cP(E)\}$ is the usual idempotent semigroup of the graph $E$. We write 
$\hat{\fX}(E)$ for the set of all finite and infinite paths on $E$, equipped with the cylinder topology, of which a basis consists of the subsets of the form
\[
Z_\beta=\{\theta\in\hat{\fX}(E)\colon \beta\ge \theta\}
\]
indexed by the finite paths $\beta$ in $E$. 
Consider the closed subspace $\fX(E)\subset\hat{\fX}(E)$ consisting of all infinite paths and all paths ending at either a sink or an infinite emitter. An action of $\cS(G,E,\phi)$ on $\hat{\fX}(E)$ is defined as follows. An element $\alpha g\beta^*$ acts through the homeomorphism 
\[
Z_\beta\to Z_\alpha,\,\, \beta\gamma\mapsto \alpha g(\gamma).
\]
Here $g(\gamma)$ is as in  Lemma \ref{lem:epextend}. Remark that the above action leaves $\fX(E)$ invariant. 

\begin{rmk}\label{rmk:lpa-sparse}
If $G$ and $\phi$ are trivial, then the action of $\cS(E)$ on $\fX(E)$ has sparse fixed points in the sense of Definition \ref{defi:sparse}. Indeed, let $\alpha\beta^\ast \in \cS(E) \setminus \cE(\cS(E))$ and suppose that there exists $\beta\theta \in \fX(E)$ such that $\beta\theta = \alpha \beta^\ast \beta\theta = \alpha\theta$. Then there is a finite path $\gamma$
such that either $\alpha = \beta\gamma$ or $\beta = \alpha\gamma$. Since $\alpha \beta^\ast$ is not idempotent $\gamma$ has positive length and thus $\theta = \gamma^\infty = \gamma\gamma\gamma\cdots$. This shows that there exists at most one path that is fixed by $\alpha\beta^\ast$.
\end{rmk}

It is shown in \cite{ep}*{Section 8} that $\hat{\fX}(E)$ is $\cS(G,E,\phi)$-equivariantly homeomorphic to the \emph{spectrum} of the idempotent subsemigroup $\cE\subset\cS(G,E,\phi)$ and $\fX(E)$ to its \emph{tight spectrum} (\cite{exel}*{Definitions 10.1 and 12.8}). Thus the germ groupoids 
$$
\cG_u(G,E,\phi)=\cS(G,E,\phi)\ltimes \hat{\fX}(E)\,\text{ and }\,\cG(G,E,\phi)=\cS(G,E,\phi)\ltimes \fX(E)
$$
are respectively the \emph{universal} and the \emph{tight} or \emph{EP}-groupoid of $(G,E,\phi)$ in the sense of \cite{pater} and \cite{exel}. 

The \emph{Cohn} algebra of $(G,E,\phi)$ over a commutative ground ring $k$ is the semigroup algebra $C(G,E,\phi)=k[\cS(G,E,\phi)]$, with the $0$ element of the semigroup identified with that of the algebra. The \emph{$EP$-algebra} of $(G,E,\phi)$ is the Steinberg algebra $L(G,E,\phi)=\A(\cG(G,E,\phi))$. Next assume a $1$-cocycle
\[
c:G\times E^1\to \cU(k)
\]
taking values in the group of invertible elements is given. Then 
\[
\phi_c:G\times E^1\to \cU(k[G]),\, \phi_c(g,e)=c(g,e)\phi(g,e)
\]
is a $1$-cocycle. The data given by $G$,$E$, $\phi$ and $c$, which we abbreviate as $(G,E,\phi_c)$, is what we call a \emph{twisted EP-tuple}. It is shown in \cite{eptwist}*{Lemma 2.3.1} that $c$ extends uniquely to a $1$-cocycle $c:G\times \cP(E)\to \cU(k)$ satisfying
\begin{equation}\label{eq:cpalabra}
c(g,v)=1,\,\text{ and } c(g,\alpha\beta)=c(g,\alpha)c(\phi(g,\alpha),\beta)
\end{equation}
for all concatenable paths $\alpha,\beta$. Consider the pointed inverse semigroup $\cU(k)_\bullet=\cU(k)\cup\{0\}$. The extended map $c$ gives rise to a semigroup $2$-cocycle $\omega:\cS(G,E,\phi)\times \cS(G,E,\phi)\to\cU(k)_\bullet $ (see \cite{eptwist}*{Formula (2.4.5)}), which in turn induces a groupoid $2$-cocycle
$\bomega\colon\cG_u(G,E,\phi)^{(2)}\to\cU(k)$, 
\begin{equation}\label{eq:bomega}
\bomega([s,t(x)],[t,x])=\omega(s,t).     
\end{equation}

The same formula also defines a $2$-cocycle on $\cG(G,E,\phi)$, which we also call $\bomega$. We write 
\[
\cG_u(G,E,\phi_c)=(\cG_u(G,E,\phi),\bomega),\,\,\, \cG(G,E,\phi_c)=(\cG(G,E,\phi),\bomega)
\]
for the groupoids above equipped with the cocycles induced by $c$.
The \emph{twisted Cohn algebra} of $(G,E,\phi_c)$ is the twisted semigroup algebra $C(G,E,\phi_c)=k[\cS(G,E,\phi),\omega]$ of \cite{eptwist}. The \emph{twisted EP algebra} of $(G,E,\phi_c)$ is the twisted Steinberg algebra $L(G,E,\phi_c)=\A(\cG(G,E,\phi_c))$ which by \cite{eptwist}*{Section 3.4 and Proposition 4.2.2} is isomorphic to the quotient of $C(G,E,\phi_c)$ by the ideal $\cK(G,E,\phi_c)$ generated by the elements 
\begin{equation}\label{eq:qv}
qvg:=vg-\sum_{s(e)=v}ee^*vg=vg-\sum_{s(e)=v}e\phi_c(g,g^{-1}(e))g^{-1}(e)^*\, (v\in\reg(E)).
\end{equation}
Hence we have an algebra extension
\begin{equation}\label{ext:cohnext}
0\to \cK(G,E,\phi_c)\to C(G,E,\phi_c)\to L(G,E,\phi_c)\to 0.
\end{equation}
In fact it is shown in \cite{eptwist}*{Proposition 3.2.5} that $\cK(G,E,\phi_c)$ is independent of $c$. 
By \cite{eptwist}*{Proposition 3.2.5}, we have an isomorphism
\begin{equation}\label{map:ck=mat}
(\bigoplus_{v\in\reg(E)}M_{\cP_v})\rtimes G\iso \cK(G,E,\phi_c),\quad \epsilon_{\alpha,\beta}\rtimes g\mapsto \alpha (q_{r(\alpha)}g)(g^{-1}(\beta))^*.
\end{equation}
Here $G$ acts on the ultramatricial algebra above via $g(\epsilon_{\alpha,\beta})=\epsilon_{g(\alpha),g(\beta)}$.

\smallskip

Let $\reg(E)'$ be a copy of $\reg(E)$. Recall from \cite{lpabook}*{Definition 1.5.16} that the \emph{Cohn graph} of $E$ is the graph $\tilde{E}$ with $\tilde{E}^0=E^0\sqcup \reg(E)'$, $\tilde{E}^1=E^1\sqcup\{e'\colon r(e)\in\reg(E)\}$, where $r,s:\tilde{E}^1\to \tilde{E}^0$ extend the source and range maps of $E$, $s(e')=s(e)$ and $r(e')=r(e)'$. Extend the $G$-action and the cocycles $\phi$ and $c$ to $\tilde{E}$ via $g\cdot x'=(g\cdot x)'$ and $\phi(g,e')=\phi(g,e)$, $c(g,e')=c(g,e)$. In particular formula \eqref{eq:bomega} applied to the extended cocycle $c:G\times\tilde{E}^1\to\cU(k)$ defines a groupoid cocycle $\cG(G,\bar{E},\phi)^{(2)}\to\cU(k)$ which, by abuse of notation, we also call $\bomega$. 

\begin{lem}\label{lem:cohn=leav}
Let $U=\hat{\fX}(E)\setminus\fX(E)$ and $\cG'=\cG_u(G,E,\phi)_{|U}$. 
\item[i)] The cocycle $\bomega$ is trivial on $\cG'$ and $\cK(G,E,\phi_c)\cong\A(\cG')$.
\item[ii)] $C(G,E,\phi_c)\cong\cA(\cG_u(G,E,\phi_c))$.
\item[iii)] $\cG_u(G,E,\phi_c)\cong\cG(G,\tilde{E},\phi_c)$. 
\item[iv)] $C(G,E,\phi_c)\cong \A(\cG(G,\tilde{E},\phi_c))$.
\end{lem}
\begin{proof}
The groupoid $\cG'$ is discrete because $U$ is. One checks that every element 
of $\cG'$ is a germ $\xi=[\alpha g\beta^*,\beta]$ with $r(\alpha)=g(r(\beta))$ 
and that if $\xi=[\mu h\nu^*,\nu]$ with $r(\mu)=h(\nu)$ then we must have $
\alpha=\mu$, $\beta=\nu$ and $g=h$. The triviality of $\bomega$ on $\cG'$ follows from 
this and the definition of $\omega$ \cite{eptwist}*{Formula (2.4.5)}. 
One further checks, using the latter formula and the isomorphism 
\eqref{map:ck=mat}, that 
$\chi_{[\alpha g\beta^*,\beta]}\mapsto \alpha q_{g(\alpha)} g \beta^*$ defines an 
algebra isomorphism $\cA(\cG')\iso \cK(G,E,\phi_c)$. This proves i). By definition, the non-zero elements of $\cS(G,E,\phi)$ form a basis of $C(G,E,\phi_c)$. 
Hence there is a unique linear map $\pi:C(G,E,\phi_c)\to \cA(\cG_u(G,E,\phi_c))$ mapping $\alpha g\beta^*\mapsto \chi_{[\alpha g\beta^*,Z_\beta]}$. By \cite{eptwist}*{Proposition 3.1.5}, 
$C(G,E,\phi_c)$ is generated as an algebra by the elements $vg$, $eg$ and $ge^*$ $(v\in E^0,\, e\in E^1,\, g\in G)$ subject to the relations listed therein. One checks that the images under $\pi$ of said 
generators satisfy those relations and so $\pi$ is an algebra homomorphism, and furthermore that $\pi$ restricts on $\cK(G,E,\phi_c)$ to the isomorphism of part i). Remark that $\cG(G,E,\phi_c)=\cG_u(G,E,
\phi_c)_{|\fX(E)}$, and thus $\cA(\cG(G,E,\phi_c))\cong \cA(\cG_u(G,E,\phi_c))/\cA(\cG')$, so $\pi$ induces an algebra homomorphism 
$\bar{\pi}:L(G,E,\phi_c)\to\cA(\cG(G,E,\phi_c))$. By inspection, $\bar{\pi}$ is precisely the isomorphism of \cite{eptwist}*{Proposition 4.2.2}. Hence $\pi$ is an isomorphism, proving ii). Next observe that 
$\reg(\tilde{E})=\reg(E)$, $\inf(\tilde{E})=\inf(E)$ and $\sink(\tilde{E})=\sink(E)\cup\reg(E)'$. Hence the infinite paths and the paths ending in infinite emitters in $\hat{\fX}(E)$ and $\fX(\tilde{E})$ are the same, as are those in either space that end in a vertex of 
$\sink(E)$, while the paths in $E$ that end in $\reg(E)$ are in one-to-one correspondence with the paths in $\tilde{E}$ that end in $\reg(E)'$, via $v\mapsto v'$ and 
$\alpha=\alpha_1 e\mapsto \alpha'=\alpha_1 e'$. Altogether we get a bijection $a:\hat{\fX}(E)\iso\fX(\tilde{E})$. One checks that for $\beta\in\cP(E)$, $a$ sends
$Z_\beta$ to itself if $r(\beta)\notin\reg(E)$ and to $Z_{\beta}\cup\{\beta'\}$ otherwise. Hence $a$ is a homeomorphism. Extend $a$ to a map
\begin{gather}
 a:\cG_u(G,E,\phi)\to\cG(G,\tilde{E},\phi)\label{map:gutog}\\
 a[\alpha g\beta^*,\beta\gamma]=\left\{\begin{matrix} [\alpha' g(\beta')^*,\beta'\gamma]& \text{ if }s(\gamma)\in\reg(E)\\ [\alpha g\beta^*,\beta\gamma]& \text{ otherwise.}\end{matrix}\right.\nonumber
\end{gather} 
One checks that \eqref{map:gutog} is an isomorphism of topological groupoids that intertwines the corresponding groupoid cocycles, proving iii). Part iv) is immediate from ii) and iii). 
\end{proof}

In what follows we shall assume that $E$ is row-finite and that the group $G$ acts trivially on $E^0$. We shall abuse notation and write $\alpha g\beta^*$ for the image in $L(G,E,\phi_c)$ of the latter element of $\cS(G,E,\phi)$ via the projection $k[\cS(G,E,\phi),\omega]=C(G,E,\phi_c)\onto L(G,E,\phi_c)$.

\begin{lem}\label{lem:EPcoli} Let $(G,E,\phi_c)$ be a twisted EP-tuple such that $E$ is row-finite and $G$ acts trivially on $E^0$. Let $\cF$ be the set of all finite complete subgraphs of $E$, partially ordered by inclusion. Then 
\item[i)] For each $F\in\cF$, restriction of the action of $G$ and of the cocycle $\phi_c$ define a twisted EP-tuple $(G,F,\phi_c)$.
\item[ii)] The assignment $F\mapsto L(G,F,\phi_c)$ defines an $\cF$-directed system of $k$-algebras.
\item[iii)] $L(G,E,\phi_c)=\colim_{F\in\cF}L(G,F,\phi_c)$. 
\end{lem}
\begin{proof}
Because $G$ acts trivially on $E^0$ by hypothesis, it acts by permutation on $vE^1w$ for each $(v,w)\in\reg(E)\times E^0$. Hence every complete subgraph $F\subset E$ is invariant under the $G$-action. The cocycles $\phi$ and $c$ also restrict to maps on $G\times F^1$ which are again cocycles, since the cocycle condition \eqref{eq:cocy} passes down to $G$-invariant subgraphs. This proves i). Because $E$ is the filtering union of its finite complete subgraphs, $\cS(G,E,\phi)$ is the filtering union of the subsemigroups $\cS(G,F,\phi)$. Remark also that the semigroup cocycle $\omega$ restricts to a semigroup cocycle on each of these subsemigroups. Hence $C(G,E,\phi_c)=\bigcup_{F}C(G,F,\phi_c)=\colim_FC(G,F,\phi_c)$, where the union runs over the finite complete subgraphs. It is also clear that if $F\subset E$ is complete, then $\cK(G,F,\phi_c)\subset\cK(G,E,\phi_c)$ and that $\cK(G,E,\phi_c)=\bigcup_{F}\cK(G,F,\phi_c)$. Both ii) and iii) are immediate from this and exactness of filtering colimits. 
\end{proof}

\subsection{The degree zero component of \topdf{$L(G,E,\phi_c)$}{L} and the ideals \topdf{$I_v$}{Iv}}\label{subsec:l0}

Fix a twisted EP-tuple $(G,E,\phi_c)$ with $E$ row-finite and such that $G$ acts trivially on $E^0$. The algebra $L=L(G,E,\phi_c)$ is $\Z$-graded and its homogeneous component of degree zero, 
$L_0$, is the inductive union of the subalgebras 
\begin{equation}\label{eq:l0n}
L_{0,n}=\mspan\{\alpha g\beta^*\colon |\alpha|=|\beta|\le n,\, r(\alpha)=r(\beta)\},\, (n\ge 0).
\end{equation}

For each vertex $v\in E^0$ let $\iota_v:k[G]\to L(G,E,\phi_c)$ be the algebra homomorphism that sends an element $g\in G$ to the generator $vg\in L(G,E,\phi_c)$. Set 
\begin{gather}\label{I_v}
I_v=\ker(\iota_v),\,\, I=\bigoplus_{v\in E^0}I_v,\,\, 
R_v=\im(\iota_v)\cong k[G]/I_v,\,\, R=\bigoplus_{v\in E^0}R_v.
\end{gather}

By \cite{eptwist}*{Lemma 8.5} we have an isomorphism
\begin{gather}\label{map:8.5}
\bigoplus_{v\in \reg(E)}M_{\cP_{v,n}}R_v\oplus \bigoplus_{v\in \sink(E)}\bigoplus_{0\le  j\le n}M_{\cP_{v,j}}R_v
\iso  L_{0,n}
\end{gather}
that maps $\epsilon_{\alpha,\beta}g\mapsto \alpha g\beta^*$.

For each $v\in E^0$, let $k[G]v$ be a copy of $k[G]$. Let $n\ge 0$; put
\begin{equation}\label{eq:cmn}
\cM(G,E,\phi_c)_n=\bigoplus_{v\in \reg(E)}M_{\cP_{v,n}}k[G]v\oplus \bigoplus_{v\in \sink(E)}\bigoplus_{0\le  j\le n}M_{\cP_{v,j}}k[G]v
\end{equation}

Remark that for the matrix units $\epsilon_{\alpha,\beta}\in\cM(G,E,\phi_c)_n$ we have $r(\alpha)=r(\beta)$. Define a $k$-linear map
\begin{gather}
\jmath_n: \cM(G,E,\phi_c)_n\to \cM(G,E,\phi_c)_{n+1}\label{map:jmath}\\
\jmath_n(\epsilon_{\alpha,\beta}g)=\left\{\begin{matrix}
 \sum_{s(e)=r(\alpha)}\epsilon_{\alpha g(e),\beta e}\phi_c(g,e)& r(\alpha)\in\reg(E)\\
 \epsilon_{\alpha,\beta}& r(\alpha)\in\sink(E).\end{matrix}\right.\nonumber
\end{gather}

Put $\jmath_{\le n}=\jmath_n\circ\cdots\circ \jmath_0$, 
\[
I(n)=\ker(\jmath_{\le n}),\, I(n)_v=I(n)\cap k[G]v\, (v\in E^0).
\]
For $v,w\in E^0$, $h\in G$ and $\alpha,\beta\in\cP_{w,n}^v$, put
\[
G_{\alpha,\beta,h}=\{g\in G\colon g(\beta)=\alpha,\, \phi(g,\beta)=h\}.
\]
\begin{prop}\label{prop:Iv}
 Let $(G,E,\phi_c)$ be a twisted EP-tuple with $E$ row-finite and such that $G$ acts trivially on $E^0$. Also let $v\in E^0$ and $n\ge 0$.
 \item[i)] $\jmath_n$ is a homomorphism of $k$-algebras.
 \item[ii)] For $x=\sum_ga_gg\in k[G]v$, we have
 \begin{gather*}
 \jmath_{\le n}(x)=\sum_{w\in E^0}\sum_{\substack{\alpha,\beta\in \cP^v_{w,n},\\ h\in G}}(\sum_{g\in G_{\alpha,\beta,h}}a_gc(g,\beta))\epsilon_{\alpha,\beta}h+\\
\sum_{w\in\sink(E)}\sum_{j=0}^{n-1}\sum_{\substack{\alpha,\beta\in \cP^v_{w,j},\\ h\in G}}(\sum_{g\in G_{\alpha,\beta,h}}a_gc(g,\beta))\epsilon_{\alpha,\beta}h.
  \end{gather*}
 \item[iii)] The projections $k[G]\to R_v$ $(v\in E^0)$ together with the isomorphism \eqref{map:8.5} induce a commutative diagram with surjective vertical maps
 \[
 \xymatrix{
 \cM(G,E,\phi_c)_n\ar[d]_{\pi_n}\ar[r]^{\jmath_n}&\cM(G,E,\phi_c)_{n+1}
 \ar[d]^{\pi_{n+1}}\\ 
 L_{0,n}\ar[r]^{\inc}&L_{0,n+1}}
 \]
 
\item[iv)] $I(n)=\bigoplus_{v\in\reg(E)}I(n)_v$ and $I_v=\bigcup_n I(n)_v$. 
\item[v)] The natural map $$\colim_n\cM(G,E,\phi_c)_n\to L(G,E,\phi_c)_0$$ is an isomorphism of $k$-algebras.
\end{prop}
\begin{proof}
\item[i)] Remark that if $r(\alpha)=r(\beta)\ne r(\alpha')=r(\beta')$, then 
$$\jmath_n(\epsilon_{\alpha,\beta})\jmath_n(\epsilon_{\alpha',\beta'})=\jmath_n(\epsilon_{\alpha',\beta'})\jmath_n(\epsilon_{\alpha,\beta})=0$$
Hence it suffices to show that the restriction of $\jmath_n$ to each summand in the decomposition \eqref{eq:cmn} preserves products. This is clear for the summands corresponding to sinks. Let $v\in \reg(E)$, $\alpha,\beta,\gamma,\delta\in\cP_{v,n}$, and $g,h\in G$. Then 

\begin{gather*}
\jmath_n(\epsilon_{\alpha,\beta}g)\jmath_n(\epsilon_{\gamma,\delta}h)=\\
=\sum_{w\in E^0}\sum_{e,f\in vE^1w}\epsilon_{\alpha g(e),\beta e}\phi_c(g,e)\epsilon_{\gamma h(f),\delta f}\phi_c(h,f)\\
=\delta_{\beta,\gamma}\sum_{w\in E^0}\sum_{e\in vE^1w}\epsilon_{\alpha gh(e),\delta e} \phi_c(g,h(e))\phi_c(h,e)\\
=\delta_{\beta,\gamma}\sum_{w\in E^0}\sum_{e\in vE^1w}\epsilon_{\alpha gh(e),\delta e}\phi_c(gh,e)\\
=\delta_{\beta,\gamma}\jmath_n(\epsilon_{\alpha,\delta}gh)\\
=\jmath_n(\epsilon_{\alpha,\beta}g\epsilon_{\gamma,\delta}h)
\end{gather*}
\item[ii)]

\begin{gather*}
\jmath_{\le m}(x)=\\
\sum_{g\in G}\sum_{w\in E^0,\alpha,\beta\in \cP^v_{m,w}}\epsilon_{g(\alpha),\beta}a_g\phi_c(g,\alpha)+
\sum_{g\in G}\sum_{w\in \sink(E)}\sum_{j=0}^{m-1}\sum_{\beta\in \cP^v_{j,w}}\epsilon_{g(\beta),\beta}a_g\phi_c(g,\beta)=\\
\sum_{w\in E^0}\sum_{\alpha,\beta\in \cP^v_{w,n},h\in G}\epsilon_{\alpha,\beta}(\sum_{g\in G_{\alpha,\beta,h}}a_gc(g,\beta))h+\\
\sum_{w\in\sink(E)}\sum_{j=0}^{m-1}\sum_{\alpha,\beta\in \cP^v_{w,j},h\in G}\epsilon_{\alpha,\beta}(\sum_{g\in G_{\alpha,\beta,h}}a_gc(g,\beta))h.
\end{gather*}

\item[iii)] Straightforward.

\item[iv)] Fix $n\ge 0$. Let $v\in\reg(E)$ and let $\jmath_{n,v}$ be the restriction of $\jmath_n$ to $M_{\cP_{v,n}}k[G]v$. It is clear from the definition of $\jmath_n$ that $\im(\jmath_n)=\bigoplus_{v\in\reg(E)}\im(\jmath_{n,v})$. Hence $\ker(\jmath_n)=\bigoplus_{v\in\reg(E)}\ker(\jmath_{n,v})$ and therefore $I(n)=\bigoplus_{v\in\reg(E)}I(n)_v$. It is also clear that $I(n)\subset I(n+1)$, and it follows from ii) that $I(n)_v\subset I_v$ for all $v$. Let $0\ne x=\sum_{g\in G}a_gg\in I_v$; we shall show that $x\in I(m)$ for some $m$. 

The fact that $x\cdot v=0$ in $L(G,E,\phi_c)$ means that the product $x\cdot v\in C(G,E,\phi_c)$ belongs to $\cK(G,E,\phi_c)$. 
Hence we have an expression
\begin{gather*}
\sum_{g \in G}a_gg\cdot v=\sum_{w\in \reg(E)} \sum_{\substack{r(\alpha)=r(\beta)=w \\ h\in G}}b^h_{\alpha,\beta}\alpha hq_{w}\beta^*\\
=\sum_{h\in G}\sum_{\substack{r(\alpha)=r(\beta)\in\reg(E)}}(b^h_{\alpha,\beta}\alpha h\beta^*-\sum_{s(e)=r(\alpha)}b^h_{\alpha,\beta}c(h,e)\alpha h(e)\phi(h,e)e^*\beta^*)\\
=\sum_{\substack{r(\alpha)=r(\beta)\in\reg(E) \\ h\in G} }b^h_{\alpha,\beta}\alpha h\beta^*-\sum_{\substack{r(\alpha)=r(\beta)=s(e)=s(f)\in\reg(E) \\ h\in G}} \Bigg(\sum_{g\in G_{f,e,h}}b^g_{\alpha,\beta}c(g,e)\Bigg)\alpha f h(\beta e)^*.
\end{gather*}
Using that the non-zero elements of $\cS(G,E,\phi)$ are linearly independent in $C(G,E,\phi_c)$, we obtain that if $x\ne 0$, then $v\in\reg(E)$ and the following identities hold
\begin{equation}\label{eq:indstep}
 b_{v,v}^g=a_g,\, b_{\alpha f,\beta e}^h=\sum_{g\in G_{f,e,h}}b_{\alpha,\beta}^gc(g,e).
\end{equation}
Now a straightforward induction argument using \eqref{eq:indstep} and \eqref{eq:cpalabra} shows that
\[
b^h_{\alpha,\beta}=\sum_{g\in G_{\alpha,\beta,h}}a_gc(g,\beta)
\]
for all paths $\alpha$, $\beta$ and all $h\in G$. Next observe that if 
$$m-1=\max\{|\beta|\colon \exists \alpha,\, \beta,\, g\text{ such that } b^g_{\alpha,\beta}\ne 0\}$$ 
we must have
$b_{\alpha,\beta}^g=0$ for all $\alpha,\beta$ of length $m$. \comment{and also for those of smaller length ending in a sink.} Now apply ii).
\item[v)] By iii), $\jmath_n(\ker(\pi_n))\subset \ker(\pi_{n+1})$. It is clear fom the definitions that 
\begin{equation}\label{eq:kerpin}
\ker(\pi_n)=\bigoplus_{v\in \reg(E)}M_{\cP_{v,n}}I_v\oplus \bigoplus_{v\in \sink(E)}\bigoplus_{0\le  j\le n}M_{\cP_{v,j}}I_v.   
\end{equation}
By iv), $I_v=0$ if $v\in\sink(E)$. Remark that $M_{\cP_v}=k^{(\cP_v)}\otimes k^{(\cP_v)}$ and that if $v\in\reg(E)$ and $r(\alpha)=r(\beta)=v$, then 
\begin{gather*}
\jmath_n(\epsilon_{\alpha,\beta}g)=\sum_{w}\sum_{e\in vE^1w}\epsilon_{\alpha g(e),\beta e}\phi_c(g,e)\\
=\sum_{w}\sum_{e\in vE^1w}\epsilon_{\alpha}\otimes \epsilon_{g(e)}\otimes\epsilon_{\beta}\otimes \epsilon_{e}\otimes \phi_c(g,e),
\end{gather*}
which permuting tensors gets mapped to $\epsilon_{\alpha,\beta}\jmath_0(vg)$. Thus upon appropriate identifications, 
$\jmath_{[n,m+n]}:=\jmath_{n+m}\circ\cdots\circ\jmath_n$ is $\jmath_{\le m}$ applied entry-wise. Next use ii) and \eqref{eq:kerpin} to deduce that for every $x\in\ker(\pi_n)$, there exists an 
$m$ such that $\jmath_{[n,n+m]}(x)=0$. It follows that $\colim_n\ker(\pi_n)=0$, which implies v). 
\end{proof}

We recall from \cite{ep}*{Section 5} that a path $\alpha$ is said to be \emph{strongly
fixed} by an element $g \in G$
if $g(\alpha) = \alpha$ and
$\phi(g,\alpha) = 1$.

\begin{coro}\label{coro:Iv}
\item[i)] $I_v=0$ for all $v\in\sink(E)$. If $I_v\ne 0$, then there exists an $n\ge 1$ such that for all $\beta\in \cP^v_n$ there is $G\owns g_\beta\ne 1$ that fixes $\beta$ strongly. 
\item[ii)] Assume that $k[G]$ is Noetherian. Then for every $v\in\reg(E)$ there exists an $n=n_v$ such that $\jmath_{\le n}$ induces an embedding $R_v\to \cM(G,E,\phi_c)_n$. 
\end{coro}
\begin{proof}
Both ii) and the first assertion of i) are immediate from part iv) of Proposition \ref{coro:Iv}. Next assume there exists $0\ne x=\sum_{g\in G}a_gg\in I_v$. Let $g_1\in G$ such that $a_{g_1}\ne 0$. Let $n$ be minimal such that $x\in I(n)_v$ and $\beta\in \cP^v_n$. Set $\alpha=g_1(\beta)$, $h=\phi(g_1,\beta)$. Then $g_1\in G_{\alpha,\beta,h}$, and since $a_{g_1}\ne 0$, by part ii) of Proposition \ref{prop:Iv}, there must exist $g_2\ne g_1\in G_{\alpha,\beta,h}$ such that $a_{g_2}\ne 0$. Then $g=g_2^{-1}g_1\ne 1$ and fixes $\beta$ strongly.
\end{proof}

\begin{rem}\label{rem:Iv}
Let $v\in \reg(E)$ and assume that there is an $n\ge 1$ and an element $1\ne g\in G$ that strongly fixes all 
\[
\beta\in \cP^v_n\cup\bigcup_{j\le n, w\in\sink(E)} P^v_{w,j}
\] 
simultaneoulsy and that $c(g,\beta)=u$ for all such $\beta$. Then $0\ne (u-g)\in I_v$. 
\end{rem}
\begin{ex}\label{ex:Iv}
Let $n\ge 2$, and let $E$ be the graph with $E^0=\{v,w\}$, $E^1=\{e_1,\dots,e_n\}$ with $s(e_i)=v$ and $r(e_i)=w$ for all $i$. Let the symmetric group $\S_n$ act on $E^1$ by permutation of subindices; let $\rho:k[\S_n]\to \End_k(k[E^1])\cong M_{\cP_{w,1}}$ be the corresponding representation. Assume, for simplicity, that $k$ is a domain. Then, by reasons of rank, $\ker(\rho)\ne 0$ for $n\ge 4$. Equip $(\S_n,E)$ with trivial $\phi$ and $c$, and let $L=L(\S_n,E)$. For $x=\sum_ga_gg$, we have
\[
\rho(x)=\sum_{g,e}a_g\epsilon_{g(e),e}=\sum_{e,f\in E^1}(\sum_{g\in G_{f,e,1}}a_g)\epsilon_{f,e}.
\]
Hence by part ii) of Proposition \ref{prop:Iv}, we have $I_v=\ker(\rho)$, which is nonzero for $n\ge 4$. Note however that there is no nontrivial element of $\S_n$ strongly fixing all the edges of $E$ simultaneously.
\end{ex}
\begin{rem}\label{rem:phibiende}
Let $(G,E,\phi_c)$ be as in Proposition \ref{prop:Iv}. Pick $v\in \reg(E)$ and $x=\sum_ga_gg\in I_v$. By \eqref{eq:kerpin} and part iii) of the proposition, we have
\[
\jmath_0(x)\in \bigoplus_{w\in E^0}M_{vE^1w}I_w.
\]
Hence by part ii) of the same proposition, 
\[
\sum_{\{g\colon g(e)=f\}}a_g\phi_c(g,e)\in I_w\,\, (\forall\, e,f\in vE^1w).
\]
\end{rem}

\subsection{Hochschild homology of Exel-Pardo algebras}\label{subsec:hhep}
Let $I$ be as in \eqref{I_v}. For $X\subset E^0$, put
\begin{gather*}
k[G]_X=\bigoplus_{v\in X}k[G]v\cong k[G]\otimes k^{(X)}, \\
I_X=I\cap k[G]_X,\, R_X=k[G]_X/I_X.
\end{gather*}
Whenever the graph $E$ is clear from the context, we shall drop it from the subscript of $R$,  and write $R$ for $R_{E^0}$, $R_{\reg}$ for $R_{\reg(E)}$ and so on. 

Make the right $k[G]_{E^0}$-module $S_m=k^{(\cP_m)}\otimes_{k^{(E^0)}}k[G]_{E^0}$ into a $k[G]_{E^0}$-bimodule with the left multiplication induced by
 \begin{equation}\label{eq:leftSm}
 vg\cdot(\alpha\otimes h)=\delta_{v,s(\alpha)}g(\alpha)\otimes\phi_c(g,\alpha)h.    
 \end{equation}
 
 Similarly, make the left $k[G]_{E^0}$-module $S_{-m}=k[G]_{E^0}\otimes_{k^{(E^0)}} k^{(\cP_m^*)}$ into a bimodule via
 \begin{equation}\label{eq:rightSm}
 (g\otimes \beta^*)\cdot h=g\phi_c(h,h^{-1}(\beta))\otimes h^{-1}(\beta)^*.
 \end{equation}
 
Let $\ell\subset k$ be a unital subring such that $k$ is flat over $\ell$. For $m\in\Z$, set $S_m^{\reg}=k[G]_{\reg(E)}S_m k[G]_{\reg(E)}$ and define a chain complex homomorphism 
 \begin{equation}\label{map:sigmam}
  \sigma_m:\HH(k[G]_{\reg(E)}/\ell^{(\reg(E))},S^{\reg}_m)\to \HH(k[G]_{E^0}/\ell^{(E^0)},S_m)   
 \end{equation}
 as follows. For $a_0,\dots,a_n\in k$ and $g_0,\dots,g_n\in G$, set
 \begin{gather*}
 \sigma_0(va_0g_0\otimes\cdots\otimes va_ng_n)=
 \\
\sum_{\scst{\begin{matrix}w\in E^0,\\ s(e)=v,\, r(e)=w,\\ (g_0\cdots g_n)(e)=e\end{matrix}}}w\phi_c(g_0,g_1\cdots g_n(e)) a_0\otimes w\phi_c(g_1,g_2\cdots g_n(e))a_{1}\otimes\cdots\otimes w\phi_c(g_n,e)a_n.
 \end{gather*}
For $m\ge 1$, if $\alpha=e_1\cdots e_m$ is a path with $s(\alpha)=v\ne w=r(\alpha)$, then the element
 \begin{equation}\label{eq:notclosed=0}
 vg_0 a_0\otimes\cdots\otimes vg_{n-1}a_{n-1}\otimes e_1\cdots e_m\otimes wg_n a_n=0    
 \end{equation}
  in $\HH(k[G]_{\reg(E)}/\ell^{(\reg(E))},S_m)$. If $v=w$ and $r(e_1)=u$, put
 \begin{gather*}
\sigma_m(vg_0 a_0\otimes\cdots\otimes vg_{n-1}a_{n-1}\otimes e_1\cdots e_m\otimes vg_n a_n)\\
 =\phi_c(g_0,g_1\cdots g_{n-1}(e_1))ua_0\otimes \cdots\otimes \phi_c(g_{n-1},e_1)ua_{n-1}\\
 \otimes e_2\cdots e_m (g_ng_0\cdots g_{n-1}(e_1))\otimes\phi_c(g_n,g_0\cdots g_{n-1}(e_1))u a_n,
 \end{gather*}
 \begin{gather*}
 \sigma_{-m}(vg_0a_0\otimes\cdots\otimes v g_{n-1}a_{n-1}\otimes g_n a_n\otimes (e_1\cdots e_m)^*) =\\
u\phi_c(g_0,g_0^{-1}(e_1))a_0\otimes\cdots \otimes u\phi_c(g_n, (g_0\cdots g_n)^{-1}(e_1))a_{n}\otimes (e_2\cdots e_r (g_0\cdots g_n)^{-1}(e_1))^*.
 \end{gather*}
 By Remark \ref{rem:phibiende}, \eqref{eq:leftSm} and \eqref{eq:rightSm} also define $R$-bimodule structures on $\bar{S}_m=k^{(\cP_m)}\otimes_{k^{(E^0)}}R$ and $\bar{S}_{-m}=R\otimes_{k^{(E^0)}} k^{(\cP_m^*)}$ for all $m\in\N_0$, so that for all $n\in\Z$, the chain map $\sigma_n$ descends to a chain map
 \[
 \bar{\sigma}_n:\HH(R_{\reg}/\ell^{(\reg(E))},\bar{S}^{\reg}_n)\to \HH(R/\ell^{(E^0)},\bar{S}_n).
 \]
Here ${\bar{S}_n}^{\reg}=R_{\reg}\bar{s}_n R_{\reg}$.

\begin{rem}\label{rem:smclosed}
For $m\ge 0$ let
\begin{equation}\label{eq:sm}
CP_m(E)=\{\alpha\in \cP(E)_m\colon s(\alpha)=r(\alpha)\},\, CP_{m}(E)^*=\{\alpha^*\colon \alpha\in CP_m(E)\}.
\end{equation}
Consider the sub-bimodules 
\[
S_m\supset S^{c}_m=k^{(CP_m(E))}\otimes_{k^{(E^0)}}k[G]_{E^0}, \,\, S_{-m}\supset S_{-m}^c=k[G]_{E^0}\otimes_{k^{(E^0)}}k^{(CP_m(E)^*)}.
\]
Remark that $S_0=S^{c}_0$ and that for $m\ne 0$, $S_m^c\subset S_m^{\reg}$. Moreover, it follows from \eqref{eq:notclosed=0} that for $m\ne 0$ the inclusion $S_m^c\subset S_m$ induces chain complex isomorphisms
\begin{equation}\label{map:smc}
 \HH(k[G]_{\reg(E)}/\ell^{(\reg(E))}, S_m^c)\cong \HH(k[G]_{\reg(E)}/\ell^{(\reg(E))}, S_m^{\reg})\cong \HH(k[G]_{E^0}/\ell^{(E^0)}, S_m).
\end{equation}
Similarly, 
\[
\HH(R_{\reg}/\ell^{(\reg(E))}, S_m^c)\cong \HH(R_{\reg}/\ell^{(\reg(E))}, S_m)\cong \HH(R/\ell^{(E^0)}, S_m).
\]
Furthermore,  $S_m$ and $S^c_m$ are also bimodules over $\prod_{v\in E^0}k[G]$, so we may consider them as $k[G]$-bimodules, by restriction along the diagonal embedding $k[G]\to\prod_{v\in E^0}k[G]$. For all $m\in\Z$, we have a $k[G]$-bimodule decompositon $S_m^c=\bigoplus_{v\in E^0}vS_m^cv$. Mapping
\[
a_0g_0\otimes\cdots a_{n-1}g_{n-1}\otimes vsv\mapsto 
a_0g_0v\otimes\cdots a_{n-1}g_{n-1}v\otimes vsv
\]
we obtain chain complex isomorphisms
\begin{gather}
\HH(k[G]_{\reg(E)}/\ell^{(\reg(E))},S_m^{c})\cong \HH(k[G]/\ell, S_m^{c})\,\, (m\ne 0),\label{map:k[G]mne0}\\
\HH(k[G]_{E^0}/\ell^{(E^0)}, S_0)\cong \HH(k[G]/\ell)^{(E^0)},\label{map:k[G]m0E}\\
\HH(k[G]_{\reg}/\ell^{(\reg(E))}, S^{\reg}_0)\cong \HH(k[G]/\ell)^{\reg(E)}.\label{map:k[G]m0reg}
\end{gather}
\end{rem}

\smallskip

\begin{rem}\label{rem:pullc}
By definition, $\phi_c(g,e)=c(g,e)\phi(g,e)$, where $\phi(g,e)\in G$ and

\goodbreak
\noindent $c(g,e)\in \cU(k)$. Hence if we set $\ell=k$, the map $\sigma_m$ becomes $k$-linear, so it is determined by its value for $a_0=\dots=a_n=1$, and we may gather all the $c$'s together into a scalar and substitute $\phi$ for $\phi_c$ everywhere. For example, if we do this with the formula for $\sigma_0$ and set all $a_i=1$, then using the cocycle equation \eqref{eq:cocy}, the term of the sum corresponding to an edge $e\in vE^1w$ becomes
\[
c(g_0\cdots g_n,e)\phi(g_0,g_1\cdots g_n(e))w\otimes\cdots\otimes\phi(g_n,e)w.
\]
\end{rem}
\begin{thm}\label{thm:hhep}
 Let $(G,E,\phi_c)$ be a twisted EP-tuple. Assume that $E$ is row-finite and that the group $G$ acts trivially on $E^0$. Let $\ell\subset k$ be a flat ring extension and let $L(G,E,\phi_c)$ be the Exel-Pardo $k$-algebra. Let
 $\HH(L(G,E,\phi_c)/\ell)=\bigoplus_{m\in\Z}{}_m\HH(L(G,E,\phi_c)/\ell)$ be the weight decomposition associated to the natural $\Z$-grading of $L(G,E,\phi_c)$. Then for every $m\in\Z$ there are natural zig-zags of quasi-isomorphisms
 \begin{gather}
 \cone(\HH(k[G]_{\reg(E)}/\ell^{(\reg(E))},S_m)\overset{1-\sigma_m}{\lra}\HH(k[G]_{E^0}/\ell^{(E^0)},S_m))\weq\label{map:conGconR}\\
  \cone(\HH(R_{\reg}/\ell^{(\reg(E))},\bar{S}_m)\overset{1-\bar{\sigma}_m}{\lra}\HH(R/\ell^{(E^0)},\bar{S}_m))\weq\label{map:conRHH}\\
  {}_m\HH(L(G,E,\phi_c)/\ell^{(E^0)}).\nonumber
 \end{gather}
\end{thm}
\begin{proof}

Part 1: proof of \eqref{map:conRHH}.

\medskip

Set $L=L(G,E,\phi_c)$ and let $L=\bigoplus_{m\in\Z}L_m$ be the $\Z$-grading. 

\medskip

\noindent{Step 1: \em $E$ finite without sources.} Pick an edge $e_v\in r^{-1}(\{v\})$ for each $v\in E^0$. Set $t_+=\sum_{v\in E^0}e_v$, $t_-=t_+^*$. Then $|t_+|=1$, $|t_-|=-1$ and
$t_-t_+=1$. Hence $\psi:L_0\to L_0$, $\psi(a)=t_+at_-$ is an isomorphism onto the corner associated to the idempotent $t_+t_-$ and thus $L$ is isomorphic to the skew Laurent polyomial 
algebra $L_0[t_+,t_-;\psi]$ of \cite{fracskewmon}. Hence by Proposition \ref{prop:skewcorner} there is a quasi-isomorphism
\begin{equation}\label{map:coneoverl}
\cone (\HH(L_0/\ell,L_m)\overset{1-\psi}{\lra} \HH(L_0/\ell,L_m))\weq {}_m\HH(L(G,E,\phi_c)/\ell)
\end{equation}
Recall that $L_0=\bigcup_nL_{0,n}$ is the increasing inductive union of the algebras \eqref{eq:l0n}.
For $m\in\Z$, set 
$$
L_{m,n}=L_{0,n}\bar{S}_mL_{0,n}.
$$ 
Thus $L_m=\bigcup_{n\ge 0}L_{m,n}$ for all $m\in\Z$. Recall from \eqref{map:8.5} that $L_{0,n}$ is a direct sum of matrix algebras, whose coefficients lie in the ring
\begin{align}\label{eq:Rn}
R_n=&\bigoplus_{v\in\reg(E)}R_v\oplus\bigoplus_{v\in\sink(E)}\bigoplus_{0\le j\le n}R_v\\
=&R_{\reg}\oplus R_{\sink}\otimesl\ell^{n+1}.\nonumber
\end{align}
For $0\le j\le n$ and $v\in\sink(E)$ we write  $R_{(v,j)}=R_v\otimes\epsilon_j$ for the $j$-th copy of $R_v$ in the direct sum above. Set $[n]=\{0,\dots,n\}$, 
$$E^0_n=\reg(E)\bigsqcup (\sink(E)\times[n]).$$
Let $\iota_n:L_{0,n}\to L_{0,n+1}$ be the inclusion map. Because $\HH$ commutes with filtering colimits and the algebra $\ell^{E^0_n}$ is separable, we have quasi-isomorphisms
\begin{gather}\label{map:hhepcoli}
\cone (\HH(L_0/\ell,L_m)\overset{1-\psi}{\lra} \HH(L_0/\ell,L_m))\\
=\colim_n\cone (\HH(L_{0,n}/\ell,L_{m,n})\overset{\iota_n-\psi}{\lra} \HH(L_{0,n+1}/\ell,L_{m,n+1}))\nonumber\\
\weq \colim_n\cone (\HH(L_{0,n}/\ell^{E^0_n},L_{m,n})\overset{\iota_n-\psi}{\lra} \HH(L_{0,n+1}/\ell^{E^0_{n+1}},L_{m,n+1})).\nonumber
\end{gather}
Put
\begin{gather*}
P_n=\mspan\{\alpha g\colon |\alpha|=n \text{ or } 0\le |\alpha|<n \text{ and } r(\alpha)\in\sink(E)\}\\
Q_n=\mspan\{g\alpha^*\colon |\alpha|=n \text{ or } 0\le |\alpha|<n \text{ and } r(\alpha)\in\sink(E)\}.
\end{gather*}
Then $P_n$ is an $(L_{0,n},R_n)$-bimodule and $Q_n$ an $(R_n,L_{0,n})$-bimodule, which correspond under the isomorphism \eqref{map:8.5} to the direct sums of the obvious bimodules of row and column vectors. In particular we have bimodule isomorphisms
\[
P_n\otimes_{R_n}Q_n\cong L_{0,n}, \, Q_n\otimes_{L_{0,n}}P_n\cong R_n.
\]
Assume that $n\ge m\ge 0$. Regard the $k$-module
\[
T_{m,n}=\bigoplus_{(v,j)\in\sink(E)\times[n-m]} R_v\otimes\epsilon_{v,j}=R_{\sink}\otimes_{k^{\sink(E)}}k^{\sink(E)\times[n-m]}
\]
as an $R_n$-bimodule with left and right multiplication defined as follows. For $x\in R_n$, $v\in\sink(E)$ and $0\le i\le n$, let $x_{(v,i)}\otimes\epsilon_i$ be the component of $x$ in $R_{(v,i)}=R_v\otimes\epsilon_i$; if $z\otimes\epsilon_{v,j}\in R_v\otimes\epsilon_{v,j}$, set
\[
x\cdot (z\otimes\epsilon_{v,j})=x_{(v,j+m)}z\otimes\epsilon_{v,j},\,\,\,\, (z\otimes\epsilon_{v,j})\cdot x=zx_{(v,j)}\otimes\epsilon_{v,j}. 
\]
For $0\le m$ consider the algebra homomorphisms $\pi^l_{m,n},\pi^r_{m,n}:R_n\to R$, defined as follows
\begin{gather*}
\pi^l_{m,n}(x)_v=\left\{\begin{matrix} x_v & \text{ if } v\in\reg(E)\\ x_{(v,0)} & \text{ if }v\in\sink(E)\text{ and } m=0.\end{matrix}\right.\\
\pi^r_{m,n}(x)_v=\left\{\begin{matrix} x_v & \text{ if } v\in\reg(E)\\ x_{(v,n)} & \text{ if }v\in\sink(E).\end{matrix}\right.\\
\end{gather*}
In the following paragraph we will regard the $R$-bimodule $S_m$ it as $R_n$-bimodule where an element $x$ acts via $\pi^l_{m,n}(x)$ on the left and via $\pi^r_{m,n}(x)$ on the right. To regard $S_{-m}$ as an $R_n$ module, we switch the roles of $\pi^r_{m,n}$ and $\pi^l_{m,n}$; we use the former for the left multiplcation and the latter for right multiplication. 
One checks that we have isomorphisms of $R_n$-bimodules
\begin{gather}
L_{m,n}=\bigoplus_{v\in\reg(E)}\mspan\{\alpha vg\beta^*\colon |\alpha|=m+n, |\beta|=n, 
r(\alpha)=r(\beta)=v\}\label{eq:lmn}\\
\oplus \bigoplus_{\scst{v\in\sink(E),\, 0\le j\le n}}\mspan\{\alpha vg\beta^*\colon |\alpha|
=m+j, |\beta|=j, r(\alpha)=r(\beta)=v\}\nonumber\\
\cong \bigoplus_{v\in\reg(E)}\bar{S}_{m+n}\otimes_{R_v}\bar{S}_{-n}\oplus\bigoplus_{\scst{v\in\sink(E)}}\bigoplus_{j=0}^{n}\, 
\bar{S}_{m+j}\otimes_{R_v}\bar{S}_{-j}\nonumber
\\
\cong P_n\otimes_{R_n}(\bar{S}_m\oplus T_{m,n})\otimes_{R_n}Q_n.\nonumber
\end{gather}
Similarly, we write $T_{-m,n}$ for the same $k$-module $R_{\sink}\otimes_{k^{\sink(E)}}k^{\sink(E)\times [n-m]}$, but where now 
$R_n$ acts on $R_v\otimes \epsilon_{v,j}$ via $R_{v,j}$ on the left and via $R_{v,j+m}$ on the right, and we have an isomorphism
\[
L_{-m,n}\cong P_n\otimes_{R_n}(\bar{S}_{-m}\oplus T_{-m,n})\otimes_{R_n}Q_n.
\]
Hence for all $m\in \Z$ there is a trace quasi-isomorphism \cite{loday}*{Definition 1.2.1}
\begin{equation}\label{map:trace}
\tr:\HH(L_{0,n}/\ell^{E_n^0}, L_{m,n})\weq \HH(R_n/\ell^{E^0_n},\bar{S}_m\oplus T_{m,n}).   
\end{equation}

Grouping the summands corresponding to regular vertices together in one summand and 
those corresponding to sinks on the other as in \eqref{eq:Rn} and \eqref{eq:lmn}, we get a decomposition
$L_{m,n}=L_{m,n}^{\reg}\oplus L_{m,n}^{\sink}$, and 
the trace map is homogeneous with respect to these decompositions. For $m\ne 0$, we have
$$\HH(R_{\sink}\otimes\ell^{n+1}/\ell^{\sink(E)}\otimesl\ell^{n+1},T_{m,n})=\HH(R_{\sink}/\ell^{\sink(E)},\bar{S}_m)=0.$$
Hence 
$\HH(L_{0,n}/\ell^{E_n^0}, L_{m,n})$ decomposes into the direct sum of a complex quasi-isomorphic to zero and a copy of $\HH(L_{0,n}^{\reg}/\ell^{\reg(E)},L^{\reg}_{m,n})$, and the trace is a quasi-isomorphism
\[
\tr:\HH(L_{0,n}^{\reg}/\ell^{\reg(E)}, L^{\reg}_{m,n})\weq \HH(R/\ell^{\reg(E)},\bar{S}_m).
\]
For $m\ge 0$, the latter map sends
\begin{gather}\label{map:tr}
\tr(\alpha_0 g_0a_0\beta^*_0\otimes\cdots \otimes\alpha_{n-1} g_{n-1}a_{n-1}\beta_{n-1}^*\otimes \alpha_n e_1\cdots e_m a_ng_n\beta^*_n)\\
=(\prod_{i=0}^n\delta_{\beta_i,\alpha_{i+1}})\cdot g_0a_0\otimes\cdots\otimes g_{n-1}a_{n-1}\otimes e_1\cdots e_m g_na_n.\nonumber
\end{gather}
A similar formula holds for $m<0$. Observe that $\iota_n$ restricts to the 
obvious inclusion $L_{0,n}^{\sink}\subset L_{0,n+1}^{\sink}$ and is induced by the second Cuntz-Krieger relation
\[
vg=\sum_{s(e)=v}g(e)\phi(g,e)e^*
\]
on $L_{0,n}^{\reg}$. Using this together with the explicit formula \eqref{map:tr} and its analog for $m<0$, we obtain that for $m\ne 0$ the following diagrams commute
\begin{gather}\label{diag:dos}
\xymatrix{\HH(L_{0,n}/\ell^{E^0_n},L_{m,n})\ar[d]^{\wr\tr}\ar[r]^(.45){\iota_n}&\HH(L_{0,n+1}/\ell^{E_{n+1}^0},L_{m,n+1})\ar[d]^{\wr\tr}\\
        \HH(R_{\reg}/\ell^{\reg(E)},\bar{S}_m)\ar[r]^{\bar{\sigma}_m}& \HH(R_{\reg}/\ell^{\reg(E)},\bar{S}_m)}\\
        \xymatrix{\HH(L_{0,n}/\ell^{E^0_n},L_{m,n})\ar[dr]^{\tr}_{\sim}\ar[rr]^{\psi}&&\HH(L_{0,n+1}/\ell^{E_{n+1}^0},L_{m,n+1})\ar[dl]_{\tr}^{\sim}\\
        &\HH(R_{\reg}/\ell^{\reg(E)},\bar{S}_m)&}\label{diag:dos1}      
\end{gather}

Let Hence it follows from Remark \ref{rem:smclosed} and Lemma \ref{lem:locacono} that we have a quasi-isomorphism
\begin{gather}\label{map:locaconohhep}
\cone(\HH(R_{\reg}/\ell^{\reg(E)},\bar{S}_m)\overset{1-\bar{\sigma}_m}{\lra}\HH(R/\ell^{E^0},\bar{S}_m))\weq\\
\cone(\HH(R_{\reg}/\ell^{\reg(E)},\bar{S}_m)[\sigma_m^{-1}]\overset{1-\bar{\sigma}_m}{\lra}\HH(R/\ell^{E^0},\bar{S}_m)[\sigma_m^{-1}])\weq\\
\colim_n\cone (\HH(L_{0,n}/\ell,L_{m,n})\overset{\psi-\iota_n}{\lra} \HH(L_{0,n+1}/\ell,L_{m,n+1}))\cong\nonumber\\
\colim_n\cone (\HH(L_{0,n}/\ell,L_{m,n})\overset{\iota_n-\psi}{\lra} \HH(L_{0,n+1}/\ell,L_{m,n+1})).\nonumber
\end{gather}

As a preliminary to the case $m=0$, observe that for all $n\ge 0$ we have a direct sum decomposition
\[
\HH(R_n/\ell^{E^0_n})=\HH(R_{\reg}/\ell^{\reg(E)})\oplus \HH(R_{\sink}/\ell^{\sink(E)})\otimesl\ell^{[n]},
\]
We use the decomposition above to define chain homomorphisms $f_n, g_n:\HH(R_n/\ell^{E^0_n})\to \HH(R_{n+1}/\ell^{E^0_{n+1}})$ as follows. 
On the summand $\HH(R_{\reg}/\ell^{\reg(E)})$, $f_n$ restricts to $\bar{\sigma}_0$ and $g_n$ to the identity map. Both $f_n$ and $g_n$ restrict to maps
$\HH(R_{\sink}/\ell^{\sink(E)})\otimesl\ell^{[n]})\to \HH(R_{\sink}/\ell^{\sink(E)})\otimesl\ell^{[n+1]}$ 
and as such have the following matricial forms
\begin{gather*}
f_n=\sum_{j=0}^n\epsilon_{j,j}, \,\, g_n=\sum_{j=0}^{n+1}\epsilon_{j+1,j}.
 \end{gather*}
One checks that the following diagrams commute
\begin{gather}\label{diag:dosprima1}
 \xymatrix{\HH(L_{0,n}/\ell^{E^0_n})\ar[d]^{\tr}_{\wr}\ar[r]^{\psi}&\HH(L_{0,n+1}/\ell^{E_{n+1}^0})\ar[d]_{\tr}^{\wr}\\
        \HH(R_n/\ell^{E^0_n})\ar[r]^{g_n}& \HH(R_{n+1}/\ell^{E^0_{n+1}})}\\
\xymatrix{\HH(L_{0,n}/\ell^{E^0_n})\ar[d]^{\wr\tr}\ar[r]^(.45){\iota_n}&\HH(L_{0,n+1}/\ell^{E_{n+1}^0})\ar[d]^{\wr\tr}\\
        \HH(R/\ell^{E^0_n})\ar[r]^{f_n}& \HH(R/\ell^{E^0_{n+1}})}\label{diag:dosprima}             
\end{gather}

Hence the trace map induces a quasi-isomorphism $\cone(\iota_n-\psi)\weq \cone(f_n-g_n)$, $\cone(1-\bar{\sigma}_0)\subset\cone(f_n-g_n)$ is a subcomplex, and $\cone(f_n-g_n)/\cone(1-\bar{\sigma}_0)$ is the cone of the map
\begin{gather*}
h_n:\HH(R_{\sink}/\ell^{\sink(E)})\otimesl\ell^{[n]}\to \HH(R_{\sink}/\ell^{\sink(E)})\otimesl\ell^{[n+1]}\\
h_n=\sum_{j=0}^n(\epsilon_{j,j}-\epsilon_{j+1,j}).
\end{gather*}
Since $\colim_nh_n$ is an isomorphism, its cone is contractible, and thus we have a zig-zag of quasi-isomorphisms as follows
\begin{gather*}
\colim_n\cone(\HH(L_{0,n}/\ell^{E^0_n})\overset{\iota_n-\psi}{\lra}\HH(L_{0,n}/\ell^{E^0_n}))\weq\\
\colim_n\cone(\HH(R_n/\ell^{E^0_n})\overset{f_n-g_n}{\lra}\HH(R_{n+1}/\ell^{E^0_{n+1}})\\
\overset{\sim}{\longleftarrow}\cone(\HH(R_{\reg})\overset{1-\bar{\sigma}_0}{\lra}\HH(R)).    
\end{gather*}
Summing up, we obtain, for all $m\in\Z$, a natural zig-zag of quasi-isomorphisms
\begin{multline*}
\cone (\HH(R_{\reg}/\ell^{\reg(E)},\bar{S}_m)\overset{1-\bar{\sigma}_m}\lra \HH(R/\ell^{E^0},\bar{S}_m))\weq
{}_m\HH(L(G,E,\phi)).    
\end{multline*}

\noindent{Step 2: \em $E$ finite.} One can get from any finite graph $E$ to another finite graph $E'\subset E$ such that any sources of $E'$ are also sinks,  through iterations of the source elimination move $E\mapsto E_{\setminus v}$ described in
\cite{lpabook}*{Definition 6.3.26}. The algebra $L'=L(G,E',\phi_c)$ embeds into $L$ as 
the corner associated to the homogeneous idempotent $1_{E'}=\sum_{v\in E'}v$, which is 
a full idempotent \cite{fas13}*{Proposici\'on 6.11} (see also 
\cite{flow}*{Proposition 1.14}). Hence the inclusion $L'\subset L$ induces a grading-
preserving quasi-isomorphism $\HH(L(E')/\ell^{E'})\to \HH(L(E)/\ell^{E})$. Remark that 
the source elmination process may eliminate vertices which are not sources of the 
original graph, but become ones after iterating the process. However those vertices 
that lie in a closed path of the original graph remain untouched. Hence, by 
Remark \ref{rem:smclosed}, for $R'=\bigoplus_{v\in E'}R_v$, we 
have $\HH(R'/\ell^{E'^0},\bar{S}_m(E'))=\HH(R/\ell^{{E}^0},\bar{S}_m(E))$ for all $m\ne 0$. It remains to 
show that if $v\in\sour(E)\setminus\sink(E)$, then for $F=E_{\setminus v}$ and 
$R"=\bigoplus_{v\in F}R_v$, the inclusion 
is a quasi-isomorphism
\begin{multline*}
    \cone(\HH(R"_{\reg(F)}/\ell^{\reg(F)})\overset{1-\bar{\sigma}_0}{\lra}\HH(R"/\ell^{F^0}))
    \\
    \weq \cone(\HH(R_{\reg(E)}/\ell^{\reg(E)})\overset{1-\bar{\sigma}_0}{\lra}\HH(R/\ell^{E^0})).
\end{multline*}
In fact the map above is injective, and its cokernel is the cone of the identity map of 
$\HH(R_v)$, which is contractible.

\noindent{Step 3: \em $E$ row-finite.} This case follows from Lemma \ref{lem:EPcoli} and the fact that the Hochschild complex commutes with filtering colimits. 
\medskip

Part 2): proof of \eqref{map:conGconR}.
\medskip
For all $m\in\Z$ we have a commutative diagram with vertical surjections
\[
\xymatrix{\HH(k[G]_{\reg(E)}/\ell^{(\reg(E))},S_m)\ar[r]^(.55){1-\sigma_m}\ar[d]&
\HH(k[G]_{E^0}/\ell^{(E^0)},S_m)\ar[d]\\
\HH(R_{\reg}/\ell^{(\reg(E))},\bar{S}_m)\ar[r]^{1-\bar{\sigma}_m}&\HH(R/\ell^{(E^0)},\bar{S}_m)}
\]
The kernel $K^m$ of both vertical maps is the same and is spanned in dimension $n$ 
by the elementary tensors 
$x_0\otimes\cdots\otimes x_{n-1}\otimes \alpha\otimes x_n$ if $m\ge 0$ and $x_0\otimes\cdots\otimes x_{n}\otimes \alpha^*$ if $m<0$, with at least one 
$x_i\in I=\bigoplus_{v\in \reg(E)}I_v$ (recall $I_v=0$ if $v\in\sink(E)$). In particular 
$\sigma_m$ resticts to an endomorphism of $K^m$. We shall show that this 
endomorphism is locally nilpotent, and thus that $1-\sigma_m:K^m\to K^m$ is an 
isomorphism, from which \eqref{map:conGconR} will follow.

As was made clear in the proof of the first part, $\bar{\sigma}_0$ is the composite of the trace map and the chain homomorphism induced by the inclusion $\inc:L_{0,0}\subset L_{0,1}$. It is also clear that its lift $\sigma_0$ factors as $\tr\circ \jmath_0$. Recall from Remark \ref{rem:smclosed} that
$\HH(R_{\reg}/\ell^{\reg(E)},\bar{S}_m)\cong\HH(R/\ell^{(E^0)}, \bar{S}_m)$ and similarly with $k[G]$ and $S_m$ substituted for $R$ and $\bar{S}_m$. The morphism $\inc$ also induces a chain map $\inc^{\reg}:\HH(R_{\reg}/k^{(\reg(E))},\bar{S}_m)\to \HH(L_{0,1}^{\reg}/k^{(\reg(E))},P^{\reg}_1\otimes_{R_{\reg}}\bar{S}_m\otimes_{R_{\reg}}Q^{\reg}_1)$, and again $\bar{\sigma}_m=\tr\circ\jmath_0$. Similarly, for $n\in\N_0$
\[
\hat{P}_n=\bigoplus_{v\in\reg(E)}k^{(\cP_{n,v})}\otimes k[G] \text{ and }\hat{Q}_n=\bigoplus_{v\in\reg(E)}k[G]\otimes k^{(\cP_{n,v}^*)}
\]
and $\cM^{\reg}_n=\bigoplus_{v\in\reg(E)}M_{\cP_{n,v}}k[G]$, $\jmath_n$ induces a chain map 
\begin{multline*}
\jmath^{\reg}_n:\HH(\cM_n^{\reg}/\ell^{(\reg(E))},\hat{P}_n\otimes_{k[G]_{\reg(E)}}S^{\reg}_m\otimes_{k[G]_{\reg(E)}}\hat{Q}_n)\to \\
\HH(\cM_{n+1}^{\reg}/\ell^{(\reg(E))},\hat{P}_{n+1}\otimes_{k[G]_{\reg(E)}}S^{\reg}_m\otimes_{k[G]_{\reg(E)}}\hat{Q}_{n+1})
\end{multline*}
and $\sigma_m=\tr\circ \jmath^{\reg}_0$. There are also trace maps from the homology of $\cM_{n+1}$ to that of $\cM_n$, and also between their $\reg$-summands, and we have
$\tr\circ \jmath_{n+1}=\jmath_n\circ\tr$, by naturality. Hence
$\sigma_m^n$ factors through $\jmath_{\le n}$. By Proposition \ref{prop:Iv} this implies that $\sigma_m$ is locally nilpotent on $K^m$, completing the proof. 
\end{proof}
 
\subsection{Twisted homology of Exel-Pardo groupoids}\label{subsec:hep}

Let $(G,E,\phi_c)$ be a twisted Exel-Pardo tuple; recall we write $\cG(G,E,\phi_c)$ for its tight groupoid, together with the associated groupoid cocycle $\bomega$ induced by $c$. In this section we abbreviate
\[
\cG=\cG(G,E,\phi_c).
\]

Let $L=L_k(G,E,\phi_c)$, 
\[
L\supset\cD=\mspanl\{\alpha\alpha^*\colon \alpha\in \cP(E)\}
\]
the diagonal $\ell$-subalgebra. Remark that the $k$-algebra isomorphism $L\iso \cA_k(\cG)$ 
mapping $\alpha g\beta^*\mapsto \chi_{[\alpha g\beta^*,Z_\beta]}$ of \cite{eptwist}*{Proposition 4.2.2} sends $\cD$ isomorphically 
onto $\cC_c(\cG^{(0)},\ell)$. Hence as explained in Section \ref{sec:hh-stein}
we have a monomorphism of chain complexes
\begin{equation}\label{map:hook}
\H(\cG,k/\ell)\hookrightarrow {}_0\HH(L/\cD).
\end{equation}
And if furthermore $\cG$ is Hausdorff then restriction of functions defines a chain map
\begin{equation}\label{map:res}
\res:{}_0\HH(L/\cD)\to \H(\cG,k/\ell)
\end{equation}
that is left inverse to \eqref{map:hook}. 

\begin{lem}\label{lem:gamepath}
Let $(G,E,\phi_c)$ be a twisted EP-tuple, $L=L(G,E,\phi_c)$, $n\ge 0$, $\alpha_0,\beta_0,\dots,\alpha_n,\beta_n\in \cP(E)$, $g_0,\dots,g_n\in G$, $a_0,\dots,a_n\in k$,
\[
\xi=\alpha_0a_0g_0\beta_0^*\otimes\cdots\otimes\alpha_ng_n a_n\beta^*_n\in L^{\otimes^{n+1}_\cD}
\]
and $\natural\xi$ its image in $(L^{\otimes^{n+1}_\cD})_\natural$
\item[i)] If $\xi\ne 0$, then there exist paths $\gamma_0,\cdots,\gamma_{n+1}\in \cP(E)$, $h_0,\dots,h_n\in G$ and $b_0,\dots,b_n\in k$ such that
\begin{equation}\label{eq:nuevoxi}
\xi=\gamma_0b_0h_0\gamma_1^*\otimes\gamma_1 b_1 h_1\gamma_2^*\otimes\cdots\gamma_nb_nh_n\gamma_{n+1}^*.  
\end{equation}

\item[ii)] If $0\ne \natural\xi$ has total degree $|\xi|=0$, then there are paths $\mu_0,\dots,\mu_n$, $f_0,\dots,f_n\in G$ and $c_0,\dots,c_n\in k$ such that 
\[
\natural \xi=\natural \mu_0c_0f_0\mu_1^*\otimes\cdots\otimes \mu_nc_nf_n\mu_0^*.
\]
\end{lem}
\begin{proof}
We begin by noticing that if $\alpha,\beta,\gamma$ and $\delta$ are paths in $E$, $g,h\in G$ and $a,b\in k$, then for $\otimes=\otimes_{\cD}$ we have 
\begin{gather}\label{eq:pasabeta}
\theta:=\alpha ag\beta^*\otimes\delta bh\gamma^*=\alpha ag\beta^*\delta\delta^*\otimes\delta bh\gamma^*\\
=\alpha a g\beta^*\otimes\beta\beta^*\delta bh\gamma^*.\nonumber
\end{gather}
Thus if $\theta\ne 0$  we must either have $\beta=\delta\beta'$ or $\delta=\beta\delta'$. In both cases we can use the identities \eqref{eq:pasabeta} to rewrite $\xi$ as a tensor in which both middle paths coincide. Indeed, in the first case we have
\begin{gather*}
\theta=\alpha ag\beta^*\otimes\beta(\beta')^*bh\gamma^*\\
=\alpha ag\beta^*\otimes \beta (c(h,h^{-1}(\beta'))b)\phi(h,h^{-1}(\beta'))(\gamma h^{-1}(\beta'))^*,
\end{gather*}
and in the second
\begin{gather*}
\theta=(\alpha g(\delta'))(ac(g,\delta'))\phi(g,\delta')\delta^*\otimes\delta bh\gamma^*.
\end{gather*}
Hence in the situation of i) the fact that $\xi\ne 0$ implies that $\beta_i$ and $\alpha_{i+1}$ are comparable for all 
$0\le i\le n-1$, and we shall show how one can use the procedure above to rewrite $\xi$ as in \eqref{eq:nuevoxi}. As a 
first step, we compare $\beta_0$ and $\alpha_1$; if they are equal, we pass to the second step. Otherwise we use the 
procedure above to replace either $\beta_0$ or $\alpha_1$ by whichever of them is longer (i.e. has higher length), and 
modify either $\beta_0$ or $\alpha_1$ and their accompanying coefficients accordingly. In the second step we repeat the procedure at the second $\otimes$; 
if $|\beta_1|\ge|\alpha_2|$ we proceed exactly as before and pass over to the next $\otimes$. If instead $|\alpha_2|>|
\beta_1|$ the above procedure will make us modify again the newly acquired $\alpha_1$, replacing it by a longer path, 
which will in turn force us to change $\alpha_0$, so that in the new rewriting of $\xi$, $\beta_0=\alpha_1$ and $
\beta_1=\alpha_2$. Following in this way, after at most $n(n+1)/2$ steps we end up with an elementary tensor where all 
the $\beta_i=\alpha_{i+1}$ for $0\le i\le n-1$. This proves i). In the situation of ii), the hypothesis that $\natural 
\xi\ne 0$ implies that $\gamma_0$ and $\gamma_{n+1}$ are comparable; since in addition we are assuming that $|\xi|=0$, they must have the same length. All this together implies that $\gamma_0=\gamma_{n+1}$, completing the proof. 
\end{proof}

Let $T=(G,E,\phi)$ be an EP-tuple, $g\in G$, and $\gamma\in\cP(E)$. Recall that $g$ fixes $\gamma$ strongly if $g(\gamma)=\gamma$ and $\phi(g,\gamma)=1$. For example, every path is strongly fixed by the trivial element $1\in G$. The triple $T$ is called \emph{pseudo-free} if $1$ is the only element of $G$ that fixes a path strongly. In other words, $T$ is pseudo-free whenever
\[
g(\gamma)=\gamma\text{ and } \phi(g,\gamma)=1\Rightarrow g=1.
\]

\begin{rem}\label{rem:e*unitary} It was shown in \cite{ep}*{Proposition 5.8} that an EP-triple $(G,E,\phi)$ is pseudo-free if and only if $\cS=\cS(G,E,\phi)$ is \emph{$E^*$-unitary}. This means that if $s,p\in \cS$ and $p^2=p\ne 0$, then $sp=p$ implies $s^2=s$. 
\end{rem}

\begin{lem}\label{lem:pseudo}
Let $(G,E,\phi_c)$ be a twisted Exel-Pardo tuple  such that $\cG$ is Hausdorff. Let 
$(\alpha_0,\dots,\alpha_n)$ be a tuple of paths in $E$ such that 
$r(\alpha_i)=r(\alpha_{i+1})$ for all $i$. Also let $g_0,\dots,g_n\in G$ and $a_0,
\dots,a_n\in k$ such that $0\ne a_0\otimes\dots\otimes a_n\in k^{\otimesl n+1}$. 
Consider the element
\begin{equation}\label{eq:elxi}
\xi=a_0\alpha_0g_0\alpha_1^*\otimes a_1\alpha_1g_1\alpha_2^*\otimes\cdots\otimes a_n\alpha_ng_n\alpha_0^*\in {}_0\HH(L/\cD)_n.
\end{equation}
\item[a)] If $g_0\cdots g_n=1$, then $\xi\in \H(\cG,k/\ell)$. 
\item[b)] The following are equivalent.
\item[i)] $(G,E,\phi)$ is pseudo-free.
\item[ii)] $\res(\xi)\ne 0$ implies that $g_0\cdots g_n=1$.
\item[c)] The elements $\xi$ as above such that $g_0\cdots g_n=1$ generate $\H(\cG,k/\ell)$ as an abelian group. 
\end{lem}
\begin{proof}
The function $f:\cG_{\cyc}^n\to k^{\otimesl n+1}$ corresponding to $\xi$ is supported on the following subset of $\cG_{\cyc}^n$
\begin{equation}\label{eq:suppxi}
\{([\alpha_0g_0\alpha_1^*,\alpha_1 g_1\cdots g_n(\theta)],\dots,  [\alpha_{n-1}g_{n-1}\alpha_{n}^*, \alpha_ng_n(\theta)] , [\alpha_ng_n\alpha_0^*,\alpha_0\theta])\colon g_0\cdots g_n(\theta)=\theta\}.
\end{equation}
The product of the coordinates of the element above is
\[
\eta=[\alpha_0g_0\cdots g_n\alpha_0^*,\alpha_0\theta].
\]
The element $\eta$ is in $\cG^0$ if and only if there is a finite path $\gamma\ge\theta$ such that for $g=g_0\cdots g_n$, the following identity holds in $S(G,E,\phi)$ 
\begin{equation}\label{eq:eunitary}
\alpha_0\gamma\gamma^*\alpha_0^*=\alpha_0g\alpha_0^*\alpha_0\gamma\gamma^*\alpha_0^*.
\end{equation}
Left-multiplying by $\alpha_0^*$ and right multiplying by $\alpha_0\gamma$ and using that $s(\gamma)=r(\alpha)$, we get
\begin{equation}\label{eq:strofix}
\gamma=g\cdot \gamma 
\end{equation}
which implies that $\gamma$ is strongly fixed by $g$. Conversely, left multiplying \eqref{eq:strofix} by $\alpha_0$ and right-multiplying it by $(\alpha_0\gamma)^*$ recovers \eqref{eq:eunitary}. If $g=1$, \eqref{eq:strofix} holds for 
$\gamma=r(\alpha_0)$, proving a). The converse holds if and only if there are no strongly fixed paths, which implies that 
$(G,E,\phi)$ is pseudo-free. This proves b). Next consider the subgroup $M$ spanned by the elements of c). In view of a), $M\subset \H(\cG,k/\ell)$. To prove the other inclusion it suffices to show that for a general element \eqref{eq:elxi}, 
$\res(\xi)\in M$. As above, we set $g=g_0\cdots g_n$; we may assume $g\ne 1$. Now $\res(\xi)$ is the constant function $a_0\otimes\cdots\otimes a_n$ on $Y=\supp(\xi)\cap \Gamma(\cG_{\cyc},\cG^0)$, which, by what we have just seen, consists of those elements 
\begin{equation}\label{eq:ygcyc}
y=([\alpha_0g_0\alpha_1^*,\alpha_1 g_1\cdots g_n(\theta)],\dots,  [\alpha_{n-1}g_{n-1}\alpha_{n}^*, \alpha_ng_n(\theta)])\in Y
\end{equation}
for which there is  $\gamma\ge \theta$ with $\gamma$ strongly fixed by $g$. Because $\cG$ is Hausdorff by assumption, there are finitely many paths, say $\gamma_1,\dots,\gamma_l$, starting at $r(\alpha_0)$ and which are minimal among those strongly fixed by $g$ \cite{ep}*{Theorem 12.2}. Hence we can write $Y=\bigsqcup_{i=1}^l Y_i$ where $Y_i$ consists of those elements $y$ of the form \eqref{eq:ygcyc}
 with $\gamma_i\ge \theta$. Then if $y\in Y_i$ we can write $\theta=\gamma_i\theta_i$, and for  
 \[
 x_j=\alpha_{j+1}(g_{j+1}\cdots g_n)(\gamma_i\theta_i)
 \]
 we have
\begin{gather*}
    y_j=[\alpha_jg_j\alpha^*_{j+1},x_j]\\
    =[\alpha_jg_j\alpha^*_{j+1}(\alpha_{j+1}(g_{j+1}\cdots g_n)(\gamma_i))(\alpha_{j+1}(g_{j+1}\cdots g_n)(\gamma_i))^*,x_j]\\
    =[\alpha_j(g_j\cdots g_n)(\gamma_i)\phi(g_j,(g_{j+1}\cdots g_n)(\gamma_i))(\alpha_{j+1}(g_{j+1}\cdots g_n)(\gamma_i))^*,x_j ].    
\end{gather*}
For $0\le j\le n$, put
\[
\xi_{i,j}=a_j\alpha_j(g_j\cdots g_n)(\gamma_i)\phi(g_j,(g_{j+1}\cdots g_n)(\gamma_i))(\alpha_{j+1}(g_{j+1}\cdots g_n)(\gamma_i))^*.
\]
Consider the element 
\[
\xi_i=\xi_{i,0}\otimes\cdots\otimes \xi_{i,n}\in \HH(L(G,E,\phi_c)).
\]
Then $\xi_i$ is supported at $Y_i$ where it is constantly equal to $a_0\otimes\cdots\otimes a_n$; thus $\xi_i\in \H(\cG,\bomega,k/\ell)$ and $\res(\xi)=\sum_{i=1}^l \xi_i$. Moreover, using the cocycle condition and the fact that  $g_0\cdots g_n$ fixes $\gamma_i$ strongly, we obtain
\begin{gather*}
\phi(g_0,(g_1\cdots g_n)(\gamma_i))\phi(g_1,(g_2\cdots g_n)(\gamma_i))\cdots \phi(g_n,\gamma_i)\\
=\phi(g_0g_1,(g_2\cdots g_n)(\gamma_i))\phi(g_2,(g_3\cdots g_n)(\gamma_i))\cdots \phi(g_n,\gamma_i)\\
=\phi(g_0\cdots g_n,\gamma_i)=1.
\end{gather*}
\end{proof}

\begin{thm}\label{thm:hep}
Let $(G,E,\phi_c)$ be a twisted EP-tuple where $E$ is row-finite and $G$ acts trivially on $E^0$. Assume that the underlying untwisted EP-tuple $(G,E,\phi)$ is pseudo-free. 
Let $A=A_E\in \Z^{(\reg(E)\times E^0)}$ be the reduced adjacency matrix. For $v,w\in E^0$, let $\tau\in \Hom(\H(G,k/\ell)\otimesl\ell^{(\reg(E))},\H(G,k/\ell)\otimesl\ell^{(E^0)})$ be the matrix of chain homomorphisms with entries
 \begin{gather*}\label{map:tau}
 \tau_{v,w}:\H(G,k/\ell)_n\to\H(G,k/\ell)_n,\\
 \tau_{v,w}(a)=A_{w,v}a,\, \text{and for } n\ge 1,\nonumber\\
 \tau_{v,w}(a_0\otimes g_1a_1\otimes\cdots\otimes g_na_n)=\nonumber \\
 \sum_{s(e)=v\, r(e)=w}a_0c(g_1\cdots g_n,e)^{-1}\otimes\phi_c(g_1,g_2\cdots g_n(e))a_1\otimes\cdots\otimes \phi_c(g_n,e)a_n.\nonumber
 \end{gather*}
 Recall that $\cG$ is the tight EP-groupoid $\cG(G,E,\phi)$ equipped with the groupoid $2$-cocyle
 $\bomega:\cG\times\cG\to \cU(k)$ induced by $c$. Let $k/\ell$ be a flat ring extension and 
 $\H(\cG,k/\ell)$ the complex for relative twisted groupoid homology.  
 Then there is a natural zig-zag of quasi-isomorphisms
 \[
 \cone(I-\tau)\weq \H(\cG,k/\ell).
 \]
\end{thm}

\begin{proof}
Let $L=L(G,E,\phi_c)$ and let $\cD$ 
be the diagonal $\ell$-subalgebra. By part c) of Lemma \ref{lem:pseudo}, $\H(\cG,k/\ell)$ is the subcomplex of 
${}_0\HH(L/\cD)$ given in degree $n$ by
 \[
 \mspanl\{a_0\alpha_0(g_1\cdots g_n)^{-1}\alpha_1^*\otimes a_1\alpha_1g_1\alpha_2^*\otimes\cdots\otimes 
 a_n\alpha_ng_n\alpha_0^*\}.
 \]
 One checks that the map 
 \begin{gather*}
 \jmath:\H(G,k/\ell)\to \HH(k[G]/\ell),\\ \jmath(a_0\otimes a_1g_1\otimes\cdots\otimes a_ng_n)=a_0(g_1\cdots g_n)^{-1}\otimes a_1g_1\otimes\cdots\otimes a_ng_n
 \end{gather*}
 fits into a commutative diagram as follows, where the composite of the vertical maps is the identity
\[
\xymatrix{
\ell^{\reg(E)}\otimesl \H(G,k/\ell)\ar[d]^{1\otimes\jmath}\ar[r]^{I-\tau}& \ell^{E^0}\otimesl \H(G,k/\ell)\ar[d]^{1\otimes\jmath}\\
\ell^{\reg(E)}\otimesl\HH(k[G]/\ell)\ar[d]^{1\otimes\res}\ar[r]^{I-\sigma_0}& \ell^{E^0}\otimesl\HH(k[G]/\ell)\ar[d]^{1\otimes\res}\\
\ell^{\reg(E)}\otimesl \H(G,k/\ell)\ar[r]^{I-\tau}& \ell^{E^0}\otimesl \H(G,k/\ell)
}
\]
In particular the cone of $I-\tau$ is a direct summand of the cone of $I-\sigma_0$. 
Because we are assuming that $(G,E,\phi)$ is pseudo-free, $\cG$ is Hausdorff, so the map \eqref{map:res} is defined and thus $\H(\cG,k/\ell)$ is a 
direct summand of ${}_0\HH(L/\cD)$. We will show that the zigzag of quasi-isomorphisms of 
Theorem \ref{thm:hhep} induces one between these two direct summands. We start by 
considering the case when $E$ is finite without sources. It follows from the explicit 
formula of Remark \ref{rem:elkappa} that the map 
$$
\theta:\cone(\HH(L_0,L)\overset{1-\psi}{\lra} \HH(L_0,L))\to \HH(L)
$$ 
descends to a map
$$
\theta:\cone(\HH(L_0/\cD)\overset{1-\psi}{\lra}\HH(L_0/\cD))\to \HH(L/\cD)
$$ 
and that the restriction of the latter to $\cone(I-\bar{\sigma}_0)$ composed with the projection $\pi:\cone(I-\sigma_0)\onto\cone(I-\tau)$ fits into a commutative diagram

\begin{equation}\label{map:contauhomo}
\xymatrix{\cone(I-\tau)\ar[r]^{\pi\circ\theta}\ar[d]^{1\otimes\jmath}& \H(\cG,k/\ell)\ar[d]^{\inc}\\
\cone(I-\sigma_0)\ar[r]^{\pi\circ\theta}&\HH(L_0/\cD)}
\end{equation}
Using Lemma \ref{lem:pseudo} again and pseudo-freeness, we obtain that $\res\circ\theta\circ\pi=\theta\circ\pi\circ\res$. 
Hence $\theta\circ\pi$ is a retract of a quasi-isomorphism and therefore a quasi-isomorphism. Next assume that $E$ is finite, let 
$v\in\sour(E)\setminus\sink(E)$ and $F=E_{\setminus v}$ the source elimination graph. Set 
$L'=L(G,F,\phi_c)$, $\cG'=\cG(G,F,\phi_c)$ 
and $\cD'\subset L'$ the diagonal $\ell$-subalgebra. Then $\HH(L'/\cD')$ is a subcomplex of 
$\HH(L/\cD)$ that restricts to an inclusion between the twisted homology complexes relative 
to $k/\ell$, and is compatible with restriction maps. Hence the inclusion is a quasi-
isomorphism $\H(\cG',k/\ell)\weq \H(\cG,k/\ell)$. Let $\tau'=\tau_F:
\H(G,k/\ell)\otimes\ell^{\reg(F)}\to \H(G,k/\ell)\otimes\ell^{F^0}$. Then $\coker(\cone(1-
\tau')\to \cone(1-\tau))$ is the cone of an identity morphism, and so the inclusion $
\cone(1-\tau')\subset\cone(1-\tau)$ is a quasi-isomorphism. This proves the theorem for all 
twisted EP-tuples with finite underlying graph. The general case, for twisted EP-tuples 
over row-finite graphs, follows from Lemma \ref{lem:EPcoli} and the fact that homology 
commutes with filtering colimits. 
\end{proof}

\begin{coro}\label{coro:h0hep}
 Let $(G,E,\phi_c)$ and $\cG$ be as in Theorem \ref{thm:hep}. Write $\otimes=\otimes_\Z$. Then 
 \[
 H_0(\cG,k/\ell)=\BF(E)\otimes k.
 \]    
 In particular, the group above depends only on $k$ and $E$ and not on the flat extension $k/\ell$.
\end{coro}

\begin{coro}\label{coro:hep}
Let $(G,E,\phi_c)$ be as in Theorem \ref{thm:hep}, $\cG=\cG(G,E,\phi_c)$ and $\cG_u=\cG_u(G,E,\phi_c)$ the tight and the universal groupoid of $\cS(G,E,\phi)$, equipped with the groupoid cocycles induced by $c$. Also let $U$ and $\cG'={\cG_u}_{|U}$ be as in Lemma \ref{lem:cohn=leav}. Consider the chain maps $\iota:\H(\cG',k/\ell)\to\H(\cG_u,k/\ell)$ and $p:\H(\cG_u,k/\ell)\to\H(\cG,k/\ell)$ induced by the inclusion $U\subset \hat{\fX}(E)$ and the restriction map. Then there is an isomorphism of triangles in the derived category of chain complexes of $k$-modules
\begin{equation}\label{map:triang}
\xymatrix{\H(G,k/\ell)^{(\reg(E))}\ar[d]^{\wr}\ar[r]^{I-\tau}&\H(G,k/\ell)^{(E^0)}\ar[d]^{\wr}\ar[r]&\H(\cG,k/\ell)\ar[d]^{=}\\
\H(\cG',k/\ell)\ar[r]^\iota&\H(\cG_u, k/\ell)\ar[r]^p&\H(\cG,k/\ell).}
\end{equation}
\end{coro}
\begin{proof}
The inclusion $G\times\reg(E)\to\cG'$, $(g,v)\mapsto [vg,v]$ is a homomorphism of discrete groupoids, and the induced algebra homomorphism 
$$k[G]^{(\reg(E))}\to \cK(G,E,\phi_c)\cong\bigoplus_{v\in\reg(E)}M_{\cP_v}k[G]$$ is the full corner embedding $vg\mapsto \epsilon_{v,v}g$. By Morita invariance, the latter embedding induces a quasi-isomorphism $\HH(k[G],k/\ell)^{(\reg(E))}\weq \bigoplus_{v\in\reg(E)}\HH(M_{\cP_v}k[G])$ which one checks commutes with the inclusion and restriction maps to and from the respective groupoid homology complexes. Hence it restricts to a quasi-isomorphism between the latter complexes; this is the first vertical map of \eqref{map:triang}. By Lemma \ref{lem:cohn=leav}, $\cG_u=\cG(G,\tilde{E},\phi_c)$. By Theorem \ref{thm:hep}, $\H(\cG_u,k/\ell)$ is quasi-isomorphic to the cone of 
\begin{equation}\label{map:conocohn}
\xymatrix{\H(G,k/\ell)^{(\reg(E))}\ar[rr]^(.4){\begin{bmatrix}1-\tau^{\reg}\\ -\tau^{\reg'}\\ -\tau^{\sink} \end{bmatrix}}&& \H(G,k/\ell)^{(\reg(E)\sqcup\reg(E)'\sqcup\sink(E))}}.
\end{equation}
The projection 
\[
\H(G,k/\ell)^{(\reg(E)\sqcup\reg(E)'\sqcup\sink(E))}\to \H(G,k/\ell)^{(\reg(E))},\, (x_1,x_2',y)\mapsto x_1-x_2
\]
defines a surjection $\pi$ from the cone of \eqref{map:conocohn} onto the cone of the identity. Hence the cone of 
\eqref{map:conocohn} is equivalent to $\ker(\pi)\cong\H(G,k/\ell)^{(E^0)}$. Thus we obtain a quasi-isomorphism
$\H(G,k/\ell)^{(E^0)}\weq \H(\cG_u,k/\ell)$; this is the vertical map in the middle of \eqref{map:triang}. Next we check 
commutativity of the left square; that of the right square is clear. Let $v\in\reg(E)$. An elementary tensor 
$\xi:a_0v\otimes a_1g_1v\otimes\cdots\otimes a_n g_nv\in \H(G,k/\ell)$ goes in $\H(\cG')$ to the elementary tensor $\xi'$ that is 
obtained upon replacing $g_iv$ by $\chi_{[g_iv,v]}$ everywhere. Under the isomorphism $\cA(\cG_u)\cong C:=C(G,E,\phi_c)$ of Lemma 
\ref{lem:cohn=leav}, $\xi'$ is mapped to the elementary tensor 
$$
q\xi:=a_0(g_1\cdots g_n)^{-1}q_v\otimes a_1g_1q_v\otimes\cdots\otimes a_ng_nq_v\in \H(\cG_u,k/\ell).
$$ 
Put $g_0=(g_1\cdots g_n)^{-1}$. For each subset $A\subset[n]=\{0,\dots,n\}$, let $\xi(A)_i=a_ig_iv$ if $v\notin A$ and 
$\xi(A)_i=a_igm_v=a_ig\sum_{s(e)=v}ee^*=\sum_{s(e)=v}a_ic(g,e)g(e)\phi(g,e)e^*$ if $i\in A$. Set $\xi(A)=\xi(A)_0\otimes\cdots\otimes \xi(A)_n$. Remark that $\xi(\emptyset)=\xi$; apply \eqref{eq:pasabeta} repeatedly to obtain that for $A\ne\emptyset$
\[
\xi(A)=\sum_{e\in\cP^v_1}a_0e\phi_c(g_0,g_1\cdots g_n(e))(g_1\cdots g_n)(e)^*\otimes\cdots\otimes a_ng_n(e)\phi_c(g_n,e)e^*.
\]
Now use that $q_v=v-m_v$ and bilinearity 
of $\otimes_{\cD}$ to obtain 
\begin{gather*}
\tr(q\xi)=\sum_{A\subset [n]}(-1)^{|A|}\tr(\xi(A))=\xi+\\
\sum_{\emptyset\ne A\subset [n]}(-1)^{|A|}\sum_{e\in\cP^v_1}\tr(a_0e\phi_c(g_0,g_1\cdots g_n)(g_1\cdots g_n)(e)^*\otimes\cdots\otimes a_ng_n(e)\phi_c(g_n,e)e^*)\\
=\xi+(\sum_{i=1}^{n+1}(-1)^{i}\binom{n+1}{i})\sum_{e\in\cP^v_1}a_0\phi_c(g_0,g_1\cdots g_n)r(e)\otimes\cdots\otimes a_n\phi_c(g_n,e)r(e)\\
=\xi-\tau(\xi).
\end{gather*}
\end{proof}
\begin{rem}\label{rem:ortega}
In \cite{homology-katsura}, Eduard Ortega computes the integral homology of Katsura groupoids $\cG_{A,B}$ associated to a pair of square matrices. Since the latter are Exel-Pardo groupoids, Theorem \ref{thm:hep} also computes $H_*(\cG_{A,B},\Z)$, recovering Ortega's result in the pseudofree case.
\end{rem}

\subsection{\topdf{$K$}{K}-theory of twisted Exel-Pardo algebras}\label{subsec:kep}

Let $\ell$ be a commutative, unital ring. Let $\cT$ be category and $\cH:\ahal\to\cT$ a functor. We say that $\cH$ is \emph{homotopy invariant} if for every $A\in \ahal$, $\cH$ sends the inclusion $A\subset A[t]$ to  an isomorphism $\cH(A)\cong \cH(A[t])$. Let $X$ be an infinite set and $x\in X$; we say that $\cH$ is \emph{$M_X$-stable} if for every $A\in \ahal$, $\cH$ sends the corner inclusion $A\to M_XA$, $a\mapsto \epsilon_{x,x}a$ to an isomorphism. $M_X$-stability turns out to be independent of the choice of the element $x\in X$ \cite{kkh}*{Lemma 2.4.1}. We say that $\cH$ is \emph{excisive}
if $\cT$ is triangulated and every algebra extension 
\[
(\cE)\,\,\, 0\to A\to B\to C\to 0
\]
is mapped to a distinguished triangle
\[
\cH(C)[1]\overset{\partial_{\cE}}{\lra}\cH(A)\to \cH(B)\to \cH(C)
\]
where $[1]$ is the inverse suspension and the $\partial_{\cE}$ satisfy certain naturality conditions, as detailed in \cite{kk}*{Section 6.6}.
Let $I$ be a set and $\cT$ an additive category. We say that $\cH$ is \emph{$I$-additive} if first of all direct sums of cardinality $\le \sharp I$ exist in $\cT$ and second of all the map
\[
\bigoplus_{j\in J}\cH(A_j)\to \cH(\bigoplus_{j\in J}A_j)
\]
is an isomorphism for any family of algebras $\{A_j:j\in J\}\subset\ahal$ with $\sharp J\le\sharp I$. 
Now let $E$ be a graph and $\cT$ a triangulated category. We say that a functor $\H:\ahal\to\cT$ is \emph{$E$-stable} if it is $M_X$-stable with respect to a set $X$ of cardinality $\sharp(E^0\coprod E^1\coprod\N)$. 

Let $k$ be a commutative unital $\ell$-algebra and let $j:\ahal\to kk$ be the universal homotopy invariant, $E$-stable and excisive functor $j:\ahal\to kk$ constructed in \cite{kk}. Let $(v,w)\in \reg(E)\times E^0$ be such that $vE^1w\ne\emptyset$. Consider the homomorphism of algebras
\[
\jmath_{v,w}:k[G]\to \cM_{vE^1w}k[G], \, \jmath_{v,w}(g)=
\sum_{e\in vE^1w}\epsilon_{g(e),e}\phi_c(g,e).
\]
For any choice of $e\in vE^1w$, the homomomorphism $\inc_e:k[G]\to M_{vE^1w}k[G]$, 
$g\mapsto \epsilon_{e,e}k[G]$ yields the same $kk$-isomorphism $\epsilon:=j(\inc_e)$. 
Put
\begin{gather}
\Phi\in kk(k[G],k[G])^{\reg(E)\times E^0},\nonumber \\
\Phi_{v,w}=\epsilon^{-1}\circ j(\jmath_{v,w})\in kk(k[G],k[G]).\label{map:elPhi}
\end{gather} 
Let $\Phi^t\in kk(k[G],k[G])^{E^0\times \reg(E)}$ be the transpose of $\Phi$. 
If $\cH:\ahal\to \cT$ is homotopy invariant, $E$-stable and excisive, then by universal
property, we have $\cH=\bar{\cH}\circ j$ for some triangle functor 
$\bar{\cH}:kk\to \cT$; we shall abuse notation and write
$\cH(\Phi)^t=\bar{\cH}(\Phi)^t\in \cT(k[G],k[G])^{E^0\times\reg(E)}$. If in addition  
$\cH$ is $E^0$-additive, then by row-finiteness of $E$, $\cH(\Phi)^t$ defines a 
homomorphism in $\cT$
\[
\cH(\Phi)^t:\cH(k[G])^{(\reg(E))}\to \cH(k[G])^{(E^0)}.
\]

In particular this happens when $\cH=j$ and $E$ is finite.

Finally let $\cM$ be a stable simplicial model category, $\cT=\Ho\cM$ the homotopy category and $[\,]:\cM\to \cT$ the localization functor. We say that a functor $H:\ahal\to\cM$ is \emph{finitary} if the canonical map
$\hoco_n H(A_n)\to H(\colim_nA_n)$ is a weak equivalence for every inductive system of algebras $\{A_n\to A_{n+1}:n\in\N\}$. We say that a functor $\cH:\ahal\to \cT$ is \emph{finitary} if there is a functor $H:\ahal\to\cM$ such that $\cH=[H]$ and such that $H$ is finitary.
\begin{nota}\label{nota:kk}
  For $A,B\in\ahal$ and $n\in\Z$, we write
$$
kk_n(A,B)=\hom_{kk}(j(A),j(B)[n]),\,\, kk(A,B)=kk_0(A,B).
$$
\end{nota}
\begin{ex}\label{ex:kh}
Weibel's homotopy algebraic $K$-theory \cite{kh} gives a functor $KH$ from $\ell$-algebras to the homotopy category of spectra, that is homotopy invariant, excisive, stable, additive, and finitary. Its homotopy groups can be expressed in terms of bivariant $K$-theory; 
we have $KH_n(A)=kk_n(\ell,A)$ for all $A\in\ahal$ and all $n\in\Z$ \cite{kk}*{Theorem 8.2.1}. There is a natural map of spectra $K(A)\to KH(A)$ which is $n+1$-connected whenever the map
\[
K_n(A)\to K_n(A[t_1,\dots, t_m]) 
\]
induced by the inclusion is an isomorphism for all $m$ \cite{kh}*{Proposition 1.5}. In this case we say that $A$ is \emph{$K_n$-regular}. By a theorem of Vorst \cite{vorst}*{Corollary 2.1(ii)}, $K_n$-regularity implies $K_{n-1}$-regularity. The ring $A$ is \emph{$K$-regular} if it is $K_n$-regular for all $n\in\Z$. Recall that a unital ring $R$ is called \emph{regular} if every (right) $R$-module has finite projective dimension; Noetherian regular rings are $K$-regular by \cite{quihik1}*{Corollary to Theorem 8}. This applies, for example, to $k[\Z^n]$ whenever $k$ is Noetherian regular. The ring $R$ is (right) \emph{coherent} if the category of finitely presented right $R$-modules is abelian, and \emph{supercoherent} if in addition $R[t_1,\dots, t_m]$ is coherent for every $m\ge 1$. Regular supercoherent rings are $K$-regular by \cite{gerfree}*{Theorem 3.1}. This applies, for example, to the group algebra $k[\F_n]$ of the free noncommutative group over a Noetherian regular ring $k$ \cite{gerfree}*{Theorem 1.13}. 
\end{ex}
\begin{thm}\label{thm:kep}
Let $(G,E,\phi_c)$ be a twisted EP-tuple with $E$ row-finite such that $G$ acts trivially on $E^0$. Let $\cT$ be a triangulated category and $\cH:\ahal\to\cT$ an excisive, homotopy invariant, $E$-stable and $E^0$-additive functor. Let $\Phi$ be as in \eqref{map:elPhi}. 
Then the Cohn extension of \eqref{ext:cohnext} induces the following distinguished triangle in $\cT$
\[
\cH(k[G])^{(\reg(E))}\overset{I-\cH(\Phi^t)}\lra \cH(k[G])^{(E^0)}\to \cH(L_k(G,E,\phi_c)).
\]
If furthermore $\cH$ is finitary, then we may substitute 
$\cH(R_{\reg})$ for $\cH(k[G])^{(\reg(E))}$ and $\cH(R)$ for $\cH(k[G])^{(E^0)}$ in the triangle above. 
\end{thm}    
\begin{proof}
Put $T=(G,E,\phi_c)$, $\cK=\cK(T)$, $C=C(T)$, $L=L(T)$. For $(v,w)\in \reg(E)\times E^0$ consider the following elements of $C$
\[
m_{v,w}=\sum_{s(e)=v,\, r(e)=w}ee^*,\, m_v=\sum_wm_{v,w},\, qv=v-m_v.
\]
Observe that if $v\in\reg(E)$, then $q_v=v-m_v\in\cK$ is the element of \eqref{eq:qv}, while if $v\in\sink(E)$, $m_v=0$ and $qv=v$. By \cite{eptwist}*{Proposition 6.2.5}, the algebra homomorphism $q:k[G]^{(\reg(E))}\to\cK$, $vg\mapsto gq_v$ is a $kk$-isomorphism and thus, by the additivity  hypothesis, it induces an isomorphism $\cH(k[G])^{(\reg(E))}\to\cH(\cK)$. 
By \cite{eptwist}*{Theorem 6.3.1}, the algebra inclusion $\iota:k[G]^{(E^0)}\to C$ is a $kk$-isomorphism too, and so induces an isomorphism 
$\cH(k[G])^{(E^0)}\to \cH(C)$,  again by additivity. Let $\hat{\cK}=\langle q_v\colon v\in E^0\rangle\triqui C$.  By \cite{eptwist}*{6.3.4}, the map
\begin{equation}\label{map:Ktilde}
\bigoplus_{w\in E^0}M_{\cP_w}k[G]\to \hat{\cK},\, \epsilon_{\alpha,\beta}wg\mapsto \alpha gq_w\beta^*    
\end{equation}
is an isomorphism of $k$-algebras. By the argument of \cite{eptwist}*{Proposition 6.2.5}, the map $\tilde{q}:k[G]^{(E^0)}\to \hat{\cK}$, $vg\mapsto gq_v$ is a $kk$-equivalence. 
The proof of \cite{eptwist}*{Theorem 6.3.1} considers the algebra homomorphism $\xi:C\to C$, $\xi(vg)=gm_v$, $\xi(e)=em_{r(e)}$, $\xi(e^*)=m_{r(e)}e^*$ and shows that the quasi-homomorphism $(\id,\xi):C\to C\truqui\cK$ followed by the inverse of $\tilde{q}$, is $kk$-inverse to $\iota$. A computation shows that
\begin{equation}\label{eq:xiqv}
\xi(gq_v)=\sum_{w}\sum_{s(e)=v,r(e)=w}g(e)\phi_c(g,e)q_we^*.
\end{equation}
If $v\in\reg(E)$, then under the isomorphism \eqref{map:Ktilde}, \eqref{eq:xiqv} corresponds
to the image of $\jmath_0(vg)$. Similarly the restriction of $\tilde{q}$ to $k[G]^{(\reg(E))}$ corresponds to 
a sum of corner inclusions. The first assertion of the theorem now follows by $E$-stability, additivity and 
excisivness of $\cH$. To prove the last assertion of the theorem,  we proceed as follows. Recall that  
$I=\bigoplus_{v\in\reg(E)}I_v=\ker(k[G]^{(E^0)}\to R)=\ker(k[G]^{(\reg(E))}\to R_{\reg})$. Hence it suffices  
to show that $I-\cH(\Phi^t)$ induces an isomorphism on $\cH(I)$. By Proposition \ref{prop:Iv},   
$I_v=\bigcup_nI(n)_v$ is an increasing union of ideals such that $\jmath_{v,w}(I_v(n))\subset M_{vE1w} I_w(n)$ 
and $\jmath_{\le n}$ vanishes on $I(n)=\bigcup_{v\in\reg(E)}I_v(n)$. It follows that $\Phi^t$ induces a 
nilpotent endomorphism of $\cH(I(n))$. Thus $I-\cH(\Phi)^t$
induces an automorphism of $\cH(I(n))$ for each $n$ whence $I-\cH(\Phi)^t:\cH(I)\to\cH(I)$ is an isomorphism, 
since $\cH$ is finitary.
\end{proof}

\begin{coro}\label{coro:kep}
Put $L=L_k(G,E,\phi_c)$. For $n\in\Z$ we have a long exact sequence
\[
KH_{n+1}(L)\to KH_n(k[G])^{(\reg(E))}\overset{I-\Phi^t}\lra KH_n(k[G])^{(E^0)}\to KH_n(L).
\]
If furthermore both $k[G]$ and $L$ are $K$-regular, then we may substitute $K$ for $KH$ in the sequence above.
\end{coro}
Let $(G,E,\phi_c)$ be a twisted EP-tuple. As before, we assume that $G$ acts trivially on $E^0$. Then for $(v,w)\in \reg(E)\times E^0$, each element $g\in G$ defines a permutation $\sigma_{v,w}(g)$ of the set $vE^1w$. For each $g\in G$, Consider the matrices $B(g), C(g)\in\cU(k[G])_{\ab}^{(\reg(E)\times E^0)}$,
\[
B_{v,w}(g)=\prod_{e\in vE^1w}\phi(g,e),\,\, C_{v,w}(g)=\sg(\sigma_{v,w}(g))\prod_{e\in vE^1w}c(g,e).
\]
Consider the matrix of homomorphisms
\[
D=\begin{bmatrix} A & 0\\ C & B\end{bmatrix}.
\]
Put
\[
D^t=\begin{bmatrix} A^t & C^t\\ 0 & B^t\end{bmatrix}.
\]
Observe that $D^t$ defines a group homomorphism 
\begin{equation}\label{map:D*}
D^t:\cU(k)^{(\reg(E))}\oplus G_{\ab}^{(\reg(E))}\to \cU(k)^{(E^0)}\oplus G_{\ab}^{(E^0)}    
\end{equation}
Recall from \cite{hanbu}*{Conjecture 1.11} that the \emph{Farrell-Jones conjecture} for the $K$-theory of the group algebra $k[G]$ of torsion free group over a regular Noetherian ring $k$ says that the assembly map
\[
BG\land \K(k)\to \K(k[G])
\]
is an equivalence. Here we abuse notation and write $BG$ for the suspension spectrum of the classifying space of $G$. There is a first quadrant spectral sequence
\[
H_p(G,K_q(k))\Rightarrow H_{p+q}(G,\K(k)).
\]
If, for example, $k$ is a field or a principal ideal domain, then $K_0(k)=\Z$ and $K_1(k)=\cU(k)$, and the conjecture implies that 
\begin{equation}\label{eq:k01neg}
 K_0(k[G])=\Z,\, K_1(k[G])=\cU(k)\oplus G_{\ab},\, K_n(k[G])=0\,\forall n<0, 
\end{equation}
and that there is a surjection
\begin{equation}\label{map:k2ontoh2}
K_2(k[G])\onto H_2(G,\Z).
\end{equation}
\begin{thm}\label{thm:k01ep}
Let $(G,E,\phi_c)$ be a twisted EP-tuple with $E$ row-finite, such that $G$ acts trivially on 
$E^0$. Let $k$ be a field or a PID. Assume that $G$ is torsionfree and satisfies the Farrell-Jones conjecture and that $L(G,E,\phi_c)$ is $K_1$-regular. Let $D^t$ be as in 
\eqref{map:D*}. Then
\item[i)] $K_0(L(G,E,\phi_c))=\BF(E)$.
\item[ii)] There is an exact sequence
\[
0\to \coker(I-D^t)\to K_1(L(G,E,\phi_c))\to \ker(I-A_E^t)\to 0.
\]    
\end{thm}
\begin{proof}
Put $L=L(G,E,\phi_c)$. Because by assumption $G$ is torsionfree and satisfies the $K$-theoretic Farrell-Jones conjecture, $k[G]$ is $K$-regular, so we may substitute $K_n(k[G])$ for $KH_n(k[G])$ in the sequence of Corollary \ref{coro:kep}, and the identities \eqref{eq:k01neg} hold. Since we are moreover assuming that $L$ is $K_1$-regular, we obtain an exact sequence
\begin{equation}
\xymatrix{\cU(k)^{(\reg(E))}\oplus G_{\ab}^{(\reg(E))}\ar[rr]^{I-\Phi^t}&&\cU(k)^{(E^0)}\oplus G_{\ab}^{(E^0)}\ar[r]&K_1(L)\ar[d]\\
 0& K_0(L)\ar[l]&\Z^{(E^0)}\ar[l] &\Z^{\reg(E)}\ar[l]_{I-\Phi^t}}    
\end{equation}
Next observe that if $a\in k$, and $(v,w)\in\reg(E)\times E^0$ then $\jmath_{v,w}(av)=(\sum_{e\in vE^1w}\epsilon_{e,e})aw$. In particular, $\Phi_{v,w}$ sends $[1]\in K_0(k)=\Z$ to $A_{v,w}$, and, if $a$ is invertible, 
\[
\Phi_{v,w}(a)=\det(\jmath_{v,w}(av))=a^{A_{v,w}}.
\] 
Similarly, 
\begin{gather*}
\jmath_{v,w}(vg)=\sigma_{v,w}(g)\circ\sum_{e\in vE^1w}\epsilon_{e,e}c(g,e)\phi(g,e),
\end{gather*}
and thus for the class $[g]\in G_{\ab}\subset K_1(k[G])$, we have
\begin{gather*}
\Phi_{v,w}([g])=\det(\jmath_{v,w}(vg))=(\sg(\sigma_{v,w}(g))\prod_{e\in vEw}c(g,e), [\prod_{e\in vE^1w}\phi(g,e)])\\
=(C_{v,w}(g), B_{v,w}([g])).
\end{gather*}
\end{proof}
Next we specialize to the case $G=\Z$. Denote $\Z$ multiplicatively and let $x$ be a generator, so that $\ell[\Z]=\ell[x,x^{-1}]$. Set $\sigma=(x-1)\ell[x,x^{-1}]$; we have $\ell[\Z]=\ell\oplus\sigma$. By \cite{kk}, 
$\sigma$ represents the suspension in $kk$. Hence writing $\sigma^i=\sigma^{\otimesl^i}$, we have
\[
kk(\sigma^i k, \sigma^j k)=KH_{i-j}(k).
\]
In particular, upon permuting summands, we may identify any element of
\goodbreak

\noindent $kk(k[\Z],k[\Z])^{E^0\times \reg(E)}$ with a matrix 
\[
\begin{bmatrix}
  X&Y\\ Z& W  
\end{bmatrix}
\]
where each of the blocks has size $E^0\times \reg(E)$, the coefficients of $X$ and $W$ are in $KH_0(k)$, and those of $Y$ and $Z$ are in $KH_1(k)$ and $KH_{-1}(k)$, respectively.
The theorem below generalizes to general twisted $EP$-tuples over the group $\Z$, the result proved in \cite{eptwist} for twisted Katsura tuples. 

\begin{thm}\label{thm:kepz}
Assume that $G=\Z$ in Theorem \ref{thm:kep} above. Then under the identification above, $\Phi^t$ identifies with multiplication by 
\[
\bar{D}^t=\begin{bmatrix} A^t& C^t(x)\\
0& B^t(x)
\end{bmatrix} 
\]
In particular there is a long exact sequence
\[
\xymatrix{KH_{n+1}(L)\ar[r]& KH_{n}(\ell)^{(\reg(E))}\oplus KH_n(\ell)^{(\reg(E))}\ar[d]^{I-\bar{D}^t}\\
          KH_n(L)& KH_{n}(\ell)^{(E^0)}\oplus KH_n(\ell)^{(E^0)}\ar[l]}
\]
If furthermore, both $k$ and $L$ are $K$-regular, then we may substitute $K$ for $KH$ in the sequence above. 
\end{thm}
\begin{proof}
Immediate from the calculations of the proof of Theorem \ref{thm:k01ep}. 
\end{proof}

\subsection{The Dennis trace}\label{subsec:dtrace}

Let $k/\ell$ be a flat ring extension, $n\ge 0$, $D_n:K_n(\cA_k(\cG))\to HH_n(\cA_k(\cG)/\ell)$ the Dennis trace and $\res:HH_n(\cA_k(\cG)/\ell)\to H_n(\cG,k/\ell)$ the restriction map. Put 
\[
\olD_n=\res\circ D_n:K_n(\cA_k(\cG))\to H_n(\cG,k/\ell).
\]
\begin{lem}\label{lem:conmutrace}
Let $T=(G,E,\phi_c)$ be a twisted EP-tuple. Assume as above that $E$ is row-finite and that $G$ acts trivially on $E$. Further assume that $T$ is pseudo-free and that $k[G]$ is regular supercoherent. Let $k/\ell$ be a flat ring extension. Then for $L=L(G,E,\phi_c)$ and $\cG=\cG(G,E,\phi)$ there is a commutative diagram with exact rows
\[
\xymatrix{K_{n+1}(L)\ar[d]^{\olD_{n+1}}\ar[r]& K_n(k[G])^{(\reg(E))}\ar[d]^{\olD_n}\ar[r]^{I-\Phi^t}&K_n(k[G])^{(E^0)}\ar[d]^{\olD_n}\ar[r]&K_n(L)\ar[d]^{\olD_n}\\
H_{n+1}(\cG,k/\ell)^{(\reg(E))}\ar[r]& H_n(G,k/\ell)^{(E^0)}\ar[r]^{I-\tau}&H_n(G,k/\ell)\ar[r]&H_n(\cG,k/\ell)}
\]
\end{lem}
\begin{proof}
By Corollary \ref{coro:kep}, the exact sequence at the top of the diagram is the excision sequence associated to the Cohn extension \ref{ext:cohnext}; by Corollary \ref{coro:hep} also the bottom sequence comes from the Cohn extension. Hence the diagram commutes by naturality of the Dennis trace. 
\end{proof}
Recall that if $u\in\cU(k)$, then 
\[
D_1(u)=d\log(u):=u^{-1}du\in \Omega^1_{k/\ell}=HH_1(k/\ell).
\]
In the next proposition we consider the $E^0\times\reg(E)$-matrix of homomorphisms $d\log(C)$ with $d\log(C)_{v,w}:G^{(\reg(E))}\to (\Omega^1_{k/\ell})^{(E^0)}$, $d\log(C)_{v,w}(g)=d\log(C_{v,w})$. We put
\begin{equation}
    \ul{M}^t=\begin{bmatrix}
       A^t&d\log C\\
       0&B^t
    \end{bmatrix}
\end{equation}
\begin{prop}\label{prop:Dcomm}
Let $T=(G,E,\phi_c)$, $k/\ell$, $\cG$ and $L$ be as in Theorem \ref{thm:k01ep}. Assume further that $T$ is pseudo-free. Then

\item[i)] $D_0:K_0(L)=\BF(E)\to \BF(E)\otimes_{\Z}k=H_0(\cG^{\bomega},k/\ell)\subset HH_0(L/\ell)$ is the scalar extension map. In particular $\olD_0$ induces an isomorphism $K_0(L)\otimes k\iso H_0(\cG^{\bomega},k/\ell)$.
\item[ii)] We have a commutative diagram with exact rows, where $K_1(k[G])=\cU(k)\oplus G_{\ab}$, $H_1(G,k/\ell)=\Omega^1_{k/\ell}\oplus 
G_{\ab}\otimes k$, and the maps labelled $\iota$ come from scalar extensions
\[
\xymatrix{K_1(k[G])^{(\reg(E))}\ar[d]^{d\log\oplus \iota}\ar[r]^{I-D^t}& K_1(k[G])^{(E^0)}\ar[d]^{d\log\oplus \iota}\ar[r]&K_1(L)\ar[d]^{ \olD_1}\ar[r]^{\partial}&\ker(I-A^t)\ar[d]^{\iota}\\
H_1(G,k/\ell)^{(\reg(E))}\ar[r]^{I-\ul{M}^t}&H_1(G,k/\ell)^{(E^0)}\ar[r]& H_1(\cG^{\bomega},k/\ell)\ar[r]^{\partial'}& \ker(k\otimes (I-A^t))}
\]
\end{prop}
\begin{proof}
Let $L'=L(E)$ be the Leavitt path algebra. By Theorem \ref{thm:k01ep}, the inclusion 
$\inc: L'\subset L$ induces an isomorphism at the $K_0$ level. In particular every element of 
$K_0(L)$ is a linear combination of classes of vertices. Assertion i) follows from the fact that $D_0$ maps a vertex to its class in $HH_0(L/\ell)$, which 
lies in $H_0(\cG^{\bomega},k/\ell)$, and the latter $k$-module equals $\BF(E)\otimes k$ by Corollary 
\ref{coro:h0hep}. If $R$ is an $\ell$ algebra, and $u\in\cU(R)$, then 
$D_1(u)\in HH_1(R/\ell)$ is the class of the cycle $d\log(u)=u^{-1}\otimes u$. It follows from this that the 
two leftmost vertical maps are induced by $\olD_1$. Hence in view of Theorem \ref{thm:k01ep} and Lemma 
\ref{lem:conmutrace} it only remains to show that the map $\tau$ of \eqref{map:tau} identifies with $\ul{M}^t$. It is clear that this is 
the case when restricted to the summand involving $\Omega^1_{k/\ell}$. It remains to prove that $\tau$ of \eqref{map:tau} and 
$\ul{M}^t$ agree on the other summands. First observe that if $a,x,y$
are commuting elements in an $\ell$-algebra $R$, with $x,y\in\cU(R)$, then 
\begin{equation}\label{eq:dlog}
 ad\log(xy)=ad\log(x)+ad\log(y)-b(a(xy)^{-1}\otimes y\otimes x)\in \HH(R/\ell)_1.   
\end{equation}
Next let $a\in k$, $g\in G$, and $(v,w)\in\reg(E)\times E^0$. Using \eqref{eq:dlog} at the third and fifth steps, we get
\begin{gather*}
    [\tau_{v,w}(a g)]=[\sum_{e\in vE^1w}ac(g,e)^{-1}\otimes c(g,e)\phi(g,e)]\\
    =[\sum_{e\in vE^1w} \res(ad\log(c(g,e)\phi(g,e)))]\\
    =[\sum_{e\in vE^1w}\res(ad\log(c(g,e))+ad\log(\phi(g,e)))]\\
    =[\sum_{e\in vE^1w}ad\log(c(g,e))+a\otimes\phi(g,e)]\\
    =[ad\log(C_{v,w}(g))]+[aB_{v,w}(g)].
\end{gather*}
\end{proof}

\begin{coro}\label{coro:matui}
In the setting of Proposition \ref{prop:Dcomm}, further assume that the twisting cocycle $c$ is trivial and that $k/\Z$ is flat. There is an exact sequence
\[
0\to \cU(k)\otimes\BF(E)\otimes k \to K_1(L)\otimes k\overset{\olD_1}{\lra}H_1(\cG,k)\to 0.
\]
\end{coro}

\begin{rem}\label{rem:compali}
In \cite{xlispectra}, Xin Li associated a permutative category $\fB\cG$ to any ample groupoid $\cG$ and showed that for any $\Z$-module $M$, $H_*(\cG,M)$ is the homology of the connective $K$-theory spectrum $\K(\fB\cG)$ with coefficients in $M$. There is an assembly map $\K(\fB\cG)\land \K(k)\to \K(\cA(\cG))$ and Li conjectures in \cite{xlinotes} that the latter is an equivalence whenever $k$ is regular Noetherian and $\cG$ is \emph{torsionfree}, that is, when $\cG_x^x$ is torsionfree for all $x\in\cG^{(0)}$. Reasoning as in \eqref{eq:k01neg}, we get that if $k$ is as in the corollary, $\cG$ is torsionfree, and the conjecture holds for $\cG$ and $k$, then there is an exact sequence
\[
H_0(\cG, \cU(k))\to K_1(\cA(\cG))\to H_1(\cG,\Z)\to 0. 
\]
\end{rem}

\section{Discretization}\label{sec:disc}
\numberwithin{equation}{section}
Let $\cS$ be a pointed inverse semigroup and $\cE=\cE(\cS)\subset \cS$ the subsemigroup of idempotent elements. Regard $
\{0,1\}$ as an idempotent semigroup under multiplication. A \emph{semicharacter} on $\cE$ is a nonzero homomorphism $\chi:\cE\to\{0,1\}$ of pointed semigroups. The set $\hat{\cE}$ of all semicharacters on $\cE$, equipped with with the topology of pointwise convergence is a compact Hausdorff space, and for each $p\in\cE$, the subset
\[
\hat{\cE}\supset D_p=\{\chi\colon \chi(p)=1\}
\]
is compact open. Preorder $\cE$ via $q\le p\iff qp=q$. Then the sets 
$$p_{\ge}=\{q\in\cE\colon q\le p\}$$
form a basis for the poset topology on $\cE$.
The semigroup $\cS$ acts on $\cE$ via $s\cdot p=sps^*$; this induces actions
on $\cE$  and $\hat{\cE}$ via 
\begin{gather*}
s\cdot-:s^*s_{\ge}\to ss^*_{\ge},\, s\cdot p=sps^*\\
s\cdot-:D_{s^*s}\to D_{ss^*},\, (s\cdot \chi)(p)=\chi(sps^*).
\end{gather*}
The \emph{universal groupoid} of $\cS$ is the transportation groupoid $\cG_u(\cS)=\cS\ltimes \hat{\cE}$; 
its \emph{discretization} is $\cG_d(\cS)=\cS\ltimes\cE^\times$, where $\cE^\times=\cE\setminus\{0\}$ is given the discrete topology. 
\begin{ex}\label{ex:epdisc}
Let $(G,E,\phi)$ be an Exel-Pardo tuple, $\cS=\cS(G,E,\phi)$ and $\cE=\cE(\cS)$. By \cite{ep}*{pages 1074--1075}, there is an $\cS$-equivariant homeomorphism $\hat{\cE}\cong\widehat{\fX(E)}$. Hence the universal groupoid $\cG_u(\cS)$ as defined in this section is the same as universal groupoid $\cG_u(G,E,\phi)=\cS\ltimes \widehat{\fX(E)}$ of Section \ref{subsec:ept}. Hence $\A(\cG_u(\cS))=C(G,E,\phi)$, by Lemma \ref{lem:cohn=leav}. The discrete space $\cE^\times$ is $\cS$-equivariantly isomorphic to the open subset $V=\widehat{\fX(E)}\setminus\{\theta\colon |\theta|=\infty\}$, and so $\cG_d(\cS)\cong\cG_u(\cS)_{|V}$. If $E$ is regular, then $V=U$, the open subset of the lemma, thus by the lemma, $\cG_d(\cS)\cong \cG_u(\cS)_{|U}$ and $\cA(\cG_d(\cS))\cong \cK(G,E,\phi)$ is the algebra defined also in Section \ref{subsec:ept}. For arbitrary $E$ an argument similar to that of part i) of the same lemma shows that $\A(\cG_d(\cS))=\hat{\cK}(G,E,\phi)$ is the algebra of \cite{eptwist}*{Section 6.3}, which, as explained there, is isomorphic to $\bigoplus_{v\in E^0}M_{\cP_v}k[G]$. 
\end{ex}

Let $\cS$ be any pointed inverse semigroup. Remark that every element of $\cG_d(\cS)$ can be written uniquely as $[s,s^*s]$. It follows that the characteristic functions $\chi_{[s,s^*s]}$ form a $k$-module basis of $\A(\cG_d(\cS))$. 
One checks that the $k$-linear map
\begin{equation}\label{map:rhod}
\rho_d:\A(\cG_d(\cS))\to M_{\cE^{\times}}\A(\cG_u(\cS)), \, \chi_{[s,s^*s]}\mapsto \epsilon_{ss^*,s^*s}\chi_{[s,D_{s^*s}]}
\end{equation}
is a homomorphism of algebras. Let $\cT$ be a category and $\cH:\aha\to\cT$ a functor. Assume that the restriction of $\cH$ to algebras with local units is  $M_{\cE^\times}$-stable. Then as mentioned above the isomorphism $\iota:\cH(\A(\cG_u(\cS)))\to \cH(M_{\cE^\times}\A(\cG_u(\cS)))$ resulting from applying $\cH$ to a corner inclusion $\phi\mapsto \epsilon_{p,p}\phi$ is independent of $p$. Hence we have a natural map
\begin{equation}\label{map:discinv}
\tilde{\rho}_d=\iota^{-1}\circ\cH(\rho_d)\colon \cH(\A(\cG_d(\cS)))\to\cH(\A(\cG_u(\cS))).
\end{equation}
We say that $\cH$ is \emph{discretization invariant} if \eqref{map:discinv} is an isomorphism for every inverse semigroup $\cS$.

\begin{rmk}\label{rem:xlidisc}
Xin Li showed \cite{xlinotes}*{Corollary 4.3} that a map related to \eqref{map:discinv} induces an isomorphism in groupoid homology. He also showed that the conjecture we discussed in Remark \ref{rem:compali} implies that if $\cG_u(\cS)$ is torsionfree and $k$ is regular Noetherian, then 
\begin{equation}\label{eq:kdisc}
 K_*(\A(\cG_{d}(\cS))\cong K_*(\A(\cG_u(\cS)).   
\end{equation} 
\end{rmk}

\begin{prop}\label{prop:nodisc} Hochschild homology is not discretization-invariant.    
\end{prop}
\begin{proof} 
Let $\cR_1$ be the graph consisting of one vertex and one loop, $\cS=\cS(\cR_1)$, and $C_1=C(\cR_1)$. Then by Example \ref{ex:epdisc}.
$\cA(\cG_u(\cS))=C_1$ and  $\cA(\cG_d(\cS))\cong\cK(\cR_1)\cong M_\infty(k)$. By matricial stability, $HH_*(M_\infty(k))=HH_*(k)$ is $k$ in degree $0$ and zero in positive degrees. On the other hand, using that, by \cite{lpabook}*{Theorem 1.5.18} (see also Lemma \ref{lem:cohn=leav}) $C(\cR_1)=L(\tilde{\cR}_1)$, and applying \cite{aratenso}*{Theorem 4.4} (or Theorem \ref{thm:hhep})
we obtain that $HH_n(C(\cR_1))=k^{(\N)}$ for $0\le n\le 1$ and vanishes for $n\ge 2$.
\end{proof}

Proposition \ref{prop:nodisc} implies that matricial stability and excision for algebras with local units do not suffice to guarantee discretization invariance, since $HH$ has both properties. 

\begin{prop}\label{prop:sidisc}
Let $(G,E,\phi_c)$ be an Exel-Pardo tuple, and let $\cG_u(G,E,\phi)$ and $\cG_d(G,E,\phi)$ be the universal groupoid and its discretization. Let $\cT$ be a triangulated category and $\cH:\ahal\to\cT$ an excisive, homotopy invariant, $E$-stable and $E^0$-additive functor. Then the map 
$\tilde{\rho}_d:\cH(\A(\cG_d(G,E,\phi)))\to \cH(\A(\cG_u(G,E,\phi)))$ of \eqref{map:discinv} is an isomorphism.
\end{prop}
\begin{proof}
Composing the isomorphism $\tilde{\cK}(G,E,\phi)\cong \A(\cG_d(G,E,\phi))$ of Example \ref{ex:epdisc} with that of \cite{eptwist}*{Section 6.3} we get an isomorphism 
\begin{gather*}
\left(\bigoplus_{v\in E^0}M_{\cP_v}\right)\rtimes G\iso \cA(\cG_d(G,E,\phi)),\\
\epsilon_{\alpha,\beta}\rtimes g\mapsto\chi_{[\alpha g (g^{-1}(\beta))^*]}.
\end{gather*}
By Lemma \ref{lem:cohn=leav}, we also have an isomorphism $C(G,E,\phi)\iso \A(\cG_u(G,E,\phi))$, $\alpha g\beta^*\mapsto\chi_{[\alpha g\beta^*,D_{\beta\beta^*}]}$. Moreover, we also have $M_{\cE^\times}\cong M_{\cP}$. One checks that under these isomorphisms, the map \eqref{map:rhod} becomes
\begin{gather*}
\rho'_d:\left(\bigoplus_{v\in E^0}M_{\cP_v}\right)\rtimes G\to M_{\cP}C(G,E,\phi),\\
\rho'_d(\epsilon_{\alpha,\beta}\rtimes g)=\epsilon_{\alpha,g^{-1(\beta)}}\alpha g (g^{-1}(\beta))^*.
\end{gather*}
Composing with the inclusion 
\begin{gather*}
\inc:k^{(E^0)}\rtimes G\to \left(\bigoplus_{v\in E^0}M_{\cP_v}\right)\rtimes G,\\
\inc(v\rtimes g)=\epsilon_{v,v}\rtimes g
\end{gather*}
we obtain the map
\begin{equation}\label{map:twistepsi}
v\rtimes g\mapsto \epsilon_{v,g^{-1}(v)}vg.
\end{equation}
Fix $v_0\in E^0$ and consider the matrix
\[
u\in \sum_{v\in E^0} \epsilon_{v_0,v}v.
\]
Then $u$ is an element of the multiplier algebra of $M_{\cP}(C(G,E,\phi))$ and satisfies $u^*u=1$. Thus it defines an inner endomorphism $\ad(u)$ of $M_{\cP}(C(G,E,\phi))$. One checks that $\ad(u)$ composed with \eqref{map:twistepsi} is the corner embedding $x\mapsto \epsilon_{v_0,v_0}x$. Since $\cH(\ad(u))$ is the identity map, we obtain that $\tilde{\rho}_d\circ \cH(\inc)$ coincides with the result of applying $\cH$ to the map $\phi:k^{(E^0)}\rtimes G\to C(G,E,\phi)$, $\phi(v\rtimes g)=vg$. Since by \cite{eptwist}*{Proposition 6.2.3 and Theorem 6.3.1} both $\cH(\inc)$ and $\cH(\phi)$ are isomorphisms, we conclude that $\tilde{\rho}_d$ is an isomorphism.
\end{proof}

\begin{conj}\label{conjdisc} Let $\cT$ be a triangulated category and $\cH:\ahal\to\cT$ an excisive, homotopy invariant, matricially-stable and infinitely additive functor. Then $\cH$ is discretization invariant.
\end{conj}

\begin{rem}\label{rem:semi}
The idempotent semigroup $\cE=\cE(\cS)$, with the preorder defined above is a semilattice, where the meet is the semigroup product. One may also consider actions of inverse semigroups on more general posets. In fact Xin Li proves that his conjecture implies the isomorphism \eqref{eq:kdisc} for germ groupoids of semigroup actions on general locally finite weak semilattices. The map \eqref{map:rhod} also makes sense in this more general context. Hence one could define a more stringent version of discretization invariance by requiring it holds for actions on locally finite weak semilattices. This in turn leads to a stronger version of the conjecture above. 
\end{rem}

\appendix

\section{Corner skew Laurent polynomial algebras}\label{sec:appa}
\numberwithin{equation}{section}
Let $R$ be a unital algebra and $\psi:R\to R$ a corner isomorphism. Let $S=R[t_+,t_-;\psi]$ be the corner skew Laurent polynomial ring of \cite{fracskewmon}. 

Remark that the $\Z$-grading on $S$ induces one on $\HH(R,S)$ and $\HH(S)$, which together with the chain complex grading, make them into bigraded $k$-modules. 

\begin{lem}\label{lem:kappa}
There is a natural homomorphism of bigraded  $k$-modules
\goodbreak
\noindent $\kappa:\HH(S)[-1]\to \HH(S)$ such that $1-\psi=b\kappa+\kappa b$. 
\end{lem}
\begin{proof} Let $C^{bar}(S)=(S^{\otimes \bullet +2},b')$ be the bar resolution and $s:C_n^{bar}(S)\to C_{n+1}^{bar}(S)$, $s(x)=1\otimes x$ . Let 
$\tilde{\psi}:C^{bar}(S)\to C^{bar}(S)$, 
\begin{equation}\label{map:barpsi}
\tilde{\psi}(a_0\otimes\cdots\otimes a_{n+1})=a_0t_-\otimes\psi(a_1)\otimes\cdots\otimes \psi(a_n)\otimes t_+a_{n+1}.    
\end{equation}

Let $\tilde{\kappa}_0:S\otimes S\to S\otimes S\otimes S$, be the $S$-bimodule homomorphism determined by $\tilde{\kappa}_0(1\otimes 1)=-1\otimes t_-\otimes t_+$. Define inductively 
\begin{equation}\label{map:indukappa}
\tilde{\kappa}_{n+1}(1\otimes x\otimes 1)=s(1\otimes x\otimes 1-t_-\otimes\psi(x)\otimes t_+-\tilde{\kappa}_n(b'(1\otimes x\otimes 1))). 
\end{equation}
Then $1-\tilde{\psi}=b'\tilde{\kappa}+\tilde{\kappa}b'$.     
It follows that $\kappa=\tilde{\kappa}\otimes_{S^e}S$ has the required properties. 
\end{proof}
\begin{coro}\label{coro:kappa}
Let $\iota:\HH(R,S)\to \HH(S)$ be the inclusion map. Then
$\theta:\cone(1-\psi:\HH(R,S)\to \HH(R,S))\to \HH(S)$, defined on 
$\cone(1-\psi:\HH(R,S)\to \HH(R,S))_n=\HH(R,S)_n\oplus\HH(R,S)_{n-1}$ as $\theta(x,y)=\iota(x)+\kappa(y)$ is a graded homomorphism of chain complexes.
\end{coro}
\begin{rem}\label{rem:elkappa}
It follows from the inductive formula \eqref{map:indukappa} that the map $\tilde{\kappa}$ preserves the degenerate subcomplex, and so descends to a homotopy $\tilde{\kappa}_{\norm}$ between $\psi$ and the identity of the normalized complex $C^{bar}(S)_{norm}$. A straightfoward induction argument shows that
\begin{gather*}
 \tilde{\kappa}_{\norm}(1\otimes x_1\otimes\dots\otimes x_n\otimes 1)=\\
\sum_{i=0}^n(-1)^{i+1}1\otimes x_1\otimes\dots\otimes x_i\otimes t_-\otimes \psi(x_{i+1})\otimes\dots\otimes \psi(x_n)\otimes t_+   
\end{gather*}
Hence the map
\begin{gather*}
\kappa_{\norm}:HH(S)_{\norm}[-1]\to HH(S)_{\norm},\\
\kappa_{\norm}(x_0\otimes\dots\otimes x_n)=\sum_{i=0}^n(-1)^{i+1}t_+x_0\otimes x_1\otimes\dots\otimes x_i\otimes t_-\otimes \psi(x_{i+1})\otimes\cdots\otimes \psi(x_n)
\end{gather*}
satisfies $b\kappa_{\norm}+\kappa_{\norm}b=1-\psi$.
\end{rem}
\begin{lem}\label{lem:locacono}
Let $\xi:M\to M$ be a chain complex endomorphism. Let $M[\xi^{-1}]$ be the colimit of the $\N$-directed system
\[
\xymatrix{M\ar[r]^\xi&M\ar[r]^\xi&\dots&M\ar[r]^\xi&\dots}
\]
Let $\xi':M[\xi^{-1}]\to M[\xi^{-1}]$ be map induced by $\xi$. Then the natural map
\[
\cone(1-\xi:M\to M)\to \cone(1-\xi': M[\xi^{-1}]\to M[\xi^{-1}])
\]
is a quasi-isomorphism.
\end{lem}
\begin{proof}
We may regard $M$ as a chain complex of $\Z[x]$-modules with $x$ acting via $\xi$, and $M[\xi^{-1}]=M\otimes_{\Z[x]}\Z[x,x^{-1}]$. 
The natural map of the lemma induces a map of triangles in the derived category of chain complexes
\[
\xymatrix{\ker(1-\xi)[-1]\ar[d]\ar[r]&\cone(1-\xi)\ar[d]\ar[r]&\coker(1-\xi)\ar[d]\\
\ker(1-\xi')[-1]\ar[r]& \cone(1-\xi')\ar[r]& \coker(1-\xi')}
\]
Because the vertical maps at both ends are isomorphisms of chain complexes, that in the middle is a quasi-isomorphism. 
\end{proof}
The following proposition provides an explicit quasi-isomorphism whose existence was established in \cite{aratenso}*{Proposition 3.4}.
\begin{prop}\label{prop:skewcorner}
Let $R$ be a unital $k$-algebra, $\psi:R\to R$ a corner isomorphism,  $S=R[t_+,t_-;\psi]$ the corner skew Laurent polynomial ring. Equip $S$ with its natural $\Z$-grading and $\HH(R,S)$ and $\HH(S,S)$ with the induced gradings. Then the bigraded chain homomorphism $\theta:\cone(1-\psi:\HH(R,S)\to \HH(R,S))\to \HH(S)$ of Corollary \ref{coro:kappa} is a quasi-isomorphism. 
\end{prop}
\begin{proof} Taking an appropriate colimit as in Lemma \ref{lem:locacono}, we obtain algebras $R'=R[\psi^{-1}]$ and $S'=S[\psi^{-1}]$, such that the endomorphism $\psi':R'\to R'$ induced by $\psi$ is an automorphism and $S'=R'[t_+,t_-,\psi']=R'[t,t^{-1};\psi']=R'\rtimes_\psi'\Z$ is the crossed product. Because the Hochschild complex commutes with filtering colimits, it follows from Lemma \ref{lem:locacono} that the map
$\cone(1-\psi:\HH(R,S)\to\HH(R,S))\to \cone(1-\psi':\HH(R',S')\to \HH(R',S'))$ is a quasi-isomorphism. Similarly, using Lemma \ref{lem:kappa} and again that $\HH$ commutes with filtering colimits, we get that  $\HH(S)\to \HH(S')$ is a quasi-isomorphism.
Because by construction $\theta$ comes from a map $\bar{\theta}:\cone(1-\tilde{\psi}:\cC^{bar}(R,S)\to \cC^{bar}(R,S))\to \cC^{bar}(S)$ and because the bar complexes also commute with filtering colimits, we get that $\theta'$ also comes from a map 
$\cone(1-\tilde{\psi'}:\cC^{bar}(R',S')\to \cC^{bar}(R',S'))\to \cC^{bar}(S')$. Using the fact that because $\psi'$ is an automorphism, $S'\cong k[t,t^{-1}]\otimes R'$ as right $R'$-modules, we obtain 
\begin{align*}
H_0(\cone(1-\tilde{\psi'}))=&S'\otimes_{R'} S'/\langle a_0\otimes a_1-a_0t^{-1}\otimes ta_1\rangle\\
=& k[t,t^{-1}]\otimes S'/\langle a_0\otimes a_1-a_0t^{-1}\otimes ta_1\rangle\\
=& S'
\end{align*}
Thus $\cone(1-\tilde{\psi'})$ is an $S'$-bimodule resolution of $S'$ and $\bar{\theta'}$ lifts the identity of $S'$. It follows that $\theta'$ is a homotopy equivalence. This finishes the proof. 
\end{proof}

\end{document}